\documentclass[reqno]{amsart}%
\usepackage{amsfonts}
\usepackage{amssymb,amsthm,amsmath}
\usepackage{bbm}
\usepackage[numbers,sort&compress]{natbib}
\usepackage{color}
\usepackage{graphicx}
\usepackage{tikz}
\usepackage{verbatim}
\usepackage{makecell}
\usepackage{bm}
\usepackage{diagbox}
\usepackage{mathrsfs}
\usepackage[dvips]{epsfig}
\usepackage[colorlinks=true,linkcolor=blue,citecolor=brown]{hyperref}

\def\be{\begin{equation}}
	\def\ee{\end{equation}}
\def\ba{\begin{array}}
	\def\ea{\end{array}}

\title[The Existence of Full-Dimensional KAM tori for NLKG]{Non-relativistic limit of the nonlinear Klein-Gordon equation: Uniform in time convergence of KAM almost-periodic solutions \footnote{Li is supported by the National Natural Science Foundation of China (Grant No. 125B2011).}}

\author[H. Cong]{Hongzi Cong}
\address[H. Cong]{School of Mathematical Sciences, Dalian University of Technology, Dalian 116024, China}
\email{conghongzi@dlut.edu.cn}

\author[S. Li]{Siming Li}
\address[S. Li]{School of Mathematical Sciences, Dalian University of Technology, Dalian, Liaoning 116024, China}
\email{lisiming@mail.dlut.edu.cn}

\author[X. Wu]{Xiaoqing Wu}
\address[X. Wu]{School of Mathematics and Statistics, Ningxia University, Yinchuan 750021, China}
\email{xqwu@nxu.edu.cn}

\theoremstyle{plain}
\newtheorem{theorem}{Theorem}[section]

\newtheorem{lemma}[theorem]{Lemma}

\newtheorem{claim}[theorem]{Claim}

\theoremstyle{definition}
\newtheorem{definition}[theorem]{Definition}

\newtheorem{remark}[theorem]{Remark}

\begin{document}




\begin{abstract}
 This paper investigates the long-time dynamical stability of the nonlinear Klein-Gordon equation under periodic boundary conditions in the non-relativistic limit. Building upon the finite-dimensional KAM framework recently established in \cite{bambusi2025NLKG}, we extend this approach to full-dimensional KAM dynamics. We prove the uniform existence of full-dimensional KAM tori across the entire limit path $c \in [1, \infty)$. Furthermore, by employing external parameter frequency modulation and a frequency-domain splitting strategy, we establish global-in-time uniform convergence of these full-dimensional tori to the corresponding tori of the nonlinear Schrödinger  equation as $c \to \infty$.  
\end{abstract}

\maketitle

\section{Introduction and main results}
\subsection{Introduction}

In the field of mathematical physics, the nonlinear Klein-Gordon (NLKG) equation, as a classical Hamiltonian partial differential equation (PDE), has long occupied a central position in the study of the qualitative properties of its solutions.  In various physical contexts, the speed of light parameter $c$ has a significant influence on the underlying dynamical evolution of the system, as explored in coupled models such as the Maxwell-Klein-Gordon equation under the non-relativistic limit \cite{MR2033563,MR1949296}. 
 Specifically, we investigate the one-dimensional nonlinear Klein-Gordon equation 
 \begin{equation}\label{NLKG}
 		\frac{1}{c^{2}} u_{t t}-u_{x x}+c^{2} u+V*u+ \epsilon u^3 =0,
 	\quad\mathbb{T} := \mathbb{R}/2\pi\mathbb{Z},
 \end{equation}
 where $c \in [1, +\infty)$ is the speed of light.
 From an intuitive physical perspective, as the speed of light tends to infinity ($c \to \infty$), relativistic effects gradually recede, and the system is expected to transition into the non-relativistic regime. Mathematically, by introducing an appropriate gauge transformation to modulate the phase of the highly oscillatory components, the NLKG equation can be formally reduced to the classical cubic nonlinear Schrödinger (NLS) equation
 \begin{equation}\label{NLS}
	\mathrm{i}\varphi_{t} - \frac{1}{2}\varphi_{xx} 
+ \frac{1}{2} V*\varphi +  
\frac{3}{4}\epsilon |\varphi|^2\varphi=0,
\quad x \in \mathbb{T}. 
\end{equation}
Such a rigorous mathematical transition from the relativistic framework to the non-relativistic one constitutes the well-known \emph{non-relativistic limit} problem.

The primary objective of the present paper is to investigate the global-in-time stability and limiting convergence of full-dimensional KAM tori (corresponding to almost-periodic solutions with infinitely many frequencies) transitioning from the relativistic equation to the non-relativistic one. By developing a full-dimensional KAM framework tailored for the singular limit, we address this problem under periodic boundary conditions. Specifically, we establish the uniform existence of full-dimensional tori across the entire limit path $c \in [1, \infty)$, as well as their global-in-time uniform convergence to the limiting NLS tori. To clarify the non-triviality of overcoming this regime under the non-relativistic limit, we briefly review the historical development of KAM theories in the dynamics of Hamiltonian PDEs.

%

Classical KAM theory has been extensively applied to Hamiltonian PDEs to construct quasi-periodic solutions (finite-dimensional KAM tori) and to ensure their long-time stability. The existence of quasi-periodic solutions for the NLS and nonlinear wave (NLW) equations was established in the foundational works \cite{K1987,MR1290785,W1990,P'1996}. 
Subsequently, the theory was extended to more complex settings, such as higher spatial dimensions and completely resonant regimes \cite{MR2019250,Yuan2006, MR3502160}.
A natural progression in this field has been the extension of the KAM tori to infinite dimensions, namely, searching for almost-periodic solutions.
  A core open problem is
  the existence of ``maximal tori'' with Sobolev regularity, which was posed by Kuksin in 2004 
 (see Problem 7.1 in \cite{Kuksin2004}): ``{\it Can full-dimensional KAM tori be expected for  Hamiltonian PDEs with a suitable decay, such as 
	$$
	I_n \sim|n|^{-C}\quad \text{as}\quad |n|\to \infty,
	$$ 
	where \(C>0\)?}'' This question was also highlighted in Procesi's 2022 ICM report (see {\it Q10} in \cite{MR4680374}). 
	Bourgain \cite{Bourgain2005JFA} carried out pioneering work in this direction,  obtaining tori with radius decaying as $\mathcal{O}(e^{-\sqrt{|n|}})$. Subsequently, the study of full-dimensional KAM tori has been further developed in several aspects, such as extending to Sobolev weak solutions \cite{BMP2021Poincare, BMP2023} and relaxing the spatial topological decay rates \cite{CLSY2018JDE,CY2021}. Furthermore, utilizing Pöschel's framework \cite{Poschel2002}, almost-periodic solutions were constructed for the parameter-free NLS equation \cite{BGR2024}. Recently, Cong \cite{Cong2023TheEO} through refined analytical estimates pushed the survival boundary of maximal tori to a slower decay  $(e^{-\ln^\sigma |n|})$ regime. These studies show that when the spatial topological decay rate becomes slower, the control of high-frequency resonances imposes higher requirements on analytical estimates.

Under the non-relativistic limit, traditional PDE approaches are constrained by either finite time horizons or spatial geometries. Classical energy and modulation methods guarantee convergence only over finite time intervals $[0, T]$ \cite{Tsu1984, Naj1992, MNO2002}. While scattering theory and recent advances in harmonic analysis have successfully established global-in-time uniform convergence on the whole space $\mathbb{R}^d$ \cite{arX2308}, these methods inherently rely on spatial dispersion. Consequently, they remain inapplicable on compact domains where energy is trapped. On bounded domains, cumulative phase drifts from highly oscillatory components cause standard approximations to fail as $t \to \infty$. Although Hamiltonian Normal Forms extend these approximations to almost global times \cite{MR3745401}, widespread resonances prevent convergence over the infinite time horizon.

%

To bridge this infinite time horizon on compact domains, Bambusi, Belloni, and Giuliani \cite{bambusi2025NLKG} recently applied KAM theory to establish the global-in-time convergence of finite-dimensional KAM tori for the NLKG equation under Dirichlet boundary conditions. Because the tori are finite-dimensional, they introduced a frequency alignment mechanism: by adjusting a finite number of internal actions (amplitudes), one can force the two equations to share an identical frequency vector. For sufficiently large $c$, aligning these finite frequencies controls the resonances induced by the singular parameter.

In \cite{bambusi2025NLKG}, the authors state that ``{\it it would be very interesting to extend the results of this paper to almost-periodic solutions}''. However, transitioning from the finite-dimensional quasi-periodic regime to the full-dimensional regime encounters substantial analytical obstacles. As the dimension tends to infinity ($n = \infty$), the high-frequency modes escape the dominance of $c$, generating an irreversible structural mismatch with the limiting NLS equation. Furthermore, unlike the Dirichlet boundary conditions utilized in \cite{bambusi2025NLKG}, the periodic boundary conditions investigated in the present paper impose momentum conservation ($\sum \pm n_i = 0$), which forces the normal frequencies to form a completely coupled resonant lattice network in Fourier space \cite{Yuan2006}.

Under this full-dimensional resonant regime, frequency alignment via internal parameter modulation becomes unviable. The rapid, sub-exponential amplitude decay required to preserve the spatial regularity critically shrinks the tuning space at high frequencies, precluding tractable measure estimates without severe regularity loss.
To circumvent these difficulties, classical full-dimensional KAM frameworks (e.g., Bourgain \cite{Bourgain2005JFA} and Pöschel \cite{Poschel2002}) introduce external parameters for non-resonant frequency regulation. Following this paradigm, the present paper adopts external parameters to perform the KAM iteration, avoiding the regularity loss associated with infinite-dimensional internal modulations.

However, in the singular limit $c \to \infty$, external parameters alone are insufficient due to the $c$-dependent nonlinear deformation of the normal frequencies. To ensure the uniform existence of full-dimensional tori, we establish a non-resonance condition uniformly independent of $c$. This allows us to maintain a large measure for the Cantor parameter set when solving the homological equations under the slow decay topology, thereby establishing uniform existence along the limit path (Theorem \ref{th1}).

Regarding the limiting convergence between the full-dimensional tori of the two systems, the polynomial growth of high-frequency mismatches precludes full alignment via external parameters. We therefore adopt a frequency-domain splitting strategy: low-frequency modes are aligned via frequency matching to eliminate phase drifts, while high-frequency drifts are suppressed in a lower-order topology. Consequently, we obtain a $c^{-2}$ convergence rate for the full-dimensional tori (Theorem \ref{th2}). Technical expositions are deferred to Section \ref{comments}.

%

\subsection{Main results}

We state our main results using the notations and terminologies established in \cite{Cong2023TheEO}. For $p \ge 4$, $2 < \sigma \le 3$, and $r > 1$, we define
 the Banach space
\begin{equation}\label{042501}
\mathfrak{B}_{\sigma,r,p}=\left\lbrace z=(z_{n})_{{n}\in\mathbb{Z}} \in \mathbb{C}^{\mathbb{Z}}:	\|z\|_{\sigma,r,p}=\sup_{{n}\in\mathbb{Z}}
	\left|z_{n}\right| \lfloor n \rfloor ^p
	e^{r\ln^{\sigma}\lfloor n \rfloor}<\infty\right\rbrace ,
\end{equation}
where 
\begin{equation*}\label{030401}
	\lfloor n \rfloor :=\max\{2^{10},|n|\}.
\end{equation*}
To simplify the notation, we fix $\sigma$ and $r$, and consistently denote $\|z\|_{\sigma, r, p}$ and $\mathfrak{B}_{\sigma, r, p}$ as $\|z\|_p$ and $\mathfrak{B}_p$, respectively. We introduce the action variables $I_n = |z_n|^2$. The initial data $I_n(0) = |z_n(0)|^2$ are assumed to satisfy
\begin{equation*}
	I_{n}(0)\sim  |n|^{-2p} e^{-2r \ln^{\sigma}\lfloor n \rfloor}\qquad
	\text{as}\quad |n|\to+\infty.
\end{equation*}


The infinite-dimensional Hamiltonians for equations \eqref{NLKG} and \eqref{NLS} are given by (see Section \ref{sec 2} for details):
\begin{align*}
	H(z,\bar{z})=N(z,\bar{z})+R(z,\bar{z})
\end{align*}
and
\begin{align*}
	H^{\text{NLS}}(z,\bar{z})=N^{\text{NLS}}(z,\bar{z})+R^{\text{NLS}}(z,\bar{z}),
\end{align*}
where $N$ and $N^{\text{NLS}}$ represent the integrable components, while $R$ and $R^{\text{NLS}}$ denote the corresponding small perturbations.
The integrable part are given by
\begin{align}\label{integrable}
	N(z,\bar{z})
	&=\sum_{{n}\in\mathbb{Z}}
	\lambda_{n} (V) |z_{n}|^2 
	\quad\text{with}\quad
	\lambda_{n} (V)= c \sqrt{c^2+n^2+V_n},\\
		N^{\text{NLS}}(z,\bar{z})
	&=\sum_{{n}\in\mathbb{Z}}
	\lambda_{n}^{\text{NLS}} (V) |z_{n}|^2 
	\quad\text{with}\quad
	\lambda_{n}^{\text{NLS}} (V)=\frac{1}{2}\left( n^2+V_n\right),
	\label{integrable_NLS}
\end{align}
where $ \lambda = \left(\lambda_{n} \right)_{{n}\in\mathbb{Z}}$ and $\lambda^{\text{NLS}} = \left(\lambda_{n}^{\text{NLS}} \right)_{{n}\in\mathbb{Z}}$
 are frequencies
and parameters $V=\left( V_n\right)_{n\in\mathbb{Z}}$
belong to the space 
\begin{equation}\label{potential}
	\mathcal{V}:=\left\lbrace V= \left( V_{n}  \right)_{{n}
		\in\mathbb{Z}} : {V}_{n} \in
	\left[ 0,1\right],{n}
	\in\mathbb{Z}  \right\rbrace,
\end{equation}
which we endow the product probability measure.

Noting that for $V_n = 0$, the normal frequency is $\lambda_n(V) = c \sqrt{c^2 + n^2}$. As $V_n$ varies within $[0, 1]$, $\lambda_n(V)$ shifts monotonically from $c \sqrt{c^2 + n^2}$ to $c \sqrt{c^2 + n^2 + 1}$. Accordingly, for any $c \in [1, \infty)$, we define the frequency space $\Pi$ as
\begin{equation}\label{Pi_c}
	\Pi\equiv \Pi(c):=\left\lbrace \omega=(\omega_n)_{n\in\mathbb{Z}}\in \mathbb{R}^{\mathbb{Z}}
	:\omega_n -c\sqrt{c^2+n^2} \in \left[ 0, \frac{c}{3\sqrt{c^2+n^2}}\right] ,\; n\in\mathbb{Z}
	\right\rbrace 
\end{equation}
with product probability measure.
Note that this construction ensures a well-defined correspondence between each frequency $\omega \in \Pi$ and a potential parameter $V \in \mathcal{V}$.
Analogously, we define the frequency space $\Pi^{\text{NLS}}$ as
\begin{equation}\label{Pi}
	\Pi^{\text{NLS}}:=\left\lbrace \omega^{\text{NLS}}=(\omega^{\text{NLS}}_n)_{n\in\mathbb{Z}}\in \mathbb{R}^{\mathbb{Z}}
	:\omega_n^{\text{NLS}} -\frac{1}{2}n^2 \in \left[ 0, \frac{1}{3} \right] ,\; n\in\mathbb{Z}
	\right\rbrace 
\end{equation}

To establish the existence of almost-periodic solutions for \eqref{NLKG} and \eqref{NLS}, we impose the following Diophantine non-resonance conditions:
\begin{definition}\label{Diophantine}
 For a fixed $c \in [1, \infty)$, a frequency vector $\omega = (\omega_n)_{n \in \mathbb{Z}} \in \Pi$ is said to be non-resonant if there exists a constant $\gamma > 0$ such that for any $\ell \in \mathbb{Z}^{\mathbb{Z}}$ with $3 \le |\ell| < \infty$, the following inequalities hold:\\
	\textbf{(1)} For all such $\ell$,
	\begin{equation}\label{NR_KG1}
		\left| 	\sum_{{n}\in\mathbb{Z}} \ell_{n} \cdot \omega_{n} \right|
		\geq   {\gamma} \left( 
		\prod_{|n| \leq  n_3^{*}(\ell),\, \ell_n\neq 0}  
		\frac{1}{ |\ell_{n}|^5\,\lfloor n\rfloor^6}\right)^{5};
	\end{equation}
    \textbf{(2)} Whenever $\alpha(\ell) := \sum_{n \in \mathbb{Z}} \ell_n \neq 0$, additionally
		\begin{equation}\label{NR_KG2}
		\left| 	\sum_{{n}\in\mathbb{Z}} \ell_{n} \cdot \omega_{n} \right|
		\geq   {\gamma c^2} \left( 
		\prod_{|n| \leq  n_3^{*}(\ell),\, \ell_n\neq 0}  
		\frac{1}{ |\ell_{n}|^5\,\lfloor n\rfloor^6}\right)^{6},
	\end{equation}
where 
$(n_i^*(\ell))_{i\geq1}$ is the decreasing rearrangement of
\begin{equation*}
	\left\{|n|\ \mbox{where $n$ is repeated }
	\ell_{n}   \mbox{ times}\right\}.
\end{equation*}
\end{definition}

\begin{definition}\label{Diophantine2}
A frequency vector $\omega^{\text{NLS}} = (\omega_n^{\text{NLS}})_{n \in \mathbb{Z}} \in \Pi^{\text{NLS}}$ is said to be non-resonant if there exists $\gamma > 0$ such that for any $\ell \in \mathbb{Z}^{\mathbb{Z}}$ with $3 \le |\ell| < \infty$, the following inequality holds:
	\begin{equation}\label{NR_NLS}
	\left| 	\sum_{{n}\in\mathbb{Z}} \ell_{n} \cdot \omega_{n}^{\text{NLS}} \right|
	\geq   {\gamma}\left( 
		\prod_{|n| \leq  n_3^{*}(\ell),\, \ell_n\neq 0}  
		\frac{1}{ |\ell_{n}|^5\,\lfloor n\rfloor^6}\right)^{5}.
\end{equation}
\end{definition}

\begin{remark}\label{l=1}
\textbf{(1)} Uniformity with respect to $c$: Note that while the frequencies $\omega_n$ depend on $c$, the non-resonance conditions \eqref{NR_KG1} and \eqref{NR_KG2} are defined by a lower bound independent of $c$. This uniformity is crucial, as it ensures that the existence result for the full-dimensional tori holds uniformly for all $c \in [1, \infty)$.\\
\textbf{(2)} Range of applicability: The conditions \eqref{NR_KG1}, \eqref{NR_KG2}, and \eqref{NR_NLS} are stipulated for $3 \le |\ell| < \infty$ (see Section \ref{sec 5}). For $0 < |\ell| \le 2$, the momentum conservation condition \eqref{mome} implies
\begin{equation*}\left| \sum_{n \in \mathbb{Z}} \ell_n \omega_n \right| \ge 2c^2 \quad \text{and} \quad \left| \sum_{n \in \mathbb{Z}} \ell_n \omega_n^{\text{NLS}} \right| \ge 1.
\end{equation*}\\
\textbf{(3)} Technical distinction from \cite{Bourgain2005JFA}: Unlike the conditions in \cite{Bourgain2005JFA}, our estimates \eqref{NR_KG1}, \eqref{NR_KG2}, and \eqref{NR_NLS} involve only the first three largest indices $n_i^*(\ell)$ (i.e., $|n| \le n_3^*(\ell)$). This restriction is critical as it simplifies the treatment of the homological equations and allows for closing the KAM iteration under the slow decay topology (see Lemma \ref{IL}).

\end{remark}
%

We are now in a position to state our main result concerning the existence of almost-periodic solutions for equations \eqref{NLKG} and \eqref{NLS}.
\begin{theorem}\label{th1}
	Given any $2<\sigma\le 3$, $r>1$, $\gamma'>0$, and $c\in [1,\infty)$. 	There exist two subsets $\mathcal{R} \subset \Pi$ and $\mathcal{R}^{\text{NLS}} \subset \Pi^{\text{NLS}}$ with Lebesgue measure
	$$\mbox{meas}\  {\mathcal{R}}, \mbox{meas}\  {\mathcal{R}}^{\text{NLS}}= O( \gamma'),$$
	such that for any  $\omega\in\Pi\setminus  {\mathcal{R}}$ and $\omega^{\text{NLS}}\in\Pi^{\text{NLS}}\setminus {\mathcal{R}}^{\text{NLS}}$,
	there exists a small $\epsilon_*:=\epsilon_*(\sigma,r ,\gamma')>0$ depending on $\sigma,r,\gamma'$, and for any $0<\epsilon<\epsilon_*$,
	there exist two vectors $V,V^{\text{NLS}}\in\mathcal{V}$  such that
	the equations (\ref{NLKG}) and \eqref{NLS}
	possess invariant tori $\mathcal{E}$ and $\mathcal{E}^{\text{NLS}}$, respectively,  satisfying the following conditions:\\
	(1) the amplitude $I=(I_{n})_{{n}\in\mathbb{Z}}$ of $\mathcal{E}$ and
	$I^{\text{NLS}}=(I^{\text{NLS}}_{n})_{{n}\in\mathbb{Z}}$
	 of  $\mathcal{E}^{\text{NLS}}$ are restricted as
	$$\frac{1}{4}\lfloor n \rfloor ^{-2p}e^{-2r\ln^{\sigma}\lfloor n \rfloor} 
	 \le  I_{n} , \, I^{\text{NLS}}_{n}   \le
	   4\lfloor n \rfloor ^{-2p} e^{-2r\ln^{\sigma}\lfloor n \rfloor};$$
	(2) the frequencies on $\mathcal{E}$ and $\mathcal{E}^{\text{NLS}}$ are respectively prescribed as
	$(\omega_{n})_{n\in\mathbb{Z}}$ and $
	\left(\omega^{\text{NLS}}_{n}\right) _{n\in\mathbb{Z}};$\\
	(3) the invariant torus $\mathcal{E}$ and $\mathcal{E}^{\text{NLS}}$ are linearly stable.
\end{theorem}



To characterize the limiting convergence of solutions as $c \to \infty$, we introduce the complex coordinate
\begin{equation}\label{coor_tf}
	\psi = \frac{1}{\sqrt{2}}
	\left(D_c^{\frac{1}{2}}u -\mathrm{i}D_c^{-\frac{1}{2}} v \right)
	\quad\text{with}\quad
	D_c=c^{-1}\sqrt{-\partial_{xx}+c^2+V*},
\end{equation}
under which the NLKG equation \eqref{NLKG} is rewritten as
\begin{equation}\label{NLKG2}
	i \psi_t + c^2D_c\psi + \frac{1}{4}  D_c^{-\frac{1}{2}} \left[   D_c^{-\frac{1}{2}} (\psi + \bar{\psi}) \right]^3=0.
\end{equation}
To implement the frequency-domain splitting strategy, we define the translated frequency space $\Pi' := \Pi - c^2$ and partition it into low-frequency and high-frequency components
\begin{equation}\label{260513-1}
	\widetilde{\Pi}_{ \mathcal{N}}' \equiv \widetilde{\Pi}_{ \mathcal{N}}'(c) := \left\{  (\omega_n)_{|n| \le \mathcal{N}} : \omega_n - c\sqrt{c^2 + n^2}+c^2 \in \left[ 0, \frac{c}{3\sqrt{c^2 + n^2}} \right], \; n \in \mathbb{Z} \right\}
\end{equation}
and
\begin{equation*}
\widehat{\Pi}_{\mathcal{N}}'	\equiv \widehat{\Pi}_{ \mathcal{N}}'(c) := \left\{  (\omega_n)_{|n| > \mathcal{N}} : \omega_n - c\sqrt{c^2 + n^2}-c^2 \in \left[ 0, \frac{c}{3\sqrt{c^2 + n^2}} \right], \; n \in \mathbb{Z} \right\},
\end{equation*}
where $\mathcal{N} = \mathcal{N}(c)$ is a monotonically increasing integer.
 Similarly, we partition the NLS frequency space $\Pi^{\text{NLS}}$ into $\widetilde{\Pi}_{\mathcal{N}}^{\text{NLS}}$ and $\widehat{\Pi}_{\mathcal{N}}^{\text{NLS}}$, and express any frequency vector $\omega$ as $(\widetilde{\omega}, \widehat{\omega})$ according to the threshold $\mathcal{N}$. 
 Furthermore, to characterize the initial data $I(0)$ in our main result, we define the action space as
 $$\mathcal{I}_0 = \left\{ I = (I_n)_{n\in\mathbb{Z}} : \frac{9}{20} \le I_n  \lfloor n \rfloor^{2p} e^{2r \ln^\sigma \lfloor n \rfloor}  \le \frac{13}{20}, \;   n \in \mathbb{Z} \right\}.$$
 With these notations established, we now state the main convergence theorem:
\begin{theorem}\label{th2}

	For any $2<\sigma\le 3$, $p>4$, $r>1$, and $\gamma'>0$, there exists a small constant $\epsilon_* := \epsilon_*(\sigma,p,r ,\gamma') > 0$ such that for any $0<\epsilon<\epsilon_*$ and $c\ge \epsilon^{-0.5}$, there exist an integer $\mathcal{N} = \mathcal{N}(c)$, two frequency sets $\Xi' \equiv \Xi'(c) \subset \Pi'$ and $\Xi^{\text{NLS}} \equiv \Xi^{\text{NLS}}(c) \subset \Pi^{\text{NLS}}$ satisfying 
	$\widetilde{\Xi}'_{\mathcal{N}} = \widetilde\Xi_{\mathcal{N}}^{\text{NLS}}$, and two embeddings
	\begin{equation*}
		 \Psi^{KG}, \Psi^{\text{NLS}}:
		 \mathbb{T}^\mathbb{Z} \times \mathcal{I}_0
		 \to \mathfrak{B}_{\sigma,r,p},
	\end{equation*}
	such that for all $\omega'
	\in\Xi'$,
	$\omega
	\in \Xi^{\text{NLS}}$ with $\widetilde{\omega}'= \widetilde{\omega}$ and  $I(0) \in \mathcal{I}_0 $,
	the functions
	\begin{align}
		&\psi^{KG}_{c,\omega'}(t):=\Psi^{KG}((\omega'+c^2) t, I(0)) \label{solu1}\\
		&\varphi^{\text{NLS}}_{\omega}(t) :=\Psi^{\text{NLS}}(\omega t, I(0)) \label{solu2}
	\end{align}
	solve respectively \eqref{NLKG2} and \eqref{NLS}. Furthermore, for any $ a \in [0,1]$, there exists 
	$C > 0$ s.t.
	\begin{equation}\label{solu_conver}
		\sup_{t \in \mathbb{R}} 
		\sup_{\substack{\omega' \in\Xi', \omega\in \Xi^{\text{NLS}} 
		\\	\widetilde{\omega}' = \widetilde{\omega}} }
	 \sup_{I(0) \in  \mathcal{I}_0 }
		\left\| e^{-\mathrm{i}c^2 t} \psi_{c,\omega'}^{KG}(t) - \varphi_{\omega}^{\text{NLS}}(t) \right\|_{r,p-4a} \leq \frac{C}{c^{1.5a}},
	\end{equation}
	where
		\begin{equation}\label{mes_comm}
		\text{meas}\left( {\Pi'}\setminus {\Xi'} \right)  = O(\gamma')
		\qquad\text{and}\qquad
		\text{meas}\left( {\Pi}^{\text{NLS}} \setminus \Xi^{\text{NLS}} \right)  = O(\gamma').
	\end{equation}

Moreover, consider an increasing sequence $\{c_m\}_{m\ge 1}$ fulfilling
	\begin{equation}\label{seq_cn}
	c_m \ge \max \left\lbrace  \epsilon^{-0.5}, m^2\right\rbrace,\quad\forall m\ge 1.
	\end{equation}
There exists a set $\Xi_{\infty}$
	such that for all $\omega\in \Xi_{\infty}$ and $ I(0) \in \mathcal{I}_0$,  the following limit holds:
	\begin{equation}\label{solu_conver1}
		\lim_{m \to \infty} \sup_{t \in \mathbb{R}} \left\| e^{-\mathrm{i}c_m^2 t} \psi_{c_m,\omega}^{KG}(t) - \varphi_\omega^{\text{NLS}}(t) \right\|_{r, p-4a} = 0.
	\end{equation}
	Let
		\begin{equation}\label{260513-6}
			{\Pi}'_{\infty}:=\lim_{m \to \infty}	{\Pi}'(c_m),
		\end{equation}
		then
	\begin{equation}\label{mes_comm1}
			\text{meas}\left({\Pi}'_{\infty} \setminus\Xi_{\infty} \right)  = O(\gamma')
		\qquad\text{and}\qquad
		\text{meas}\left({\Pi}^{\text{NLS}} \setminus\Xi_{\infty} \right) = O(\gamma').
	\end{equation} 
	
	\end{theorem}

\subsection{Comments and Outline of the proof.}\label{comments}
We provide several comments on the main difficulties and our technical strategies.


To establish the $c$-uniform existence of almost-periodic solutions for the NLKG equation, one must address specific structural obstructions in the infinite-dimensional phase space:


	\textbf{(i) Loss of uniform decay and dual asymptotic regimes.} 
	
	To achieve a measure estimate that is uniformly valid for all $c \in [1, \infty)$, the non-resonance conditions must hold uniformly, meaning that the lower bounds controlling the small denominators must be strictly independent of $c$. Note that the frequencies for equation \eqref{NLKG} satisfy
	$$\left| \omega_n-c\sqrt{c^2+n^2} \right| = O\left(  \frac{c}{\sqrt{c^2+n^2}}\right) \quad\text{as}\quad n\to \infty.$$For a fixed $c$, it is easy to observe that$$O\left( \frac{c}{\sqrt{c^2+n^2}} \right) = O\left(\frac{1}{n}\right).$$
	This $O(1/n)$ decay effectively breaks the rigid linear alignment of the frequencies—a property crucial for controlling resonances in wave-type equations. However, this decay rate actually depends explicitly on $c$. In the singular limit $c \to \infty$, the remainder degenerates:$$O\left( \frac{c}{\sqrt{c^2+n^2}} \right) = O(1).$$Consequently, the uniform decay with respect to $n$ is completely lost. Notably, this pathology mirrors the difficulties encountered in the one-dimensional derivative NLW equation \cite{Berti2013DNLW}, where the frequencies satisfy$$|\omega_n - n| = O(1) \quad \text{as} \quad n \to \infty.$$In such contexts, the quasi-Töplitz property is standardly employed to control frequency shifts. Nevertheless, in our framework, the explicit parameter dependence renders traditional quasi-Töplitz structures insufficient to guarantee uniform estimates.
	
	\textit{Our strategy:} We propose a unified non-resonance framework that captures the dual asymptotic behavior of the spectrum. For small $c$, we utilize the $O(1/n)$ decay to bound the remainder effects. Conversely, for large $c$, we employ the asymptotic expansion $c\sqrt{c^2+n^2} - c^2 \sim \frac{1}{2}n^2$ to recover the quadratic separation property characterizing the NLS equation. By bridging these two regimes, we construct a specialized set of non-resonance conditions that ensures uniform measure estimates across the entire parameter range, the detailed construction of which is carried out in Section \ref{sec 5}.

\textbf{(ii)Double degeneracy and momentum selection.} 

Under periodic boundary conditions, the normal frequencies are doubly degenerate. In the KAM scheme, this degeneracy generates non-integrable resonant terms $q_n \bar{q}_{-n}$, which obstruct linear stability. Indeed, as pointed out in \cite{CY2000}, ``\textit{when all the eigenvalues are double, one should not expect the quasi-periodic solutions obtained to be linearly stable}''. Conventional KAM schemes fail to effectively eliminate such widespread non-integrable resonances.

\textit{Our strategy:} Inspired by a novel transformation originally developed for the beam equation \cite{MR4546670}, we introduce a customized symplectic coordinate transformation that exploits the inherent translation invariance of the system. A profound consequence of this method is that it equips the perturbation with a stringent momentum conservation condition ($\sum \sigma_i n_i = 0$), a rigid algebraic structure that is inherently characteristic of the NLS equation. Under this precise selection rule, the dangerous degenerate terms vanish identically in the new coordinates, eliminating the major structural obstruction prior to the KAM iteration.

To establish the global-in-time uniform convergence between the full-dimensional tori of the two distinct PDEs, we address an inherent analytic obstruction arising from the singular limit:

\textbf{(iii)The infinite-dimensional alignment obstruction.}

At high frequencies, the phase mismatch between the NLKG and NLS equations exhibits a virulent polynomial divergence of $O(n^4/c^2)$. In a finite-dimensional setting, where this mismatch vanishes as $c \to \infty$, such phase drifts can typically be aligned by modulating internal frequencies. However, in our full-dimensional regime, the bounded parameter space $[0,1]$ cannot accommodate the simultaneous alignment of infinitely many diverging phases.
Attempting a complete frequency alignment for infinitely many modes is topologically unviable.

\textit{Our strategy:} Since complete phase alignment is unviable, we outsmart this parameter barrier through a non-local decoupling. We implement a frequency-domain truncation $\mathcal{N} \sim c^{0.4}$, enforcing strict parameter alignment exclusively for the low-frequency modes $|n| \le \mathcal{N}$ to render their linear phase drifts identically zero. For the unaligned high-frequency modes $|n| > \mathcal{N}$, instead of forcing an impossible alignment, we project the geometric comparison into a lower-order Sobolev space $H^{p-4}$. The deliberate loss of 4 spatial derivatives yields an exact $n^{-4}$ topological weight, which precisely neutralizes the $O(n^4/c^2)$ divergence of the frequency mismatches. This spatial index discrepancy naturally generates a powerful decay weight $\mathcal{N}^{-4} \sim c^{-1.6}$ at the truncation boundary, suppressing the infinite-dimensional error accumulation and yielding a uniform convergence rate of $O(c^{-1.5})$, the detailed proof of which is carried out in Section \ref{sec 4}.

The remainder of this paper is organized as follows. 
Section \ref{sec 2} details the structure of the NLKG Hamiltonian, highlighting its inherent correlation with the NLS model. 
Section \ref{sec 1} introduces the representation and norms of the Hamiltonian, along with estimates for Poisson brackets and vector fields. 
In Section \ref{sec 3}, we establish the core KAM theorem by deriving the homological equations. 
Section \ref{sec 4_5} details the inductive step of our iterative scheme, proving the persistence of the structural decomposition.
Section \ref{sec 4} proves the uniform existence of almost-periodic solutions (Theorem \ref{th1}) and the limiting convergence of the solutions (Theorem \ref{th2}).
Finally, Section \ref{sec 5} provides the measure estimates, verifying that the non-resonance conditions (\ref{NR_KG1})-\eqref{NR_NLS} are satisfied for almost all frequencies.


\section{Structure of the {NLKG} Hamiltonian}\label{sec 2}

To better highlight the structure of the NLKG Hamiltonian and reveal its inherent correlation with the  NLS Hamiltonian, we first recall the Hamiltonian of NLS equation, which will serve as a guiding model for the decomposition
introduced later.

\subsection{The Hamiltonian of the NLS equation}\label{sec 2.1}

Recall that the equation \eqref{NLS}, where $V*$ is the Fourier multiplier defined by
$$\widehat{V*u}(n)= V_n \widehat{u}(n)$$ 
and $\left( V_{n}  \right)_{{n}
	\in\mathbb{Z}}\in\mathcal{V}$ given by \eqref{potential} are the Fourier coefficients of $V$.
Consider the Sturm Liouville problem
\begin{equation*}
	-\partial_{x x} \phi_{n}+V * \phi_{n}=\xi_{n} \phi_{n}
\end{equation*}
with periodic boundary conditions on $\mathbb{T}$: it is well known that the eigenvalues and the eigenfunctions are respectively
$$\xi_{n}=n^{2}+V_{n}\qquad\text{and}\qquad
\phi_{n} =\sqrt{2/\pi}\, e^{\mathrm{i} nx}.$$
And the eigenvalues $\left(\xi_{n}\right)_{n\in\mathbb{Z}} $ are distinct and the eigenfunctions  $\left(\phi_{n}\right)_{n\in\mathbb{Z}}$ 
form an orthonormal basis of $L^{2}(\mathbb{T})$.

Then the Hamiltonian associated with equation \eqref{NLS} is
\begin{align*}
	H^{\text{NLS}}( \varphi, \bar{\varphi}) 
	&=\frac{1}{2} \int_{\mathbb{T}} \left|  \varphi_x \right|^2 +\left(V\ast \varphi \right)   \bar \varphi \mathrm{~d} x
	+\frac{3\epsilon}{8}\int_{\mathbb{T}}
	\left|\varphi \right|^4
	\mathrm{~d} x,
\end{align*}
and the corresponding  equations of motion are
\begin{equation*}
\varphi_t=	\mathrm{i}\partial_{\bar{\varphi}}H^{\text{NLS}},
		\quad 
		{\bar{\varphi}}_t=-\mathrm{i}\partial_{\varphi}H^{\text{NLS}}.
\end{equation*}
Expanding the solution into Fourier series
$	\varphi=\sum_{n\in\mathbb{Z}} {z_n} \phi_{n},
$
we rewrite the above Hamiltonian in the form
\begin{equation}\label{Ham_NLS1}
	H^{\text{NLS}}(z,\bar z)=N^{\text{NLS}} (z,\bar z) +R^{\text{NLS}}(z,\bar z).
\end{equation}
The integrable part $N^{\text{NLS}}$ is given by \eqref{integrable_NLS},
while the perturbation reads
\begin{align}\label{pertur_NLS}
	R^{\text{NLS}}(z,\bar z)
	&=\sum_{\bm{n}=(n_1,\cdots,n_4) \in \mathbb{Z}^4}
	R_{\bm{n}}^{\text{NLS}}\,
	z_{n_{1}}	z_{n_{2}} 	\bar{z}_{n_{3}} \bar{z}_{n_{4}}\nonumber\\
	&=:\sum_{\bm{n}\in \mathbb{Z}^4,   \bm{\sigma}=(\sigma_1,\cdots,\sigma_4)
		\in\{ \pm 1  \}^4}
	R_{\bm{n},\bm{\sigma}}^{\text{NLS}}\,
	z^{\sigma_1}_{n_{1}}\cdots z^{\sigma_4}_{n_{4}},
\end{align}
where $z_n^{+1}=z_n$, $z_n^{-1}=\bar z_n$ (similar to $\phi_n^{+1},\phi_n^{-1}$),
\begin{equation*}
	R_{\bm{n}}^{\text{NLS}}
	=\frac{3\epsilon}{8}\int_{\mathbb{T}} \phi_{n_1}\phi_{n_2}
	\bar{\phi}_{n_3} \bar{\phi}_{n_4}\, \mathrm{d} x
	\qquad\text{and}\qquad
	R_{\bm{n},\bm{\sigma}}^{\text{NLS}}
	=\frac{\epsilon}{16}\int_{\mathbb{T}} \phi_{n_1}^{\sigma_1}\cdots  \phi_{n_4}^{\sigma_4}\, \mathrm{d} x,
\end{equation*}
which implies that 
\begin{equation}\label{mome_NLS}
R_{\bm{n},\bm{\sigma}}^{\text{NLS}}=\frac{\epsilon}{16}
\qquad\text{if}\qquad
	\mu(\bm{n},\bm{\sigma})=0\quad\text{and}\quad\nu(\bm{n},\bm{\sigma})=0;
\end{equation}
and is zero otherwise,
where
\begin{equation*}
	\mu(\bm{n},\bm{\sigma}):=	\sigma_1\cdot n_1+\cdots+\sigma_4\cdot n_4
	\qquad\text{and}\qquad
	\nu(\bm{n},\bm{\sigma}):=\sigma_{1}+\cdots+\sigma_{4}.
\end{equation*}
This shows that the perturbation $R^{\text{NLS}}$ satisfies both the momentum and mass conservation conditions.

\subsection{The Hamiltonian of the NLKG equation}\label{sec 2.2}

We now turn to the NLKG equation \eqref{NLKG}.
Introducing the new variables $v:=u_t/c^2$, the equation  can be written as
\begin{equation*} 
	\begin{cases}
		u_t=c^2v,\\
		v_{t}= u_{x x}-c^{2} u-V * u - \epsilon  u^3.
	\end{cases}
\end{equation*}
As in the NLS case, the above equation  can be expressed as the following Hamiltonian
\begin{equation*} 
	H(v, u) =\frac{1}{2} \int_{\mathbb{T}}
	u_{x} ^{2}
	+c^{2} \left( v^{2}+u^{2}\right) 
	+ (V*u) u  \mathrm{~d} x
	+\frac{\epsilon}{4}\int_{\mathbb{T}} u^4 \mathrm{~d} x.
\end{equation*}
Introducing the complex variable $\psi$ defined in \eqref{coor_tf}, the above Hamiltonian becomes
\begin{align}\label{MR-Ham1}
	H(\psi, \bar{\psi})
	&=c^{2} \int_{\mathbb{T}}
	\bar{\psi} D_c \psi
	\mathrm{~d} x
	+\frac{\epsilon }{16}\int_{\mathbb{T}} \left( D_c^{-\frac{1}{2}} 
	\left( \psi+\bar{\psi} \right) \right)^4 
	\mathrm{~d} x
\end{align}
and satisfies the equations of motion
\begin{equation*}
	\psi_t=\mathrm{i} \partial_{ \bar{\psi}}H,\quad 	\bar{\psi}_t=-\mathrm{i} \partial_{\psi} H.
\end{equation*}
Passing to Fourier coordinates $(z,\bar z)=(z_n,\bar z_n)_{n\in\mathbb{Z}}$, we obtain
\begin{align}\label{Ham}
	H(z,\bar z)=N (z,\bar z) +R(z,\bar z),
\end{align}
endowed with the canonical symplectic structure
$\mathrm{i} \sum_{{n}\in\mathbb{Z}} \mathrm{d} \bar{z}_{n} \wedge 
\mathrm{d} {z}_{n}.$
The integrable part $N$ is given by \eqref{integrable}
and the perturbation is
\begin{align}\label{pertur_NLKG}
	R(z,\bar z)
	=\sum_{\bm{n} \in \mathbb{Z}^4,\, \bm{\sigma} \in\{ \pm 1  \}^4}
	R_{\bm{n},\bm{\sigma}}\,
	z^{\sigma_1}_{n_{1}}\cdots z^{\sigma_4}_{n_{4}},
\end{align}
where
\begin{equation*}
	R_{\bm{n},\bm{\sigma}}
	=\frac{\epsilon}{16}
	d_{n_{1}}\cdots d_{n_{4}}
	\int_{\mathbb{T}} \phi_{n_1}^{\sigma_1}\cdots  \phi_{n_4}^{\sigma_4}\, \mathrm{d} x,
	\quad 
	d_n= \left(\frac{c}{\sqrt{c^{2}+n^{2}+V_{n} }}\right)^{1 / 2},
\end{equation*}
which implies that 
\begin{equation}\label{mome1}
R_{\bm{n},\bm{\sigma}}=\frac{\epsilon}{16} 
d_{n_{1}}\cdots d_{n_{4}}
\quad \text{if}\quad	\mu(\bm{n},\bm{\sigma})=0;
\end{equation}
and is zero otherwise.

\subsection{Decomposition of the NLKG Hamiltonian}\label{sec 2.3}


In this subsection, we decompose the Hamiltonian \eqref{Ham} following the approach of {\cite[Subsection 2.2]{bambusi2025NLKG}}.
This decomposition allows us to separate different components of the nonlinearity,
which will be treated differently in the subsequent analysis.

%

We first introduce the action of the Gauge group $\{U_\alpha\}_{\alpha\in\mathbb{T}}$
\[ z' = U_\alpha(z) \iff z'_n = e^{\mathrm{i}\alpha} z_n, \quad n \in \mathbb{Z}. \]
\begin{definition} {\cite[Definition 2.1]{bambusi2025NLKG}}\label{def_NG}
	Let $R \in {C(\mathcal{U}; \mathbb{C})}$ for some open set $\mathcal{U} \subseteq \mathfrak{B}_{p}$.\\
	\textbf{(1)} We define the \textit{Gauge invariant part} of $R$ as
	\[ \Pi_{G} R(z) := \frac{1}{2\pi} \int_{\mathbb{T}} R(U_\alpha (z)) \mathrm{~d}\alpha. \]
	The \textit{non-Gauge invariant part} of $R$ is $\Pi_{NG} R = R - \Pi_{G} R$.\\
	\textbf{(2)} A Hamiltonian $R$ will be said to be \textit{Gauge invariant} if $R = \Pi_{G} R$, while it will be said to be \textit{non-Gauge invariant} if $\Pi_{G} R = 0$.
\end{definition}

With this preparation, we decompose the Hamiltonian \eqref{Ham} as
\[ H(z,\bar{z}) =  H^{NG}(z,\bar{z}) + H^{\text{NLS}}(z,\bar{z})+ H^{{h}}(z,\bar{z}) \;, \]
where $H^{NG}$ is non-Gauge invariant, $ H^{\text{NLS}}$ is given by \eqref{Ham_NLS1}, and the vector field of
$H^{{h}}$ locally uniformly vanishes as $h:=1/c^2 \rightarrow 0$. 
We now describe this decomposition by treating separately the linear and nonlinear parts.
\medskip

\noindent\textbf{Linear part.}
We begin with the integrable part $N$ defined in \eqref{integrable}.

 Then one has
	\begin{equation}
		\lambda_{n}=c \sqrt{c^2+n^2+V_n}=c^2+\frac{1}{2}\left( n^2+ V_n\right)+ \lambda_{n}^{h},
	\end{equation}
	where
	\begin{equation}
		\lambda_{n}^{h}=\frac{-n^4}{2\left(c+\sqrt{c^2+n^2+V_n}\right)^2 }.
	\end{equation}
Consequently,
	\[ N (z,\bar{z})= N^{NG}(z,\bar{z}) + N^{\text{NLS}}(z,\bar{z}) + N^h(z,\bar{z}), \]
	where $N^{\text{NLS}} $ is given by \eqref{integrable_NLS},
	\[ N^{NG}(z,\bar{z}) = 
	\sum_{n\in\mathbb{Z}} c^2  |z_{n}|^2
	\quad\text{and}\quad  
	N^h(z,\bar{z}) = 
	\sum_{n\in\mathbb{Z}} \lambda_{n}^{h}  |z_{n}|^2.\]

\medskip

\noindent\textbf{Nonlinear part.}
We now turn to the perturbation $R$ defined in \eqref{pertur_NLKG}.

	Then it admits the following decomposition
	\begin{equation}\label{260421-1}
		 R(z,\bar{z}) =  R^{{NG}}(z,\bar{z}) + R^{\text{NLS}}(z,\bar{z}) + R^h(z,\bar{z}) , 
	\end{equation}
	where $R^{\text{NLS}}(z,\bar{z})$ is given by \eqref{pertur_NLS},
	\begin{align*} 
		R^{NG}(z,\bar z)
		=\sum_{\bm{n}\in \mathbb{Z}^4,  \bm{\sigma} \in\{ \pm 1 \}^4}
		R_{\bm{n},\bm{\sigma}}^{NG}\,
		z^{\sigma_1}_{n_{1}}\cdots z^{\sigma_4}_{n_{4}},
	\end{align*}
	\begin{align*} 
		R^{h}(z,\bar z)
		=\sum_{\bm{n}\in \mathbb{Z}^4,  \bm{\sigma} \in\{ \pm 1 \}^4	}
		R_{\bm{n},\bm{\sigma}}^{h}\,
		z^{\sigma_1}_{n_{1}}\cdots z^{\sigma_4}_{n_{4}}.
	\end{align*}
We observe from \eqref{mome1} that
\begin{equation*}
	d_{n}= \left(\frac{c}{\sqrt{c^{2}+n^{2}+V_{n} }}\right)^{1 / 2}
	=1+ d_{n}'
	\quad\text{with}\quad
	\left| d_{n}' \right| \le \frac{n^2}{c^2+n^2},
\end{equation*}
which implies that $ d_{n}'\to 0 $ as $c\to \infty$,
then we define
		\begin{equation}\label{mome_h}
	R^h_{\bm{n},\bm{\sigma}}=\frac{\epsilon}{16}\cdot
	\left( d_{n_{1}}\cdots d_{n_{4}}-1\right) 
	\qquad\text{for}\qquad
	\mu(\bm{n},\bm{\sigma})=0;
\end{equation}
 and zero otherwise.
	Moreover, by Definition \ref{def_NG}, we set
	\begin{equation}\label{mome_NG}
	R_{\bm{n},\bm{\sigma}}^{NG} =\frac{\epsilon}{16}
		\qquad\text{for}\qquad
		\nu(\bm{n},\bm{\sigma})\neq 0
		\quad\text{and}\quad
		\mu(\bm{n},\bm{\sigma})=0;
	\end{equation}
	and zero otherwise.
which implies that
 $\Pi_{{G}} R^{NG} = 0$.
Finally, we define $R^{\text{NLS}}$ as the remaining part satisfying
\begin{equation*} 
	R_{\bm{n}, \bm{\sigma}}^{\text{NLS}} = R_{\bm{n}, \bm{\sigma}} - R_{\bm{n}, \bm{\sigma}}^{\text{NG}} - R_{\bm{n}, \bm{\sigma}}^h = \frac{\epsilon}{16} \qquad \text{for}\qquad \mu(\bm{n}, \bm{\sigma}) = \nu(\bm{n}, \bm{\sigma}) = 0;
\end{equation*}
and zero otherwise.

\section{The Norm of the Hamiltonian}\label{sec 1}

%
%

\subsection{Representation and norms of the Hamiltonian}

Consider the following Hamiltonian 
\begin{equation}\label{pertur1}
	R(z,\bar z)=\sum_{a,k,k'\in\mathbb{N}^{\mathbb{Z}}}R_{akk'} I(0)^{a}z^k\bar z^{k'},
\end{equation}
where
\begin{equation*}
	I(0)^{a}z^k\bar z^{k'}
	:=\prod_{{n}\in\mathbb{Z}}I_{n}(0)^{a_{n}}z_{n}^{k_{n}}\bar z_{n}^{k_{n}'}
\end{equation*}
are the monomials, with $R_{\bm{n}, \bm{\sigma}}$ as their corresponding coefficients.

\begin{definition}
Given a monomial $I(0)^{a}z^k\bar z^{k'}$ with $a,k,k'\in\mathbb{N}^{\mathbb{Z}}$, we say it satisfies:\\
\textbf{(1)} the \textit{momentum conservation} condition if
\begin{equation}\label{mome}
	\mu(k,k'):=\sum_{n\in\mathbb{Z}} (k_n-k_n') n = 0;
	\end{equation}
\textbf{(2)} the \textit{mass conservation} condition if
\begin{equation}\label{mass}
	\nu(k,k'):=\sum_{n\in\mathbb{Z}} (k_n-k_n') = 0.
	\end{equation}
\end{definition}

%


For any monomial $	I(0)^{a}z^k\bar z^{k'}$ with $a,k,k'\in\mathbb{N}^{\mathbb{Z}}$,
denote
\begin{align*}
	\mbox{degree}\ I(0)^{a}z^k\bar z^{k'}
	:=\sum_{{n}\in\mathbb{Z}}
	\left( 2a_{n}+k_{n}+k_{n}'\right)
\end{align*}
and
$(n_i^*(a,k,k'))_{i\geq1}$ the decreasing rearrangement of
\begin{equation*}
	\left\{|n|\ \mbox{where $n$ is repeated}
	\left(2a_{n}+k_{n}+k_{n}'\right)   \mbox{times}\right\}.
\end{equation*}

\begin{remark}\label{Re2.2}
	Consider the Hamiltonian $R(z,\bar z)$ given by (\ref{pertur1}),  we always assume that the monomials $I(0)^{a}z^k\bar z^{k'}$ satisfy the momentum conservation condition and
	$$\mbox{degree}\ I(0)^{a}z^k\bar z^{k'} \geq 4.$$
\end{remark}



We observe that the coefficients $R_{\bm{n}, \bm{\sigma}}$ and $R_{\bm{n}, \bm{\sigma}}^h$ of the Hamiltonians $R$ and $R^h$ exhibit different decay properties relative to $c$ (see \eqref{mome1} and \eqref{mome_h}). Incorporating this observation into the Hamiltonian norm established in \cite{CY2021}, we define the following  weighted $\ell^{\infty}$-norm.
\begin{definition}(\textbf{The weighted $\ell^{\infty}$ norm})\label{def_norm}
	For any $2<\sigma\le3$, $\rho>0$, consider the Hamiltonian
	$R$
	given by \eqref{pertur1}.
	We define its the weighted $\ell^\infty$ norms depending on the range of the parameter $c$:\\
	\textbf{(1)}	For $c \in [1, +\infty]$, 
	  \begin{equation}\label{norm+_1}
	  	\|R\|^{c,+}_{\sigma,\rho}:=
	  	\sup_{a,k,k'\in\mathbb{N}^{\mathbb{Z} }}
	  	\frac{   \left( \prod_{n\in\mathbb{Z}} 	  		\left\langle n/c\right\rangle^{\frac{1}{2}
		\left(2a_{n}+k_{n}+k_{n}' \right)} \right)  |R_{akk'}|}{e^{\rho\left(\sum_{{n}\in\mathbb{Z}}
	  	\left(2a_{n}+k_{n}+k_{n}'\right)
	  	\ln^{\sigma} \lfloor n \rfloor-2\ln^{\sigma} \lfloor n_1^*\rfloor\right)}};
	  \end{equation}
	  	\textbf{(2)}	For $c \in [1, +\infty)$, 
	  	\begin{equation}\label{norm+_3}
	  	\|R\|^{\ast,+}_{\sigma,\rho}:=
	  	\sup_{a,k,k'\in\mathbb{N}^{\mathbb{Z} }}
	  	\frac{   \left(  \sqrt{c}/n^*_1\right)^4
	  		\cdot |R_{akk'}|}{
	  		e^{\rho\left(\sum_{{n}\in\mathbb{Z}}
	  	\left(2a_{n}+k_{n}+k_{n}'\right) \ln^{\sigma} \lfloor n \rfloor-2\ln^{\sigma} \lfloor n_1^*\rfloor\right)}},
	  \end{equation}
	where  $\left\langle x\right\rangle  
	:=\sqrt{1+x^2}$ and  $n_1^*=n_1^*(a,k,k')$.
\end{definition}
\begin{remark}
	For the norm  \eqref{norm+_1}, when  {$c = +\infty$}, the weight factor $\langle n/c \rangle= 1$,  then it reduces to the standard weighted $\ell^\infty$ norm:
	 \begin{equation}\label{norm+_2}
		\|R\|^{{\infty},+}_{\sigma,\rho}:=
		\sup_{a,k,k'\in\mathbb{N}^{\mathbb{Z} }}
		\frac{|R_{akk'}|}{e^{\rho\left(\sum_{{n}\in\mathbb{Z}}\left(2a_{n}+k_{n}+k_{n}'\right) \ln^{\sigma} \lfloor n \rfloor-2\ln^{\sigma} \lfloor n_1^*\rfloor\right)}}.
	\end{equation}
	
	Specifically, the norms defined in \eqref{norm+_1}, \eqref{norm+_3}, and \eqref{norm+_2} are tailored for $R$, $R^h$, and $R^{\text{NLS}}$ (or $R^{{NG}}$), respectively.
	Furthermore, the index $\sigma$ is fixed throughout the KAM iteration. For brevity, we shall drop the subscript $\sigma$ and write $\|R\|^{c,+}_{\rho}$ and $\|R\|^{*,+}_{\rho}$ instead.
\end{remark}

%

It should be pointed out that 
the expression in \eqref{norm+_1} {and \eqref{norm+_3} satisfies} $$\sum_{{n}\in\mathbb{Z}}\left(2a_{n}+k_{n}+k_{n}'\right)
\ln^{\sigma} \lfloor n \rfloor-2\ln^{\sigma} \lfloor n_1^*\rfloor> 0.$$
In fact, the following lemma provides a more precise lower bound estimate. The proof is based on Lemma 2.1 in \cite{Cong2023TheEO}.
\begin{lemma}\label{005}
	Given any $a,k,k'\in\mathbb{N}^{\mathbb{Z}}$, assume that the momentum conservation condition (\ref{mome}) is satisfied.
	Then for any $2<\sigma\le3$, one has
	\begin{equation*}		\sum_{{n}\in\mathbb{Z}}\left(2a_{n}+k_{n}+k_{n}'\right)\ln^{\sigma} \lfloor n \rfloor-2\ln^{\sigma} \lfloor n_1^*\rfloor
		\geq\frac12\sum_{i\geq 3}{\ln^{\sigma}\lfloor n_i^*\rfloor}.
	\end{equation*}
\end{lemma}

\subsection{Some properties of the norm}

Next, we give some properties of the above norm, {with all corresponding proofs deferred to Appendix \ref{sec 6}. To this end, we give the definitions of Hamiltonian vector field and Poisson bracket.}

For any given Hamiltonian $R$, define its 
\textit{Hamiltonian vector field} by
$$X_R(z,\bar z)=\mathrm{i} \left(
\frac{\partial R }{\partial \bar z},
-\frac{\partial R }{\partial z} \right).$$
And the \textit{Poisson bracket} of two Hamiltonians $R$ and $F$ is given by
\begin{equation*}
	\{R,F\}=\mathrm{i}	\sum_{{n}\in\mathbb{Z}}
	\left( \frac{\partial R }{\partial \bar z_{n}}
	\frac{\partial F }{\partial z_{n}}
	-\frac{\partial R }{\partial z_{n}}
	\frac{\partial F }{\partial  \bar z_{n}}\right).
\end{equation*}


We first provide estimates for the Hamiltonian vector field, where the coefficient \eqref{260628-1} offers an improvement over that in Lemma 2.4 of \cite{Cong2023TheEO}.
\begin{lemma}(\textbf{Hamiltonian vector field})\label{063004}
	Given a Hamiltonian $R$ defined in \eqref{pertur1},
	then for any 
	 $p\ge 4$,
	$2<\sigma\le 3$ and $\rho\in(0,\rho')$ with $\rho'=3-2\sqrt{2}$, the following estimates hold:\\
	\textbf{(1)} For the mapping  $\mathfrak{B}_{p} \to \mathfrak{B}_{p}$,
	\begin{align}
			\sup_{\|z\|_{p}\le 1}\|X_{R}\|_{p}
		\leq C(\rho)
		\|R\|^{c,+}_{2\rho};\label{vector_1} 
	\end{align}	
	\textbf{(2)} For the mapping $\mathfrak{B}_{p} \to \mathfrak{B}_{p-4}$,
	\begin{equation}\label{vector_2}
		\sup_{\|z\|_{p}\le 1}\|X_{R}\|_{p-4}
		\leq h\cdot C(\rho)
		\|R\|^{*,+}_{2\rho},
	\end{equation}
	where
\begin{equation}\label{260628-1}
	C(\rho) =\exp\left\{\frac{80}{\rho'-\rho}\cdot\exp\left\{\frac{20}{\rho'-\rho}\right\}\right\}.
\end{equation}
\end{lemma}
Secondly, we will estimate the Poisson bracket.
\begin{lemma}(\textbf{Poisson bracket})\label{010}
Let $2 < \sigma \leq 3$ and $\rho > 0$. Suppose that $\delta_1, \delta_2$ satisfy$$0 < \delta_1, \delta_2 < \min \left\{ \frac{\rho}{4}, \, 3 - 2\sqrt{2} \right\}.$$Given two Hamiltonians $R_1$ and $R_2$, the following estimates hold:\\
	\begin{equation}\label{042704}
		\left\|\{R_1,R_2\}\right\|^{c,+}_{\rho}\leq  C(\delta_1,\delta_2,\sigma)
		\left\|R_1\right\|^{c,+}_{\rho-\delta_1}
		\left\|R_2\right\|^{c,+}_{\rho-\delta_2}
	\end{equation}
	and
		\begin{equation}\label{250804-3}
			\left\|\{R_1,R_2\}\right\|^{*,+}_{\rho}\leq C(\delta_1,\delta_2,\sigma)
			\left\|R_1\right\|^{c,+}_{\rho-\delta_1}
			\left\|R_2 \right\|^{*,+}_{\rho-\delta_2},
	\end{equation}
	where $$C(\delta_1,\delta_2,\sigma)=\frac{1}{\delta_2}\exp\left\{\frac{1000}{\delta_1}\cdot\exp\left\{
	\left( \frac{100}{\delta_1}\right)^{\frac{1}{\sigma-1}} \right\}\right\}.$$
\end{lemma}
%

Finally, let $\Phi := X^t_F|_{t=1}$ denote the Hamiltonian flow generated by $F$. For any Hamiltonian $R$, we define the iterated Poisson brackets recursively as
\begin{equation}\label{Poisson}
	 R^{(0)} = R, \qquad R^{(n)} = \{R^{(n-1)}, F\} \quad\text{ for }\; n \geq 1.
\end{equation}
Then, the action of the Hamiltonian flow $\Phi$ on $R$ can be represented by the following Lie series
\begin{equation}\label{flow}
	 R \circ \Phi = \sum_{n \ge 0} \frac{1}{n!} R^{(n)}.
\end{equation}
Regarding the convergence estimates of this series, one has the following result:

\begin{lemma}(\textbf{Hamiltonian flow})\label{E1}
Let $2 < \sigma \leq 3$ and $\rho > 0$. Suppose that $\delta$ satisfies$$0 < \delta < \min \left\{ \frac{\rho}{4}, \, 3 - 2\sqrt{2} \right\}.$$
 Given a Hamiltonian $F$, 
 assume further that
	\begin{equation*}
		\left\|F\right\|^{c,+}_{\rho},
		\left\|F\right\|^{*,+}_{\rho} <C^{-1}(\delta,\sigma),
	\end{equation*}
	where 
	\begin{equation*}
		C(\delta,\sigma)=	\frac{4e}{\delta}\exp\left\{\frac{2000}{\delta}\cdot\exp\left\{\left( \frac{200}{\delta}\right)^{\frac{1}{\sigma-1}}\right\}\right\}.
	\end{equation*}
	Then for any Hamiltonian $R$, one has the following estimates:
	for any $s\ge 0$,\\
	\begin{equation}\label{260413-1}
		\left\| \sum_{n \geq s}\frac{1}{n!} R^{(n)} \right\|^{c,+}_{\rho+\delta}
		\leq C(\delta,\sigma)
		\left( \|F\|^{c,+}_{\rho}\right)^s  \|R\|^{c,+}_{\rho}
	\end{equation}
	and
   \begin{equation}\label{260413-4}
   	\left\| \sum_{n \geq s}\frac{1}{n!} R^{(n)} \right\|^{*,+}_{\rho+\delta}
   	\leq C(\delta, \sigma) \min \left\{ \left( \|F\|_{\rho}^{c,+}\right) ^s \|R\|_{\rho}^{*,+}, \ \left( \|F\|_{\rho}^{*,+} \right) ^s \|R\|_{\rho}^{c,+} \right\}.
   \end{equation}
\end{lemma}

\subsection{Another norm of the Hamiltonian and its properties.
}

During the KAM iteration process, we need to eliminate the terms that affect the existence of the invariant tori (see the KAM iteration in Section \ref{sec 3}). To conveniently represent the terms to be eliminated,  we introduce some notations $J=(J_{n})_{n\in\mathbb{Z}}$ with
\begin{equation*}
J_{n}:=I_{n}-I_{n}(0),
\end{equation*}
and for any  $k\in\mathbb{N}^{\mathbb{Z}}$,
\begin{equation*}
	|k|:=\sum_{{n}\in\text{supp}\ k}|k_{n}|\quad \mbox{with}\quad \text{supp}\ k:=\left\{n\in\mathbb{Z}:k_{n}\neq 0\right\}.
\end{equation*}
 Then we can rewrite \( R \) given by \eqref{pertur1} in the following form of
\begin{align}\label{pertur'}
	R(z,\bar z)
	=\sum_{a,b,k,k'\in\mathbb{N}^{\mathbb{Z}}}R_{abkk'} I(0)^{a}J^bz^k\bar z^{k'}.
\end{align}

In analogy with Definition \ref{def_norm}, we give another Hamiltonian norm as follows.
\begin{definition}(\textbf{The weighted $\ell^{\infty}$ norm})
	\label{083103}
	For any $2<\sigma\le3$, $\rho>0$ and a Hamiltonian $R$ given in \eqref{pertur'}.
	We define the weighted $\ell^\infty$ norms depending on the range of the parameter $c$:\\
\textbf{(1)}	For $c \in [1, +\infty]$, 
	\begin{equation}\label{norm_1}
		\|R\|^{c}_{\sigma,\rho}:=
     	\sup_{a,b,k,k'\in\mathbb{N}^{ \mathbb{Z} }}
		\frac{   \left( \prod_{n\in\mathbb{Z}} 
		\left\langle n/c\right\rangle  ^{\frac{1}{2}\left(2a_{n}+2b_{n}+k_{n}+k_{n}' \right)}\right)  |R_{abkk'}| } {e^{\rho\left(\sum_{{n}\in\mathbb{Z}}			\left(2a_{n}+2b_{n}+k_{n}+k_{n}'\right)\ln^{\sigma} \lfloor n \rfloor-2\ln^{\sigma} \lfloor n_1^*\rfloor\right)}},
	\end{equation}
	\textbf{(2)}	For $c \in [1, +\infty)$, 
	\begin{equation} \label{norm_3}
		 	\|R\|^{*}_{\sigma,\rho}:=
		 \sup_{a,b,k,k'\in\mathbb{N}^{ \mathbb{Z} }}
		 \frac{ \left(   \sqrt{c}/n^*_1\right)^4  \cdot |R_{abkk'}|   } {  e^{\rho\left(\sum_{{n}\in\mathbb{Z}}			\left(2a_{n}+2b_{n}+k_{n}+k_{n}'\right)\ln^{\sigma} \lfloor n \rfloor-2\ln^{\sigma} \lfloor n_1^*\rfloor\right)}},
	\end{equation}
		where  $\left\langle x\right\rangle  
		:=\sqrt{1+x^2}$ and	$n_1^*=n_1^*(a,b,k,k')$.	
\end{definition}

{Similarly, we shall drop the subscript $\sigma$ and write $\|R\|^{c}_{\rho}$ and $\|R\|^{*}_{\rho}$ instead}. In fact, the two sets of Hamiltonian norms given by {(\ref{norm+_1})-\eqref{norm+_3}} and {(\ref{norm_1})-(\ref{norm_3})} are respectively equivalent in the following sense, {with the corresponding proof deferred to Appendix \ref{sec6.2}}.
\begin{lemma}\label{051301}
	Given any $2<\sigma\le 3$, $\rho>\delta>0$ and a Hamiltonian $R$, one has
	\begin{equation}\label{N6}
		\|R\|^{c}_{\rho}\leq
		\exp\left\{3\left( \frac{6}{\delta}\right)  ^{\frac{1}{\sigma-1}}\cdot\exp\left\{\left(\frac6{\delta}\right)^{\frac1{\sigma}}\right\} \right\}\|R\|_{\rho-\delta}^{c,+}
	\end{equation}
	and
	\begin{equation}\label{N7}
		\|R\|_{\rho}^{c,+} \leq \frac{64}{e^2\delta^2}\|R\|^{c}_{\rho-\delta}.
	\end{equation}  
	Furthermore, these estimates remain valid for the norms { $\|\cdot\|_{\rho}^{*,+}$ and $\|\cdot\|_{\rho}^{*}$}.
	
\end{lemma}

Building upon Lemma \ref{051301}, we present an alternative version of {Lemmas \ref{063004} and \ref{E1}}:
\begin{lemma}\label{vector}
	Under the assumptions of Lemma \ref{063004}, 
	then for any Hamiltonian $R$ given by \eqref{pertur'}, one has\\
	\textbf{(1)} For the mapping  $\mathfrak{B}_{p} \to \mathfrak{B}_{p}$,
\begin{align}
	\sup_{\|z\|_{p}\le 1}\|X_{R}\|_{p}
	\leq C'(\rho)
	\|R\|^c_{\rho};\label{vector_1'} 
\end{align}	
\textbf{(2)} For the mapping $\mathfrak{B}_{p} \to \mathfrak{B}_{p-4}$,
\begin{equation}\label{vector_2'}
	\sup_{\|z\|_{p}\le 1}\|X_{R}\|_{p-4}
	\leq h\cdot C'(\rho)
	\|R\|^{*}_{\rho},
\end{equation}
where 
\begin{equation*}
	C'(\rho) =\frac{64}{e^2\rho^2}\exp\left\{\frac{80}{\rho'-\rho}\cdot\exp\left\{\frac{20}{\rho'-\rho}\right\}\right\}.
\end{equation*}
\end{lemma}

\begin{lemma}\label{E1'}
	
	Under the assumptions of Lemma \ref{E1}, then for any Hamiltonian $R$ given by \eqref{pertur'}, one has the following estimates:	 for any $s\ge 0$,\\\\
	\begin{equation}\label{flow1}
			\left\| \sum_{n \geq s}\frac{1}{n!} R^{(n)} \right\|^{c}_{\rho+3\delta}
		\leq C'(\delta,\sigma)
		\left( \|F\|^{c}_{\rho}\right)^s  \|R\|^{c}_{\rho}
	\end{equation}
	and
	\begin{equation}\label{flow4}
		\left\| \sum_{n \geq s}\frac{1}{n!} R^{(n)} \right\|^{*}_{\rho+3\delta}
	\leq C'(\delta, \sigma) 
	\min \left\{ \left( \|F\|_{\rho}^{c}\right) ^s \|R\|_{\rho}^{*}, \ \left( \|F\|_{\rho}^{*} \right) ^s \|R\|_{\rho}^{c} \right\},
	\end{equation}	
	where
	\begin{equation*}
	C'(\delta,\sigma)	= \exp\left\{\frac{3000}{\delta}\cdot\exp\left\{\left( \frac{300}{\delta}\right)^{\frac{1}{\sigma-1}}\right\}\right\}.
	\end{equation*}
\end{lemma}

	\begin{lemma}[\cite{bambusi2025NLKG}, Lemma A.3]\label{claim3}
	Let $D' \subset D$ be two open domains in the Banach space $\mathfrak{B}_{p}$. Given two Hamiltonians $F_1, F_2$, let $X_{F_i}^t$ denote the Hamiltonian flows generated by $X_{F_i}$ for $i=1,2$. For $|t| \le 1$ and $p\ge 4$, the following estimates hold:
	\begin{equation}\label{260418-1}
		\sup_{z \in D'} \| X_{F_i}^t - \mbox{id} \|_p
		\le \sup_{z \in D} \| X_{F_i} \|_p, \quad  i=1,2
	\end{equation}
	and
	\begin{equation}\label{260418-2}
		\sup_{z \in D'} \| X_{F_1}^t - X_{F_2}^t \|_p
		\le \sup_{z \in D} \| X_{F_1} - X_{F_2} \|_p 
		\exp \left( \sup_{z \in D} \| D X_{F_1} \|_p \right). 
	\end{equation}
\end{lemma}


\section{KAM theorem}\label{sec 3}

\subsection{Derivation of homological equations}\label{sec 3.1}

 To state the KAM iteration lemma, we need to introduce the homological equation following the classical KAM way,
which is parallel to the one in Section 3.1 of \cite{Cong2023TheEO}.

Consider the general Hamiltonian 
$$H=N+R\qquad\text{with}\qquad
 N = \sum_{n \in \mathbb{Z}} \widetilde{\lambda}_{n} |z_{n}|^2,$$ 
where $R$ is given by \eqref{pertur'}. 
 Our goal is to construct a transformation $\Phi$ such that
 \begin{equation}\label{H_s}
 	H\circ \Phi  = N_{+} + R_{+},
 \end{equation}
where $N_{+}$ and  $ R_{+}$
denote the modified normal form and the corresponding smaller perturbation, respectively.

To begin the construction, we need to rewrite the perturbation $R$ as 
\begin{align}\label{pertur}
	R(z,\bar z)
	= R^0(z,\bar z) +R^1(z,\bar z)+R^2(z,\bar z),
\end{align}
where
\begin{eqnarray*}
	{R}^0(z,\bar z)&=&\sum_{a,b,k,k'\in\mathbb{N}^{\mathbb{Z}},|b|=0\atop\text{supp}\, k\bigcap \text{supp}\, k'=\emptyset}R_{abkk'} I(0)^{a}J^bz^k\bar z^{k'},\\
	{R}^1(z,\bar z)&=&\sum_{a,b,k,k'\in\mathbb{N}^{\mathbb{Z}},|b|=1\atop
		\text{supp}\, k\bigcap \text{supp}\, k'=\emptyset}R_{abkk'}	I(0)^{a}J^bz^k\bar z^{k'},\\
	{R}^2(z,\bar z)&=&\sum_{a,b,k,k'\in\mathbb{N}^{\mathbb{Z}},
		|b|=2\atop\rm{no \ assumption}}
	R_{abkk'} I(0)^{a}J^bz^k\bar z^{k'}.
\end{eqnarray*}

We aim to eliminate the terms $R^0$ and $R^1$ via a coordinate transformation $\Phi$, defined as the time-1 map $X_F^t \big|_{t=1}$ of a Hamiltonian vector field $X_F$. Let $F = F^0 + F^1$, where each $F^m$ takes the same form as $R^m$
\begin{eqnarray*}
	{F}^m=\sum_{a,b,k,k'\in\mathbb{N}^{\mathbb{Z}},|b|=
		m\atop\text{supp}\, k\bigcap \text{supp}\, k'=\emptyset}F_{abkk'} I(0)^{a} J^b z^k\bar z^{k'},
	\quad m=0,1.
\end{eqnarray*}
The homological equations then become
\begin{equation}\label{4.27}
	\{N,{F}\}+R^{\le 1}=\left[ R^{\le 1} \right]
	\quad\text{with}\quad
	R^{\le 1}:= R^0+R^1,
\end{equation}
 and
the resonant terms
\begin{equation}\label{051502}
	\left[ R^{\le 1}\right] =\sum_{a,b\in\mathbb{N}^{\mathbb{Z}},|b|\le 1}R_{ab00} I(0)^{a}J^b.
\end{equation}
Thus, the coefficients of the solutions $F$ for the equations (\ref{4.27}) are given by
\begin{equation}\label{051304}
	F_{abkk'}=\frac{R_{abkk'}}
	{\mathrm{i}\sum_{{n}\in\mathbb{Z}}
		\left(k_{n}-k^{'}_{n}\right)
		\widetilde \lambda_{n} }
	\quad\text{with}\quad k\neq k',\, |b|\le 1.
\end{equation}

Moreover, recalling \eqref{Poisson} and \eqref{flow},
then new Hamiltonian ${H}_{+}$ has the form
\begin{align}
	H_{+}\nonumber&=H\circ\Phi\\
	&=\nonumber N+\{N,F\}+R^{\le 1} \\
	&\nonumber\quad\, +\sum_{n\ge 2}\frac{1}{n!}N^{(n)}
	+\sum_{n\ge 1}\frac{1}{n!} \left( R^{\le 1} \right)^{(n)}
	+\sum_{n\ge 0}\frac{1}{n!}\left( R^{2}\right)^{(n)}\\
	&=:\nonumber\label{051401}N_++R_+,
\end{align}
where by using \eqref{4.27}
\begin{equation}\label{051402}
	N_+=N+\left[ R^{\le 1}\right] 
\end{equation}
and
\begin{align*}
	R_+
	=\sum_{n\ge 2}\frac{1}{n!}N^{(n)}
	+\sum_{n\ge 1}\frac{1}{n!} \left( R^{\le1}\right)^{(n)}
	+\sum_{n\ge 0}\frac{1}{n!}\left( R^{2}\right)^{(n)}.
\end{align*}
By setting $\widetilde{R} := R - [R]$, we rewrite \eqref{4.27} as $\{N, F\} = -\widetilde{R^{\le 1}}$, which implies that
\begin{equation}\label{051403}
	R_+
	=\sum_{n\ge 1}
	\left( \frac{1}{(n+1)!} \left( \widetilde{R^{\le 1}}\right)^{(n)}
	+ \frac{1}{n!} \left( R^{\le1}\right)^{(n)}\right) 
	+\sum_{n\ge 0}\frac{1}{n!}\left( R^{2}\right)^{(n)}.
\end{equation}

\subsection{KAM Iteration}\label{031501}

Next, we give the precise
set-up of the $s$-th iteration parameters: for any $2<\sigma\le 3$, $r>1$ and $p\ge 4$,
\begin{itemize}
	\item$\delta_{s}=\frac{\rho_0}{(s+4)\ln^2(s+4)}$ with 
	$\rho_0= \frac{3-2\sqrt{2}}{100}$,
	
	\item$\rho_{s+1}=\rho_{s}+3\delta_s,$
	
	\item$\epsilon_s=\epsilon_{0}^{({3}/{2})^s}$, which dominates the size of the perturbation,
	
	\item$\lambda_s=\epsilon_s^{0.01}$,
	
	\item$\eta_{s+1}=\frac{1}{20}\lambda_s\eta_s$ with $\eta_0=\lambda_0$,
	
	\item$d_{s+1}=d_{s}+\frac{1}{\pi^2(s+1)^2}$ with $d_0=0$,
	
	\item$D_{s,p}=\{(z_{n})_{n\in\mathbb{Z}}:\frac{1}{2}+d_s\leq|z_{n}|e^{r\ln^{\sigma}\lfloor n \rfloor} \lfloor n \rfloor^p
	\leq1-d_s\}$.
\end{itemize}

For the sake of brevity, we shall drop the subscript $p$ and write $D_s$ instead of $D_{s,p}$ when no confusion arises.
Fixed a parameter $V^* \in\mathcal{V}\subset\mathbb{R}^{\mathbb{Z}}$ given by \eqref{potential}, denote the complex cube of size $\lambda>0$
\begin{equation*}
	\mathcal{C}_{\lambda}\left({ {V}^*}\right)=\left\{( {V}_{n})_{n\in\mathbb{Z}}\in\mathbb{C}^{\mathbb{Z}}:
	\left| {V}_{n}-{V}^*_{n}\right|\leq \lambda,n\in\mathbb{Z}\right\}.
\end{equation*}

\begin{lemma}\label{IL}(\textbf{Iterative Lemma})
Consider the Hamiltonian $H_{s}=N_{s}+R_{s}$, which is real analytic on $D_{s}\times\mathcal{C}_{\lambda_s}\left( V_s^*\right)$, with the normal form $N_s$ given by
\begin{equation}\label{231130-1'}
N_s = \sum_{n \in \mathbb{Z}} \widetilde{\lambda}_{s, n}(V)\left|z_n\right|^2, \quad  \widetilde{\lambda}_{s, n}(V) = c \sqrt{c^2+n^2+\widetilde{V}_{s, n}(V)}.	
\end{equation}
	Suppose that $ \widetilde{V}_s=\left( \widetilde{V}_{s,n}(V) \right)_{n\in\mathbb Z}$ 
	satisfies
	\begin{equation}\label{199}
		\left\| \frac{\partial \widetilde{V}_s}
		{{\partial {V}}}-Id\right\|_{{\infty}\rightarrow {\infty}}
		<d_s\epsilon_{0}^{0.1}
	\end{equation}
	(where $\left\|\cdot\right\|_{\infty\rightarrow\infty}$ is an operator norm and $\left\|\cdot\right\|_{\infty}:=\left\|\cdot\right\|_{ \ell^{\infty}}$) and
	\begin{equation}\label{IL3}
		\widetilde{V}_{s,n}({V_s^*})
		= \frac{\omega_n^2}{c^2} -c^2-n^2=:T_n,
	\end{equation}
	where  $\omega = (\omega_n)_{n \in \mathbb{Z}}\in \Pi$ satisfying \eqref{NR_KG1}.
	Furthermore, assume that  
	$R_{s}=R^0_{s}+R^1_{s}+R^2_{s}$
	 is a small perturbation satisfying
	\begin{eqnarray}
		\label{200'}&&\left\|R^0_{s}\right\|^c_{\rho_{s}}
		\leq \epsilon_{s},\\
		\label{201'}&&	\left\|R^1_{s}\right\|^c_{\rho_{s}}
	\leq \epsilon_{s}^{0.6}, \\
		\label{202'}&&\left\|R^2_{s}\right\|^c_{\rho_{s}}
		\leq (1+d_s)\epsilon_0.
	\end{eqnarray}

	Then for 
	${V}\in\mathcal{C}_{\eta_{s}}\left({V}_{s}^* \right)$
	satisfying 
	$\widetilde{V}_{s}( { {V}})
	\in\mathcal{C}_{\lambda_{s-1}}(T)$,
	there exists a real analytic symplectic coordinate transformation
	$\Phi_{s}:D_{s+1}\rightarrow D_{s}$  satisfying
	\begin{eqnarray}
		\label{203}&&
		\sup_{z\in D_{s+1}} \left\|\Phi_{s}-\mbox{id}\right\|_{p}
		\leq \epsilon_{s}^{0.5},\\
		\label{204}&&
		\sup_{z\in D_{s+1}} 
		\left\|D\Phi_{s}-\mbox{Id}\right\|_{p\rightarrow p}
		\leq \epsilon_{s}^{0.5},
	\end{eqnarray}
	such that 
	\begin{align*}
		H_{s+1}=H_{s}\circ\Phi_{s}=N_{s+1}+R_{s+1}
	\end{align*}
	and the same assumptions
	are satisfied with `$s+1$' in place of `$s$', 
	where 
	\begin{align}
		&\left\|\widetilde{V}_{s+1}-\widetilde{V}_{s}\right\|_{\infty}	\leq\epsilon_{s}^{0.5}, \label{206}\\
		&\left\|{ {V}}_{s+1}^*-{ {V}}_{s}^*\right\|_{\infty}
		\leq 2\epsilon_{s}^{0.5}\label{205},\\
		&	\mathcal{C}_{\eta_{s+1}}\left({{V}}_{s+1}^*\right) \subset \widetilde{V}^{-1}_{s+1}\left( \mathcal{C}_{\lambda_{s}}(T)\right)  \label{IL2}.
	\end{align}
\end{lemma}

\begin{proof}
	The proof can be divided into the following three steps.
The first two steps—the estimates of the homological equation solutions and the remainder terms—are more concise than those in \cite{Cong2023TheEO}, owing to the more precise non-resonance conditions \eqref{NR_KG1} and the refined representation of the new remainder terms in \eqref{flow}. In the final step, we verify that the frequency shifts are small perturbations for the current iteration by introducing the suitable norm \eqref{norm+_1}.
	
	\textbf{Step 1. Estimate on the solutions of homological equations.}

	First, we will give the lower bound on the right hand side of the  nonresonance conditions (\ref{NR_KG1})
	after truncation. 

	In fact,  since in the step $s+1$ one only needs
	\begin{eqnarray*}
		\left\|R^0_{s+1}\right\|^c_{\rho_{s+1}}\leq \epsilon_{s+1}
		\qquad\text{and}\qquad
		\left\|R^1_{s+1}\right\|^c_{\rho_{s+1}}\leq \epsilon_{s+1}^{0.6},
	\end{eqnarray*} 
	there is a saving of a factor
	\begin{equation*}
		e^{-3\delta_{s}\left(\sum_{{n}\in\mathbb{Z}}\left(2a_{n}+2b_{n}+k_{n}+k'_{n}\right)\ln^{\sigma}\lfloor n \rfloor-2\ln^{\sigma}\lfloor n_1^*\rfloor\right)}.
	\end{equation*}
	Combining with Lemma \ref{005}, it suffices to eliminate the monomials $I(0)^{a}J^bz^k\bar z^{k'}$ in $R^0_{s}, R^1_{s}$ for which
	\begin{equation*}
		e^{-\frac32\delta_{s}\sum_{i\geq3}\ln^{\sigma}\lfloor n_i^*\rfloor}\geq\epsilon_{s+1},
	\end{equation*}
	that is,
	\begin{equation}\label{M1}
		\sum_{i\geq3}\ln^{\sigma}\lfloor n_i^*\rfloor
		\leq\frac{2(s+4)\ln^2(s+4)}{3\rho_0}
		\cdot\ln\frac{1}{\epsilon_{s+1}}=: B_s,
	\end{equation}
	where $n_i^*=n_i^*\left(a,b,k,k'\right)$. 
	
	We begin show that for any $\ell\in\mathbb{Z}^{\mathbb{Z}}$ 
	with $3\le  |\ell|<\infty$,
	\begin{equation}\label{est_NR}
		\gamma \left( 
		\prod_{|n| \leq  n_3^{*}(\ell),\, \ell_n\neq 0}  
		\frac{1}{ |\ell_{n}|^5\,\lfloor n\rfloor^6}\right)^{5} \geq \gamma\epsilon_s^{0.01}.
	\end{equation}

	More concretely, let
	\begin{equation*}\label{030701}
		 N_s=
		{B_s^{1-\frac{1}{\sigma}}}\left(  {\ln B_s}\right) ^{-\sigma}.
	\end{equation*}
	Then, we have
	\begin{eqnarray}
		&&   \left( 
		\prod_{|n| \leq  n_3^{*}(\ell),\, \ell_n\neq 0}  
		\frac{1}{ |\ell_{n}|^5\,\lfloor n\rfloor^6}\right)^{5} 
		\nonumber \\
		&\geq& \exp\left\{-150 \left( \sum_{0\le |n| \leq  N_s,\ell_{n}\neq 0}\ln\left(8+|\ell_{n}|
		\right)\,\ln\lfloor n \rfloor
		+\sum\limits_{  N_s <n\le n_3^{*}(\ell)} |\ell_{n}|
		\ln \lfloor n \rfloor \right)  \right\}\nonumber \\
		&\geq& \exp\left\{-450 N_s B_s^{\frac1{\sigma}}
		-150 B_s(\ln N_s)^{^{1-\sigma}}\right\}
		\quad  (\mbox{by } 
		\ln\left(8+|\ell_{n}|
		\right)\le 3|\ell_{n}|^{\frac{1}{3}} )\nonumber\\
		&\geq& \exp\left\{-{600B_s}{\left(\ln N_s\right)^{1-\sigma}}\right\} \label {M6}.
	\end{eqnarray}
	Furthermore, one deduces from (\ref{M1}) that
	\begin{eqnarray}
		600 B_s(\ln \mathcal N_s)^{1-\sigma}
		\leq 0.01 \ln{\epsilon_s^{-1}}
		\label{022402},
	\end{eqnarray}
	where the last inequality follows from  $2<\sigma\le 3$ and $\epsilon_0$ is sufficiently small depending on $\sigma$.
	Finally, in view of \eqref{M6} and \eqref{022402}, one finishes the proof of \eqref{est_NR}.
		Hence, by using Remark \ref{l=1}, \eqref{IL3} and \eqref{est_NR}, we have for any $\ell\in\mathbb{Z}^{\mathbb{Z}}$ with $0<|\ell|<\infty$,
	\begin{eqnarray}\label{M7}
		\left|\sum_{{n}\in\mathbb{Z}}
		\ell_{n}\cdot \widetilde{\lambda}_{s,n}(V_s^*)\right|
		=	\left| \sum_{{n}\in\mathbb{Z}}
		\ell_{n}\cdot \omega_{n}
		\right|
		\geq \gamma \epsilon_s^{0.01}.
	\end{eqnarray}

	Next, we estimate on the solutions of homological equations.
	

	By using  (\ref{051304}), (\ref{200'}), (\ref{201'}), (\ref{M7}) we get
	\begin{align*}
		&\left\|F_{s}^0\right\|^c_{\rho_s}\leq {\gamma}^{-1} \epsilon_s^{-0.01}\left\|R_{s}^0\right\|^c_{\rho_s}\le\epsilon_s^{0.95} ,\\
		&	\left\|F_{s}^1\right\|^c_{\rho_s}\leq {\gamma}^{-1} \epsilon_s^{-0.01}\left\|R_{s}^1\right\|^c_{\rho_s}\le \epsilon_s^{0.55} ,
	\end{align*}
	which implies that
	\begin{equation}\label{022501}
		\left\|F_{s}\right\|^c_{\rho_s}\leq \epsilon_s^{0.54}.
	\end{equation}

%


	Finally, we estimate the coordinate transformations  \eqref{203} and \eqref{204}.
	
	Due to (\ref{vector_1'}),
	we get
	\begin{equation}\label{est_X_F}
		\sup_{z\in D_{s}}\|X_{F_s}\|_{p}
		\leq C(\rho_s)	\cdot \left\|F_{s}\right\|^c_{\rho_s}
		\leq \epsilon_{s}^{0.52}.
	\end{equation}
Therefore, setting $\Phi_{s} = X_{F_s}^t|_{t=1}$, we apply \eqref{260418-1} to obtain
$\Phi_{s} : D_{s+1} \to D_{s}$ with
\begin{align}
	&\sup_{z\in D_{s+1}} 
	\left\|\Phi_{s}-\mbox{id}\right\|_{p}
	\leq\sup_{z\in D_{s}} \|X_{F_s}\|_{p}
	<\epsilon_{s}^{0.5},\label{203-1}\\
	&\sup_{z\in D_{s+1}}  
	\left\|D\Phi_{s}-\mbox{Id}\right\|_{p\to p}
	\leq \frac{1}{d_{s+1}-d_s}
	\sup_{z\in D_{s}} \|X_{F_s}\|_{p}
	<\epsilon_{s}^{0.5}, \label{204-1} 
\end{align}
	where 
	 using the Cauchy's estimate and the fact that $$\epsilon_{s}^{0.01}\ll\frac{1}{\pi^2(s+1)^2}=d_{s+1}-d_s.$$ 
	Therefore, we have completed the proofs of \eqref{203} and \eqref{204}.

	\textbf{Step 2. Estimate on the remainder terms.} 
	
%
	First of all, recalling (\ref{051403}) and the new remainder term $R_{s+1}$ is given by
	$$R_{s+1}=\sum_{n\ge 1}\left( 
	\frac{1}{(n+1)!} \left( \widetilde{R_s^{\le 1}}\right)^{(n)}
	+ \frac{1}{n!} \left( R_s^{\le 1}\right)^{(n)}\right) 
	+\sum_{n\ge 0}\frac{1}{n!}\left( R_s^{2}\right)^{(n)},
	$$
	where 
	$$	R_s^{(0)}=R_s	\quad\text{and}\quad
	R_s^{(n)}=\left\{R_s^{(n-1)},F_s\right\}
	\quad\text{for}\quad n\ge 1
	$$
	with $F_s=F_{s}^0+F_{s}^1$.
Indeed, since $\widetilde{R_s^{\le 1}}$ is a component of 
$R_s^{\le 1}= R_s^0+R_s^1$, it suffices to consider
	\begin{equation*}
		\sum_{n\ge 1}\frac{1}{n!} 
		\left( R_s^0+R_s^1\right)^{(n)}
		+\sum_{n\ge 0}\frac{1}{n!}
		\left( R_s^{2}\right)^{(n)}
		=:R^0_{s+1}+R^1_{s+1}+R^2_{s+1}.
	\end{equation*} 
	
	Secondly, we can find that \\
	$\bullet$ the term $R^0_{s+1}$ may come from
	\begin{align}\label{1203-1}
		&\sum_{n\ge 1} \left( R_{s}^{0}\right)^{(n)},
		\qquad \sum_{n\ge 0} \{R^1_{s},F^0_{s} \}^{(n)},
		\qquad  \sum_{n\ge 1} \{R^1_{s},F^1_{s} \}^{(n)},\nonumber\\
		& \sum_{n\ge 1} \{R^2_{s},F^0_{s} \}^{(n)},
		\qquad  \sum_{n\ge 2} \{R^2_{s},F^1_{s} \}^{(n)};
	\end{align}
	$\bullet$ the term $R^1_{s+1}$ may come from
	\begin{align}\label{1203-2}
		&\sum_{n\ge 1} \left( R_{s}^{0}\right)^{(n)},
		\qquad	\sum_{n\ge 1} \left( R_{s}^{1} \right)^{(n)},
		\qquad	\sum_{n\ge 0} \{R^2_{s},F^0_{s} \}^{(n)},\nonumber\\
		&\sum_{n\ge 1} \{R^2_{s},F^1_{s} \}^{(n)};
	\end{align}
	$\bullet$ the term $R^2_{s+1}$ may come from
	\begin{align}\label{1203-3}
		\sum_{n\ge 1}\left(  R_{s}^{0}\right)^{(n)},
		\qquad	\sum_{n\ge 1} \left( R_{s}^{1}\right)^{(n)},
		\qquad	\sum_{n\ge 0} \left( R_{s}^{2}\right)^{(n)}.
	\end{align}
	
	Finally, we estimate the term $R_{s+1}$.
	Without loss of generality, consider the term $\sum_{n\ge 1}\left( R_{s}^{0}\right) ^{(n)}$.
	By utilizing (\ref{flow1}), (\ref{200'}) and (\ref{022501}), we get
	\begin{eqnarray}\label{N11}
		\left\| \sum_{n\ge 1} \frac{1}{n!}\left( R_{s}^{0}\right) ^{(n)} \right\|^c_{\rho_s+3\delta_s}
		\leq
	C'(\delta_s,\sigma) \cdot
		\|{F_s}\|^c_{\rho_s}\cdot
		\left\| R^0_{s}\right\|^c _{\rho_s }
		\leq\epsilon_s^{1.53},
	\end{eqnarray}
	where the last inequality using the fact that
	$$\exp\left\lbrace s^{ \frac{ 1}{\sigma-1 } } \right\rbrace < \left( \frac{3}{2}\right)^s
	\quad\text{for}\quad \sigma >2.$$ 
	Following the proof of (\ref{N11}), one can deduce from (\ref{1203-1})-(\ref{1203-3}) that
	\begin{align*}
		\left\| R^0_{s+1}\right\|^c_{\rho_{s+1}}
		&\leq10\epsilon_s^{1.53}\leq \epsilon_{s+1},\\
		\left\| R^1_{s+1}\right\|^c_{\rho_{s+1}}
		&\leq 10\epsilon_s^{0.93}\leq \epsilon_{s+1}^{0.6},\\
		\left\| R^2_{s+1}\right\|^c_{\rho_{s+1}}
		&\leq(1+d_s)\epsilon_0+2\epsilon_s^{0.5}\leq(1+d_{s+1})\epsilon_0.
	\end{align*}
These estimates confirm that the assumptions (\ref{200'})-(\ref{202'}) are satisfied for $R^0_{s+1},R^1_{s+1}$ and $R^2_{s+1}$.

	\textbf{Step 3. Estimate on the frequency shift.}

In this step, we need to verify the assumptions \eqref{199}, \eqref{IL3}  for \( s+1 \), as well as the estimates \eqref{206}-\eqref{IL2}.
	
	Recall that (\ref{051402}), the new normal form $N_{s+1}$ is given by
	\begin{equation*}
		N_{s+1}=N_s+\left[ R^0_{s}\right] +\left[ R^1_{s}\right] .
	\end{equation*}
Since the term $[R^0_{s}]$ in \eqref{051502} is a constant that does not affect the Hamiltonian vector field, we define the new frequencies as
\begin{equation}\label{shift1}
	\widetilde{\lambda}_{s+1,n}
	= \widetilde{\lambda}_{s,n}+
	\sum_{a\in\mathbb{N}^{\mathbb{Z}}}R_{s,ae_{n}00} I(0)^{a},
\end{equation}
where the last term is the so-called frequency shift.

		First of all, we prove  the estimate \eqref{206}.

		On one hand, by using \eqref{norm_1}, one has for any $n\in\mathbb{Z}$,          
		\begin{equation*}
			\left|  \left( \left\langle \frac{n}{c}\right\rangle  \prod_{m\in\mathbb{Z}} 
			\left\langle \frac{m}{c}\right\rangle  ^{a_{m}}  \right) 
			R_{s,ae_{n}00} \right|\le 	
			\|R_s^1\|^c_{\rho_s}  e^{2\rho_s\left(\sum_{m\in\mathbb{Z}}
				a_{m}\ln^{\sigma} \lfloor m \rfloor+ \ln^{\sigma} \lfloor n \rfloor
				-\ln^{\sigma} \lfloor m_1^*\rfloor\right)}
		\end{equation*}
		with  $\left\langle n/c\right\rangle  
		=\sqrt{1+n^2/c^2}$ and $m_1^*=m_1^*(a,e_n)$,
		which implies from \eqref{201'} 
		that
		\begin{equation}\label{shift4}
			\left|R_{s,ae_{n}00} \right|\le  
			\left\langle \frac{n}{c}\right\rangle^{-1}  \,
			\epsilon_{s}^{0.6}\, e^{2\rho_s\sum_{m\in\mathbb{Z}}a_{m}\ln^{\sigma} \lfloor m \rfloor
			}.
		\end{equation}

		On the other hand, 
		combining the estimate 
		\begin{equation}\label{shift6}
			\left|I_m(0)\right| \le e^{-2r\ln^{\sigma} \lfloor m \rfloor} \lfloor n \rfloor^{-2p}
			\le e^{-2r\ln^{\sigma} \lfloor m \rfloor} ,
		\end{equation}
		since $I(0) \in \mathcal{I}_0$,  one obtains
		\begin{align}\label{shift2}
			\left|	\sum_{a\in\mathbb{N}^{\mathbb{Z}}}
			R_{s,ae_{n}00} I(0)^{a}\right|
			&\le  \left\langle \frac{n}{c}\right\rangle^{-1}   \epsilon_{s}^{0.6}
			\sum_{a\in\mathbb{N}^{\mathbb{Z}}}
			e^{2(\rho_s-r)\sum_{m\in\mathbb{Z}}a_{m}\ln^{\sigma} \lfloor m \rfloor }\nonumber\\
			&\le  \left\langle \frac{n}{c}\right\rangle^{-1} \epsilon_{s}^{0.6}
			\sum_{a\in\mathbb{N}^{\mathbb{Z}}}
			e^{-0.1\sum_{m\in\mathbb{Z}}a_{m}\ln^{\sigma} \lfloor m \rfloor },\nonumber
		\end{align}
		where the last inequality follows from $r-2\rho_s>0.1$.

				\begin{claim}[\cite{Cong2023TheEO}, Lemma 4.3]\label{a5}
			Given any $a\in\mathbb{N}^{\mathbb{Z}}$,  assume $\sigma>2$ and $\delta\in(0,3-2\sqrt{2})$. Then we have
			\begin{equation*}
				\sum_{a\in\mathbb{N}^{\mathbb{Z}}}e^{-\delta\sum_{n\in\mathbb{Z}}a_{n}\ln^{\sigma}\lfloor n\rfloor}
				\leq\exp\left\{\frac{18}\delta\cdot \exp\left\{\left(\frac4{\delta}\right)^{\frac1{\sigma-1}}\right\}\right\}.
			\end{equation*}
		\end{claim}	
Thus, due to \eqref{shift1} and Claim \ref{a5}, we get
	\begin{equation}\label{shift5}
			\left| 	\widetilde{\lambda}_{s+1,n}
		- \widetilde{\lambda}_{s,n} \right|=
		\left|	\sum_{a\in\mathbb{N}^{\mathbb{Z}}}
	   R_{s,ae_{n}00} I(0)^{a}\right|	\le 
	\frac{c}{\sqrt{c^2+n^2}}\cdot
	 \epsilon_{s}^{0.6} \cdot C(\sigma).
	\end{equation}
		Meanwhile, note that by \eqref{231130-1'}
		\begin{align}
			\widetilde{\lambda}_{s+1,n}
			- \widetilde{\lambda}_{s,n}
			&=c\left( \sqrt{c^2+n^2+\widetilde{V}_{s+1,n}}
			-\sqrt{c^2+n^2+\widetilde{V}_{s,n}}\right)\nonumber \\
			&=\frac{c\left( \widetilde{V}_{s+1,n}-\widetilde{V}_{s,n}\right)  }{\sqrt{c^2+n^2+\widetilde{V}_{s+1,n}}+\sqrt{c^2+n^2+\widetilde{V}_{s,n}}},
			\label{lambda-}
		\end{align}
		then one deduces from \eqref{shift5} that
		\begin{equation}\label{V-}
			\left| \widetilde{V}_{s+1,n}-\widetilde{V}_{s,n}\right|
			\le 3\epsilon_{s}^{0.6}\cdot C(\sigma)
			\le\epsilon_{s}^{0.5}.
		\end{equation}
		Therefore, we finish the proof of \eqref{206},
		which implies by induction that 
		\begin{equation}\label{053190}
			\left\|\widetilde{V}_{s+1}\right\|_{\infty}
			\leq 1+{\sum_{i=0}^{s}}{\epsilon_{i}^{0.5}}.
		\end{equation}
		\begin{remark}
			
			Based on \eqref{lambda-}, the factor $\langle n/c \rangle$ in the definition \eqref{norm_1} is essential. Its inclusion is crucial for demonstrating that $V_{s+1}$ remains a small perturbation of $V_s$, and that $\|\widetilde{V}_{s}\|_\infty$ converges as $s \to \infty$, as established in \eqref{053190}.
		\end{remark}

		
		Secondly, we prove the two assumptions (\ref{199}) and (\ref{IL3}) with $s+1$.
		
		For any $n\in\mathbb{Z}$, applying Cauchy's estimate on $\mathcal{C}_{\lambda_s\eta_s}\left(V_s^*\right)$,  one gets
		\begin{eqnarray}
			\sum_{m\in\mathbb{Z}}\left|\frac{\partial \widetilde{V}_{s+1,n}}{\partial V_{m}}-\frac{\partial \widetilde{V}_{s,n}}{\partial V_{m}}\right|
			\leq\frac{\left\|\widetilde{V}_{s+1}-\widetilde{V}_{s}\right\|_\infty}{\lambda_s\eta_s}
			\label{M17}\leq\frac{\epsilon_{s}^{0.5}}{\lambda_s\eta_s},
		\end{eqnarray}
		where the last inequality follows from \eqref{206}.
		Note that
		\begin{eqnarray*}
			\lambda_s\eta_{s}=\eta_0\prod_{i=0}^{s}\left(\frac1{20}\lambda_i\right)
			=\eta_0\prod_{i=0}^{s}\left(\frac1{20}\epsilon_i^{0.01}\right)
			\geq\epsilon_{0}^{0.06\times(\frac{3}{2})^{s}},
		\end{eqnarray*}
		then for any $n\in\mathbb{Z}$ and $V\in \mathcal{C}_{\frac{1}{10}\lambda_s\eta_s}\left(V_s^*\right)$, we deduce from (\ref{199}) and (\ref{M17}) that
		\begin{align*}
			&\quad \sum_{m\in\mathbb{Z}}\left|\frac{\partial \widetilde{V}_{s+1,n}}{\partial V_{m}}-
			\delta_{nm}\right|\\
			&\leq\sum_{m\in\mathbb{Z}}
			\left| \frac{\partial \widetilde{V}_{s+1,n}}{\partial V_{m}}-\frac{\partial \widetilde{V}_{s,n}}{\partial V_{m}} \right|
			+\sum_{m\in\mathbb{Z}}
			\left| \frac{\partial  \widetilde{V}_{s,n}}{\partial V_{m}}-\delta_{nm}\right|\\
			&\leq\epsilon_{0}^{0.05\times(\frac{3}{2})^{s}}+d_s\epsilon_{0}^{0.1}
			<d_{s+1}\epsilon_{0}^{0.1},
		\end{align*}
		where $\delta_{nm}=\mbox{id}$ when $n=m$; otherwise $\delta_{nm}=0$.
		Thus, one obtains that
		\begin{equation}\label{M20}
			\left\|\frac{\partial \widetilde{V}_{s+1}}{\partial V}-
			Id\right\|_{\infty\rightarrow\infty}
			<d_{s+1}\epsilon_{0}^{0.1},
		\end{equation}
		which verifies the assumption (\ref{199}) with $s+1$.

		Furthermore, we fix $T=	\left( T_n \right)_{n\in\mathbb{Z} }$, where $T_n$ is given by \eqref{IL3}.
	By applying the inverse function theorem to the equation
		\begin{equation}\label{M21}
			\widetilde{V}_{s+1,n}\left(V\right)=T_n,
			\qquad   V\in\mathcal{C}_{\frac{1}{10}\lambda_s\eta_s}
			\left(V_s^*\right).
		\end{equation}
	and combining it with \eqref{M20}, we deduce the existence of a unique solution $V_{s+1}^*$ such that
		\begin{equation}\label{240504-1}
			\widetilde{V}_{s+1,n}\left(V_{s+1}^* \right)=T_n,
		\end{equation}
		thereby verifying the assumption (\ref{IL3}) with $s+1$. 
		\begin{remark}
			It is worth noting that the choice of \( T_n \) in \eqref{IL3} differs from that in \cite{Bourgain2005JFA}.
			, owing to the distinct relation between the frequency $\omega = (\omega_n)_{n \in \mathbb{Z}}$ and the vector $T$. Specifically, in this paper, we have
			$\omega_n = c \sqrt{c^2+n^2+T_n}$.	
		\end{remark}

		Finally, we prove the estimates (\ref{205}) and \eqref{IL2}.

		According to \eqref{IL3} and \eqref{240504-1},
		we rewrite $V_{s+1}^*-V_s^*$ as
		\begin{equation*}\label{M22}
			V_{s+1}^*-V_s^*=\left(I-\widetilde{V}_{s+1}\right)\left(V_{s+1}^*\right)-\left(I-\widetilde{V}_{s+1}\right)\left(V_s^*\right)
			+\left(\widetilde{V}_s-\widetilde{V}_{s+1}\right)\left(V_s^*\right)
		\end{equation*}
		and then by using \eqref{206} and (\ref{M20})
		\begin{equation*}\label{M23}
			\left\|V_{s+1}^*-V_s^*\right\|_{\infty}\leq d_{s+1}\epsilon_{0}^{\frac{1}{10}}\left\|V_{s+1}^*-V_s^*\right\|_{\infty}+\epsilon_{s}^{0.5}.
		\end{equation*}
		Note that $d_{s+1}\epsilon_{0}^{0.1} \le 1/2$, then we can obtain
		\begin{equation*}
			\left\|V_{s+1}^*-V_s^*\right\|_{\infty}\leq 2\epsilon_{s}^{0.5},
		\end{equation*}
		which verifies the estimate (\ref{205}).

		Moreover, if $V\in \mathcal{C}_{\frac{\eta_{s}}{2}}\left(V_{s+1}^*\right)$, for any $n\in\mathbb{Z}$, we derive from Cauchy's estimate that
		\begin{eqnarray}\label{M11}
			\sum_{m\in\mathbb{Z}}
			\left|\frac{\partial \widetilde{V}_{{s+1},n}}{\partial V_{m}}\right|
			\leq\frac{5}{\eta_{s}}
			\left\|\widetilde{V}_{s+1}\right\|_\infty
			<\frac{10}{\eta_{s}},
		\end{eqnarray}
		where the last inequality is based on
		\eqref{053190}.
		Let $V \in \mathcal{C}_{\frac{1}{10}\lambda_s\eta_s}\left(V_{s+1}^*\right)$, then by (\ref{M11}) for any $n\in\mathbb{Z}$,
		\begin{align*}
			\left| \widetilde{V}_{s+1,n}(V)-T_n
			\right|
			&=\left| \widetilde{V}_{s+1,n}(V)-
			\widetilde{V}_{s+1,n}\left(V_{s+1}^*\right)\right|\\
			&\leq\sup_{\mathcal{C}_{\frac{1}{10}\lambda_s\eta_s}\left(V_{s+1}^*\right)}\left\|\frac{\partial  \widetilde{V}_{s+1}}{\partial V}\right\|_{ \infty\rightarrow  \infty} \left|V_n-V_{s+1,n}^*\right|\\
			&<\frac{10}{\eta_{s}}\cdot\frac{1}{10}\lambda_{s}\eta_{s}
			=\lambda_{s},
		\end{align*}
		which means
		\begin{equation*}
			\widetilde{V}_{s+1}\left(\mathcal{C}_{\frac{1}{10}\lambda_s\eta_s}\left(V_{s+1}^*\right)\right)\subseteq \mathcal{C}_{\lambda_{s}}\left(T \right)
		\end{equation*}
		and the estimate \eqref{IL2} holds.
		
		To sum up, we complete the proof of the iterative lemma.
	\end{proof}

\subsection{KAM Theorem}\label{120301}
We now proceed to prove the convergence of Lemma \ref{IL}, thereby arriving at the following theorem:
	\begin{theorem}(\textbf{KAM Theorem})\label{KAMth}
 For any $2<\sigma\le 3$, $p\ge 4$, $r>1$ and $c\in [1,\infty)$, suppose the Hamiltonian
 $H=N+R$
 	is real analytic on $D_{0}\times\mathcal{V} $,
 	where the integrable term
 	$$N(z,\bar{z})=\sum_{{n}\in\mathbb{Z}}
 	\lambda_{n} (V) |z_{n}|^2 \quad \text{with}\quad
 	\lambda_{n} (V)= c \sqrt{c^2+n^2+V_n}$$
 	  and 
 	 the perturbation $R$ satisfies $\|R\|_{0} \le C\epsilon$.
	Then, given any \(\omega= (\omega_{n})_{n \in \mathbb{Z}}  \in \Pi(c)\) satisfying the non-resonant conditions \eqref{NR_KG1},
	there exists a small $\epsilon_*:=\epsilon_*(\sigma,p,r)>0$.
	For any $0<\epsilon<\epsilon_*$, there exist \(V^* \in \mathcal{V}\) and real analytic symplectic coordinate transformation
	$$\Phi:D_*\times\{V^*\}\rightarrow D_0\times\{T\}, $$
	where
		\begin{equation}\label{V^*}
		\left\| V^*-T\right\| _{\infty}
		\leq\sum_{s=0}^{\infty}2\epsilon_{s}^{0.5}
		<\epsilon_{0}^{0.4}
	\end{equation}
	and
	\begin{equation}\label{260514-3}
		D_*=\left\{(z_{n})_{{n}\in\mathbb{Z}}:\frac{2}{3}\leq
		\left| z_{n}\right| \lfloor n \rfloor^p e^{r\ln^{\sigma}\lfloor n \rfloor}\leq\frac{5}{6}\right\}.
	\end{equation}
Furthermore, the transformation $\Phi$ satisfies 
	\begin{align}\label{symid}
	\left\| \Phi-\mbox{id}\right\|_{p},\left\| D\Phi-\mbox{Id}\right\|_{p
	\rightarrow p}\leq \epsilon_{0}^{0.4},
	\end{align}
rendering the Hamiltonian  $H$ in the form
	\begin{equation*}
		 H_*=H\circ\Phi=N_*+R_{*}^2,
	\end{equation*}
	 where
	\begin{align}\label{cvg1}
		N_*=\sum_{{n}\in\mathbb{Z}}
		\omega_{n} |z_{n}|^2
		\quad&\text{and}\quad
		\left\| R^2_{*}\right\|^c_{4\rho_0}\leq\frac{7}{6}\epsilon_0.
	\end{align}
Moreover, by choosing $I_n(0) \in \mathcal{I}_0$, we establish the existence of invariant tori $\mathcal{T}$ with frequency $\omega$ for the Hamiltonian vector field $X_{H_*}$.
	\end{theorem}

\begin{proof}
	Firstly, applying  Lemma \ref{IL} with $s=0$, set
	\begin{equation*}
		V_0^*=T,\hspace{12pt}\widetilde{V}_0=\mbox{id},
		\hspace{12pt}\epsilon_0=C\epsilon,
	\end{equation*}
	then the conditions (\ref{199})-(\ref{202'}) with $s=0$ are satisfied. 
	Next,  by Lemma \ref{IL} again, we obtain a sequence of domains $$D_{s}\times\mathcal{C}_{\eta_{s}}\left(V_{s}^*\right) \subset\cdots\subset
	D_{0}\times\mathcal{C}_{\eta_{0}}\left(V_{0}^*\right)$$
	 and a sequence of transformations
	\begin{equation*}
		\Phi^s=\Phi_1\circ\cdots\circ\Phi_s:
		D_{s}\times\mathcal{C}_{\eta_{s}}\left(V_{s}^*\right)\rightarrow D_{0}\times
		\mathcal{C}_{\eta_{0}}(V^{*}_{0}),
	\end{equation*}
	such that $H\circ\Phi^s=N_s+R_s$ for $s\geq1$. 
	Since the coordinate transformations
	satisfy (\ref{203})-\eqref{204} 
	and the frequency shift meet (\ref{206})-(\ref{IL2}). 
	we conclude from the iteration parameters in Section \ref{031501} that $V_s^*$ converges to a limit $V^*$ satisfying \eqref{V^*}
	and $\Phi^s$ converges uniformly 
	to $\Phi$  satisfying \eqref{symid}.
	Finally, from \eqref{200'}-\eqref{202'} in Lemma \ref{IL}, one gets  \eqref{cvg1}.  
\end{proof}

		\section{Persistence of the Structural Decomposition}\label{sec 4_5}
		
		This section establishes the inductive step for our iterative scheme. We aim to show that the structural decompositions of the NLKG and NLS Hamiltonians persist under symplectic transformations. 
		Following the approach in Section 3.1 of \cite{bambusi2025NLKG}, we verify this structural persistence, which is essential for the convergence analysis.
	
		 We begin with a preliminary claim on gauge invariance:
			\begin{claim}[\cite{bambusi2025NLKG}, Lemma 3.6]\label{claim1}
			Let $G$ be a smooth Gauge invariant Hamiltonian, denote by $\phi_G^t$ the flow generated by its Hamiltonian vector field. We assume $\phi_G^t$ to be defined for $t \in I$, for some interval $I \subseteq \mathbb{R}$, on some set $\mathcal{U}$. Then\\
				\textbf{(1)} $U_{-\alpha} \phi_G^t (U_\alpha z) = \phi_G^t (z)$ for all $z \in \mathcal{U}, t \in I$.\\
				\textbf{(2)} If $F$ is a Gauge invariant Hamiltonian, then, for any fixed $t \in I$, $F \circ \phi_G^t$ is a Gauge invariant Hamiltonian.\\
				\textbf{(3)} If $F$ is a non-Gauge invariant Hamiltonian, then, for any fixed $t \in I$, $F \circ \phi_G^t$ is a non-Gauge invariant Hamiltonian.
		\end{claim}
		
		\subsection{Structural decomposition of perturbations}\label{sec 3.2}

		We begin by considering the NLKG and NLS Hamiltonians
		\begin{equation*} 
			H = N + R \qquad \text{and} \qquad H^{\text{NLS}} = N^{\text{NLS}} + R^{\text{NLS}},
		\end{equation*}
		where
		\begin{align}\label{231130-1}
			N=\sum_{{n}\in\mathbb{Z}}
			\widetilde\lambda_{n} \,|z_{n}|^2
			\quad\text{and}\quad
			N^{\text{NLS}}=\sum_{{n}\in\mathbb{Z}}
			\widetilde\lambda^{\text{NLS}}_{n} \,|z_{n}|^2,
		\end{align}
		 $R$ and $R^{\text{NLS}}$ are small perturbations. Applying the method from Section \ref{sec 3.1}, we can construct two transformations $\Phi$ and $\Phi^{\text{NLS}}$ such that
		\begin{equation}\label{260503-3}
			H \circ \Phi = N_{+} + R_{+} 
			\qquad\text{and}\qquad
			 H^{\text{NLS}} \circ \Phi^{\text{NLS}} = N_{+}^{\text{NLS}} + R_{+}^{\text{NLS}},
		\end{equation}
		where $N_{+}$ and $R_{+}$ are given by \eqref{051402} and \eqref{051403}, and $N_{+}^{\text{NLS}}$ and $R_{+}^{\text{NLS}}$ are defined analogously.

		We next show that the transformation $\Phi$ in \eqref{260503-3} preserves the structural decomposition of the perturbation. Specifically, assuming $R$ admits the decomposition defined in \eqref{260421-1}:
		\begin{equation}\label{260503-2}
			R = R^{\text{NLS}} + R^{\text{NG}} + R^h,
		\end{equation}
		we aim to prove that the transformed perturbation $R_+$ maintains this structure.

		To be more precise, 
		let the components $R^{\text{NLS}}$, $R^{NG}$ and $R^h$ admit an expansion analogous to that of $R$ in \eqref{pertur}. 
		Based on \eqref{051304} and \eqref{260503-2}, we define
		\begin{equation}\label{F_app}
			F_{abkk'}^{app}:=\frac{R_{abkk'}^{\text{NLS}}}
			{\mathrm{i}\sum_{{n}\in\mathbb{Z}}
				\left(k_{n}-k^{'}_{n}\right)
				\widetilde \lambda_{n} },
		\end{equation}
		\begin{equation}\label{F_NG}
			F_{abkk'}^{{NG}}=\frac{R_{abkk'}^{{NG}}}
			{\mathrm{i}\sum_{{n}\in\mathbb{Z}}
				\left(k_{n}-k^{'}_{n}\right)
				\widetilde \lambda_{n} },
		\end{equation}
		\begin{equation}\label{F_h1}
			F_{abkk'}^{h1}:=\frac{R_{abkk'}^{h}}
			{\mathrm{i}\sum_{{n}\in\mathbb{Z}}
				\left(k_{n}-k^{'}_{n}\right)
				\widetilde \lambda_{n} }
		\end{equation}
		for any $a,b,k,k' \in \mathbb{N}^{\mathbb{Z}}$ with $|b|\le 1$ and
		$k\neq k'$.
		Furthermore, we define 
		\begin{equation}\label{F^h}
			F^{h}:=F^{{NG}}+ F^{app} + F^{h1}- F^{\text{NLS}},
		\end{equation}
		where the coefficients of $F^{\text{NLS}}$ are given by
		\begin{equation*}
			F_{abkk'}^{\text{NLS}}:=\frac{R_{abkk'}^{\text{NLS}}}
			{\sum_{{n}\in\mathbb{Z}}
				\left(k_{n}-k^{'}_{n}\right)
				\widetilde \lambda^{\text{NLS}}_{n} }.
		\end{equation*}
		Thus, we have $F= F^{\text{NLS}}+F^{h}$, where
		the coefficients of $F$ are given by \eqref{051304}.

		Next, we proceed to prove that the new perturbation given by \eqref{051403} admits the structure
		\begin{equation}
			R_+=R^{\text{NLS}}_++R^{{NG}}_++R^{h}_+.
		\end{equation}
		
		Analogous to  \eqref{051403}, we define $R_+^{\text{NLS}}$ as
		\begin{equation}\label{051403-4}
			R^{NLS}_+:=
			\sum_{n\ge 1} \left( \frac{1}{(n+1)!} \left( \widetilde{R^{NLS,\le 1}}\right)^{(n\star)}+ \frac{1}{n!}
			\left( R^{NLS,\le 1}\right)^{(n\star)}\right) 
			+\sum_{n\ge 0}\frac{1}{n!}\left( R^{NLS,2}\right)^{(n\star)},
		\end{equation}
		where 
		\begin{equation*}
			R^{(0\star)}:=R
			\quad\text{and}\quad
			R^{(n\star)}:=\left\{R^{(n-1)}, F^{\text{NLS}} \right\},\quad n\ge 1.
		\end{equation*}
		
		Meanwhile, in accordance with Definition \ref{def_NG} and Claim \ref{claim1}, $R_+^{\text{NG}}$ is defined by
		\begin{equation}\label{051403-2}
			R^{NG}_+:=
			\sum_{n\ge 1} \left( \frac{1}{(n+1)!} \left( \widetilde{R^{NG,\le 1}}\right)^{(n\star)}+ \frac{1}{n!}
			\left( R^{NG,\le 1}\right)^{(n\star)}\right) 
			+\sum_{n\ge 0}\frac{1}{n!}\left( R^{NG,2}\right)^{(n\star)}.
		\end{equation}
		It is worth noting, based on the property $(3)$ in Claim \ref{claim1}, that $R^{NG}_{+}$ is non-Gauge invariant.

		To close the decomposition, we define the remainder term as
		\begin{equation}\label{260503-4}
			R^{h}_+ :=R_+-R^{\text{NLS}}_+-R^{NG}_+.
		\end{equation}
		By observing that
		\begin{equation}\label{R^h+}
			{R}^{(n)}-\left({R^{NLS}}\right)^{(n\star)} -
			\left( {R^{NG}}\right)^{(n\star)}= \left( {R^{h}}\right)^{(n\star)} + {R}^{(n\divideontimes)},
		\end{equation}
		where
		\begin{equation*}
			R^{(0\divideontimes)}:=R
			\quad\text{and}\quad
			R^{(n\divideontimes)}:=\left\{R^{(n-1\divideontimes)}, F^{h} \right\},\quad n\ge 1.
		\end{equation*}
		Thus, it follows from \eqref{260503-4} and \eqref{R^h+} that 
		\begin{align}\label{051403-3}
			R^{h}_+&:=\sum_{n\ge 1} \frac{1}{(n+1)!}
			\left( \left( \widetilde{R^{h,\le 1}}\right)^{(n\star)} +\left( \widetilde{R^{\le 1}} \right)^{(n\divideontimes)}\right) 
			+\frac{1}{n!}	\left( \left( R^{h,\le 1}\right)^{(n\star)} +\left( R^{\le 1} \right)^{(n\divideontimes)}\right) 
			\nonumber\\
			& \quad\;
			+\sum_{n\ge 0}\frac{1}{n!}
			\left( \left( R^{h,2}\right)^{(n\star)}
			+\left( R^{2} \right)^{(n\divideontimes)}\right).
		\end{align}

		\subsection{Persistence of the decomposition and inductive estimates}\label{03150111}
		
To establish the convergence of the iterative scheme for $H_s$ and $H_s^{\text{NLS}}$, we first examine the persistence of their structural decompositions under symplectic transformations.
 Building upon the framework of Lemma \ref{sec 3.1}, we  have the following lemma.

		\begin{lemma}\label{Decomp}(\textbf{Iterative persistence of the decomposition})
			Consider the two Hamiltonians
			\begin{equation*}
			H_{s}=N_{s}+R_{s}	
			\qquad\text{and}\qquad
			H^{\text{NLS}}_{s}=N^{\text{NLS}}_{s}+R^{\text{NLS}}_{s},
			\end{equation*}
			which are real analytic on $D_{s}\times\mathcal{C}_{\lambda_s}\left( V_s^*\right)$
			and $D_{s}\times{\mathcal{C}_{\lambda_s}\left( V_s'\right)}$, respectively.
		Under the same assumptions as in Lemma \ref{IL},
		we further assume that the perturbation $R_s$ admits the structural decomposition
			\begin{equation}\label{260509-1}
				R_{s}=R_{s}^{\text{NLS}}+R_{s}^{{NG}}+R_{s}^{h}
			\end{equation}
		and its components satisfy$$R^{\text{NLS}}_{s} = R^{\text{NLS},0}_{s} + R^{\text{NLS},1}_{s} + R^{\text{NLS},2}_{s}$$
		with similar degree-based decompositions holding for $R^{\text{NG}}_{s}$ and $R^{h}_{s}$. 
Suppose these components are small perturbations satisfying the following bounds:
\begin{align}
	\label{200}
	\|R^{\text{NLS},0}_{s}\|^{\infty}_{\rho_{s}}, \quad 
	\|R^{\text{NG},0}_{s}\|^{\infty}_{\rho_{s}}, \quad 
	\|R^{h,0}_{s}\|^{*}_{\rho_{s}} &\leq \epsilon_{s}, \\
	\label{201}
	\|R^{\text{NLS},1}_{s}\|^{\infty}_{\rho_{s}}, \quad 
	\|R^{\text{NG},1}_{s}\|^{\infty}_{\rho_{s}}, \quad 
    \|R^{h,1}_{s}\|^{*}_{\rho_{s}} &\leq \epsilon_{s}^{0.6}, \\
	\label{202}
	\|R^{\text{NLS},2}_{s}\|^{\infty}_{\rho_{s}}, \quad 
	\|R^{\text{NG},2}_{s}\|^{\infty}_{\rho_{s}}, \quad 
    \|R^{h,2}_{s}\|^{*}_{\rho_{s}} &\leq (1+d_s)\epsilon_0.
\end{align}
			
		Then, for any $V \in \mathcal{C}_{\eta_{s}}(V_{s}^*)$ and $V^{\text{NLS}} \in \mathcal{C}_{\eta_{s}}(V_{s}')$, there exist two real analytic symplectic coordinate transformations $\Phi_{s}, \Phi^{\text{NLS}}_{s} : D_{s+1,p} \to D_{s,p}$
		such that
			$$H_{s+1} = H_{s} \circ \Phi_{s} = N_{s+1} + R_{s+1},$$$$H^{\text{NLS}}_{s+1} = H^{\text{NLS}}_{s} \circ \Phi^{\text{NLS}}_{s} = N^{\text{NLS}}_{s+1} + R^{\text{NLS}}_{s+1}.$$
		Moreover, the difference $\Psi_{s}^h := \Phi_{s} - \Phi^{\text{NLS}}_{s} : D_{s+1,p} \to D_{s,p-4}$ satisfies 
			\begin{equation}\label{trans_h}
				\sup_{z\in D_{s+1,p}} 
				\left\|\Psi_s^h\right\|_{p-4}\le h \epsilon_{s}^{0.5}.
			\end{equation}
		In addition, the perturbation $R_{s+1}$ preserves the  decomposition \eqref{260509-1}, and the estimates \eqref{200}-\eqref{202} hold for $s+1$.
		\end{lemma}
		
		\begin{proof}
			The proof can be divided into the following two steps.
		
			
			\textbf{Step 1. Estimates on the difference map $\Psi_s^h$.}

		Building on the lower bound estimate for \eqref{NR_KG1} established in Step 1 of Lemma \ref{IL}, it remains to derive analogous lower bounds for the non-resonance conditions \eqref{NR_KG2} and \eqref{NR_NLS}.

		Following the proof of \eqref{M7} and utilizing the non-resonance conditions \eqref{NR_KG2}, we deduce that for any $\ell \in \mathbb{Z}^{\mathbb{Z}}$ with $0 < |\ell| < \infty$ and $\alpha(\ell) \neq 0$,
		\begin{equation}\label{M8}
			\left| \sum_{n \in \mathbb{Z}} \ell_{n} \cdot \widetilde{\lambda}_{s,n}(V_s^*) \right|
			= \left| \sum_{n \in \mathbb{Z}} \ell_{n} \cdot \omega_{n} \right|
			\geq \gamma c^{2} \epsilon_s^{0.01}.
		\end{equation}
		Similarly, by applying the same proof alongside the condition \eqref{NR_NLS}, we obtain that for any $\ell \in \mathbb{Z}^{\mathbb{Z}}$ with $0 < |\ell| < \infty$,
		\begin{equation}\label{M9}
			\left| \sum_{n \in \mathbb{Z}} \ell_{n} \cdot \widetilde{\lambda}^{\text{NLS}}_{s,n}(V_s') \right|
			= \left| \sum_{n \in \mathbb{Z}} \ell_{n} \cdot \omega^{\text{NLS}}_{n} \right|
			\geq \gamma \epsilon_s^{0.01}.
		\end{equation}


			Next, we proceed to estimate the coordinate transformation \eqref{trans_h}.
			
			Recalling \eqref{F^h}, one has $F_{s}=F^{\text{NLS}}_{s}+F^{h}_{s}$, where
			\begin{equation}\label{260510-2}
				F^{h}_{s}=F^{NG}_{s}+F^{h1}_{s}
				+\left( F^{app}_{s}-F^{\text{NLS}}_{s}\right),
			\end{equation}
			where the coefficients of \( F_s^{app} \), \( F_s^{NG} \) and \( F_s^{h1} \) are defined in \eqref{F_app}-\eqref{F_h1}.
			On one hand, following the proof of \eqref{022501}, we can deduce from \eqref{M7}, \eqref{200}, \eqref{201} and \eqref{M8}  
			that
			\begin{equation}\label{F^NG+h1}
				\left\|F^{NG}_{s}\right\|^{*}_{\rho_s}, \; \left\|F^{h1}_{s}\right\|^{*}_{\rho_s} \leq  \epsilon_s^{0.54}.	
			\end{equation}
			On the other hand, we calculate the coefficients of the term \( F_s^{app} - F_s^{\text{NLS}} \) as
			\begin{align}\label{260510-1}
				\left|F_{abkk'}^{app} - F^{\text{NLS}}_{abkk'}\right| &= \left|\frac{1}{\mathrm{i}\sum_{{n}\in\mathbb{Z}}
					\left(k_{n}-k^{'}_{n}\right)
					\widetilde \lambda_{n} } - \frac{1}{\mathrm{i}\sum_{{n}\in\mathbb{Z}}
					\left(k_{n}-k^{'}_{n}\right)
					\widetilde \lambda^{\text{NLS}}_{n} }\right| \left|R_{abkk'}^{\text{NLS}} \right|\nonumber\\
				&\le \frac{ \left({n_1^{*} }/{ \sqrt{c} }\right)^4 \left(\sum_{{n}\in\mathbb{Z}}
					k_{n}+k^{'}_{n}\right)}{\left| \sum_{{n}\in\mathbb{Z}}
					\left(k_{n}-k^{'}_{n}\right)
					\widetilde \lambda_{n} \right|
					\cdot \left| \sum_{{n}\in\mathbb{Z}}
					\left(k_{n}-k^{'}_{n}\right)
					\widetilde \lambda^{\text{NLS}}_{n} \right| }   \left|R_{abkk'}^{\text{NLS}} \right|,
			\end{align}
			where the last inequality follows from the fact that
			$$ \widetilde \lambda_{s,n} -	\widetilde \lambda^{\text{NLS}}_{s,n}  - c^2
			 \le \frac{n^4 }{ c^2+n^2 }
			\le \left( \frac{n_1^{*} }{ \sqrt{c} }\right)^4  $$
			with $n_1^*=n_1^*(k,k')$.
			Therefore, 
			following the proof of \eqref{022501} again, we can deduce from \eqref{M7}, \eqref{200}, \eqref{201} and \eqref{M9} that
			\begin{equation}\label{F^app}
				\left\| F_s^{app} - F^{\text{NLS}}_s \right\|^{*}_{\rho_s'}
				\leq   {\gamma}^{-2} \epsilon_s^{-0.02}	\|R^{\text{NLS}}_{s}\|^{\infty}_{\rho_{s}}
				\leq \epsilon_s^{0.54},
			\end{equation}
			where $\rho_s'=\rho_s+(\rho_{s+1}-\rho_s)/2$.
			It should be noted that the  factor $
			(\sum_{n} k_n + k'_n)$ in \eqref{260510-1} is absorbed by the increased regularity of the norm $\|\cdot\|_{\rho'_s}^*$.

			By utilizing \eqref{260510-2}, \eqref{F^NG+h1} and \eqref{F^app}, we obtain
			\begin{equation}\label{est_X_Fh}
				\left\| F_s^{h} \right\|^{*}_{\rho_s'}
				\leq \left\|F^{NG}_{s}\right\|^{*}_{\rho_s}+  \left\|F^{h1}_{s}\right\|^{*}_{\rho_s} +
				\left\| F_s^{app} - F^{\text{NLS}}_s \right\|^{*}_{\rho_s'}
				\le  \epsilon_s^{0.53}.
			\end{equation}
			Moreover, invoking  \eqref{vector_2'} in Lemma \ref{vector} yields
			\begin{equation}\label{X_app}
				\sup_{z\in D_{s,p}}
				\left\| X_{F^{h}_{s}}\right\|_{p-4}
				\leq h\cdot	C'(\rho_s) 
				\left\| F_s^{h} \right\|^{*}_{\rho_s'}
				\le h \epsilon_s^{0.52}.
			\end{equation}
			%
			%
		Applying the estimate \eqref{260418-2} in Lemma \ref{claim3}, we obtain $\Psi_{s+1}^h : D_{s+1,p} \to D_{s,p-4}$ such that
			\begin{equation}\label{trans_h-1}
				\sup_{z\in D_{s+1,p}}\left\|\Psi_s^h\right\|_{p-4}
				\leq\sup_{z\in D_{s,p}}	\left\|X_{F_s}-X_{F^{\text{NLS}}_s} \right\|_{p-4}
				=\sup_{z\in D_{s,p}}
				\left\|X_{ F_s^{h}} \right\|_{p-4}
				\le h\epsilon_{s}^{0.5}.
			\end{equation}
			Thus, the estimate \eqref{trans_h} is established.

			\textbf{Step 2. Estimate on the new perturbations.} 
			
%
		 
		First, we provide the explicit expressions for $R^{\text{NLS}}_{s+1}$, $R^{NG}_{s+1}$ and $R^{h}_{s+1}$.
		
	$\bullet$	In analogy to \eqref{051403}, we define $R^{\text{NLS}}_{s+1}$ as
			$$R^{\text{NLS}}_{s+1}=\sum_{n\ge 1}\frac{1}{(n+1)!} 
			\left( \widetilde{R_s^{NLS, \le 1}} \right)^{(n\star)}
			+\frac{1}{n!} \left( R_s^{\text{NLS},\le1}\right)^{(n\star)}+\sum_{n\ge 0}\frac{1}{n!}\left( R_s^{\text{NLS},2}\right)^{(n\star)}.
			$$

			$\bullet$ From \eqref{051403-2}, we define $R^{NG}_{s+1}$ as
			\begin{equation*}
				R^{NG}_{s+1} =
				\sum_{n\ge 1}\frac{1}{(n+1)!} 
				\left( \widetilde{R_s^{NG,\le 1}} \right)^{(n\star)}+ \frac{1}{n!}
				\left( R_s^{NG, \le 1}\right)^{(n\star)}
				+\sum_{n\ge 0}\frac{1}{n!}\left( R^{NG,2}\right)^{(n\star)}.
			\end{equation*}

			$\bullet$ From \eqref{051403-3}, we define $R^{h}_{s+1}$ as
			\begin{align}
				R^{h}_{s+1}&:=\sum_{n\ge 1} \frac{1}{(n+1)!}
				\left( \left( \widetilde{R_s^{h,\le 1}}\right)^{(n\star)} +\left( \widetilde{R_s^{\le 1}} \right)^{(n\divideontimes)}\right) 
				+\frac{1}{n!}	\left( \left( R_s^{h,\le 1}\right)^{(n\star)} +\left( R_s^{\le 1} \right)^{(n\divideontimes)}\right) 
				\nonumber\\
				& \quad\;
				+\sum_{n\ge 0}\frac{1}{n!}
				\left( \left( R_s^{h,2}\right)^{(n\star)}
				+\left( R_s^{2} \right)^{(n\divideontimes)}\right). 
			\end{align}

		Next, we verify that the estimates \eqref{200}-\eqref{202} hold for $s+1$. 

		By setting $F^{\text{NLS}}_s = F^{\text{NLS}, 0}_s + F^{\text{NLS}, 1}_s$ and $F^{h}_s = F^{h, 0}_s + F^{h, 1}_s$,
		 it follows from \eqref{200}, \eqref{201}, \eqref{M9} and the proof of \eqref{est_X_Fh}
		 that the following uniform estimates
		\begin{equation}\label{260511-1}
				\| F^{\text{NLS},0}_s \|^{\infty}_{\rho_s} \leq \epsilon_s^{0.95},\quad
			\| F^{\text{NLS},1}_s \|^{\infty}_{\rho_s} \leq \epsilon_s^{0.55},\quad
			\| F^{\text{NLS}}_s \|^{\infty}_{\rho_s} \leq \epsilon_s^{0.54}
		\end{equation}
		and
		\begin{equation}\label{260511-2}
			\| F^{h,0}_s \|^*_{\rho_s'} \leq \epsilon_s^{0.94},\quad
			\| F^{h,1}_s \|^*_{\rho_s'} \leq \epsilon_s^{0.54},\quad
			\| F^{h}_s \|^*_{\rho_s'} \leq \epsilon_s^{0.53}.
		\end{equation}

	Regarding the remainder terms $R^\sigma_{s+1}$ with $\sigma \in \{\text{NLS,  } NG,h\}$, we verify that the inductive assumptions \eqref{200}-\eqref{202} are satisfied for each component $R^{\sigma, i}_{s+1}$ ($i=0,1,2$) by following Step 2 in Lemma \ref{IL}:
		
		$\bullet$ For $R^{\text{NLS}}_{s+1}$ and $R^{NG}_{s+1}$: By combining \eqref{flow1} (with $c=\infty$) and \eqref{260511-1}, it follows that the estimates \eqref{200}-\eqref{202} hold for the components
		 $R^{NG,i}_{s+1}$ and $R^{\text{NLS},i}_{s+1}$.

		$\bullet$ For $R^{h}_{s+1}$: By invoking \eqref{flow4} along with \eqref{260511-1} and \eqref{260511-2} respectively, we conclude that the assumptions \eqref{200}-\eqref{202} are also satisfied for the components $R^{h, i}_{s+1}$.

		To sum up, we complete the proof of the iterative lemma.
			\end{proof}

	\section{Proof of the main theorem}\label{sec 4}

	\subsection{Proof of Theorem \ref{th1} }
	Consider the NLKG Hamiltonian (\ref{Ham}), then from  \eqref{mome1} one has 
	\begin{equation*}
		\|R\|^c_{\rho_{0}}\leq \|R\|^c_{0} 
	   	\le C\epsilon:= \epsilon_0
	\end{equation*}
	with constant $C>0$, which implies the assumptions (\ref{200})-(\ref{202}) with $s=0$ in Lemma \ref{IL} are satisfied.
	In conclusion, by applying Theorem \ref{KAMth}, we have that
	$ \mathcal{E}=\Phi(\mathcal{T}) $
	is the desired invariant tori for the Hamiltonian (\ref{Ham}) of equation \eqref{NLKG}.
	 Moreover, we deduce that the tori \( \mathcal{E} \) 
	is linearly stable. This conclusion follows from the fact that the Hamiltonian \(N_*\) 
	given by \eqref{cvg1} is normal forms of order $2$ around the invariant tori.

 Regarding the NLS Hamiltonian \eqref{Ham_NLS1} of equation \eqref{NLS}, it follows from \eqref{mome_NLS} that 
 $$\|R^{\text{NLS}}\|^{\infty}_{\rho_{0}} \leq \|R^{\text{NLS}}\|^{\infty}_{0} \le \epsilon_0.$$ Thus, by applying a similar KAM machinery, one can construct the linearly stable invariant torus $\mathcal{E}^{\text{NLS}} = \Phi^{\text{NLS}}(\mathcal{T}^{\text{NLS}})$. This result for the NLS equation  has been established in \cite{Cong2023TheEO}.

	\subsection{Proof of Theorem \ref{th2}}
\begin{proof}
	\textbf{Step 1. Global approximation of solutions to equations \eqref{NLKG} and \eqref{NLS}.}

	In this step, we aim to establish the convergence estimate \eqref{solu_conver}. Once \eqref{solu_conver} is obtained, the asymptotic behavior \eqref{solu_conver1} follows immediately by considering the sequence $\{c_m\}_{m\ge 1}$ defined in \eqref{seq_cn} and passing to the limit as $m \to \infty$.
	
 According to  Theorem \ref{KAMth}, there exists a near-identity transformation $\Phi$ such that
\begin{equation}\label{2605161113}
H_* = H \circ \Phi = N_* + R_*^2, \quad N_* = \sum_{n \in \mathbb{Z}} \left( \omega'_n+{c^2}\right)  I_n,
\end{equation}
where $\omega'= \omega-c^2 \in\Pi'$ and the perturbation  $R_*$ is quadratic in the action $J$. It is easy to see that for a given initial action { $I(0) \in \mathcal{I}_0 $} and the frequency $\omega'\in\Pi'$, the explicit solution of the Hamiltonian  (\ref{2605161113}) can be given by $z^*(t)=(z^*_n(t))_{n \in \mathbb{Z}}$ with
\begin{equation*} 
	z^*_n(t) = \sqrt{I_n(0)} e^{\mathrm{i}\left( \omega'_n+{c^2}\right)  t},
\end{equation*}
which is the corresponding linear flow on the KAM torus.
Consequently, $$\psi^{ {KG}}_{c,\omega'}(t):=\Psi^{{KG}}((\omega'+c^2) t, I(0)) =\Phi(z^*(t), \bar{z}^*(t))$$
 represents the solution to the NLKG equation \eqref{NLKG} in the original coordinates.  Then
 {$$ e^{-\mathrm{i}c^2 t}  \psi^{{KG}}_{c,\omega'}(t)=\Psi^{{KG}}(\omega' t, I(0)) \label{solu1}=\Phi(z'(t), \bar{z}'(t)),$$
 	where $z'(t)=(z'_n(t))_{n \in \mathbb{Z}}$ with
 	\begin{equation}\label{260514-1}
 		z'_n(t) = \sqrt{I_n(0)} e^{\mathrm{i}
 			\omega'_n  t}.
 \end{equation}}

Similarly, we have analogous conclusions for the NLS equation \eqref{NLS}, i.e., for the same $I(0) \in\mathcal{I}_0$ and frequency $\omega\in\Pi^{\text{NLS}}$, the explicit solution corresponding to the Hamiltonian $H_*^{\text{NLS}}$ is $z(t)=(z_n(t))_{n \in \mathbb{Z}}$ with
\begin{equation}\label{260514-2}
	z_n(t) = \sqrt{I_n(0)} e^{\mathrm{i} \omega_n t}, 
\end{equation} 
and thus 
$$\varphi^{\text{NLS}}_{\omega}(t) :=\Psi^{\text{NLS}}(\omega t, I(0))=\Phi^{\text{NLS}}(z(t), \bar{z}(t))$$
represents the solution to  the NLS equation \eqref{NLS} in the original coordinates.


%
We take
\begin{equation}\label{260513-4}
	\mathcal{N}\equiv \mathcal{N}(c) =\left[\gamma^{1/15} c^{2/5} \right]  + 1.
\end{equation}
Then if the frequency $	{\omega}'$ and $ {\omega}$
satisfy
\begin{equation*} 
	\widetilde{\omega}' = \widetilde{\omega},
\end{equation*} 
i.e.,
\begin{equation}\label{same_omega}
	\omega_n'=\omega_n\quad \text{when}\quad |n|\le \mathcal{N},
\end{equation}
%
we claim that the following inequality holds
\begin{align}\label{est_p-4_1}
	&\left\|\Phi(z', \bar{z}') -\Phi^{\text{NLS}}(z, \bar{z})\right\|_{p-4} \leq 10 h^{0.7}.
\end{align}

	In fact, by the triangle inequality we have
	\begin{align*}
	&\quad\,\left\|\Phi(z', \bar{z}') -\Phi^{\text{NLS}}(z, \bar{z})\right\|_{p-4} \nonumber\\
	&\le \left\|\Phi(z', \bar{z}') -\Phi^{\text{NLS}}(z', \bar{z}')\right\|_{p-4} +\left\|\Phi^{\text{NLS}}(z', \bar{z}') -\Phi^{\text{NLS}}(z, \bar{z})\right\|_{p-4}.
\end{align*}

	\textbf{Estimate the first term.} 
	Based on the estimate \eqref{trans_h} in Lemma \ref{Decomp}, we have
	\begin{equation}\label{maps1}
			\left\|\Phi- \Phi^{\text{NLS}} \right\|_{p-4}
			\le\sum_{s=0}^{\infty} \left\| \Psi_s^h \right\|_{p-4} 
			\leq\sum_{s=0}^{\infty}h\epsilon_{s}^{0.5}
			<h\epsilon_{0}^{0.4}.
	\end{equation}
Consequently, in view of  \eqref{maps1}, it follows that
		\begin{align}\label{same_theta}
		\quad\, 	\left\|\Phi(z', \bar{z}') -\Phi^{\text{NLS}}(z', \bar{z}')\right\|_{p-4}
		\le 
		\left\|\Phi- \Phi^{\text{NLS}} \right\|_{p-4}
		\le  h \epsilon_{0}^{0.4}.
	\end{align}
%

	\textbf{Estimate the second term.} In view of \eqref{symid}, we obtain 
	\begin{align}\label{same_trans}
		\quad\, 	\left\|\Phi^{\text{NLS}}(z', \bar{z}') -\Phi^{\text{NLS}}(z, \bar{z})\right\|_{p-4}\le
		\left\|D\Phi^{\text{NLS}} \right\|_{p-4}\cdot
		\left\|z'- z \right\|_{p-4}
		\le  4	\left\|z'- z \right\|_{p-4}.
	\end{align}

	Now we estimate the remaining factor $	\left\|z'- z \right\|_{p-4}$. Note that from \eqref{same_omega} that
	\begin{align*}
		\left|z'_n-z_n\right|=\left|\sqrt{I_n(0)} \left(e^{\mathrm{i}\omega'_n t}-e^{\mathrm{i} \omega_n t}\right)\right|=0\quad\text{for}\quad |n|\le \mathcal{N}.
	\end{align*}
Hence in view of (\ref{042501}), (\ref{260514-1}) and (\ref{260514-2}), we obtain
	\begin{align}\label{260512-1}
	  \left\|z'- z \right\|_{p-4}
		&= \sup_{|n|\ge \mathcal{N}+1}
		\left| \sqrt{I_n(0)} \left( e^{\mathrm{i} \omega'_n t }-
		e^{\mathrm{i}\omega_n t} \right) \right| \lfloor n \rfloor^{p-4} e^{r\ln^{\sigma}\lfloor n \rfloor}\nonumber\\
		&\le 2\mathcal{N}^{-4} \sup_{|n|\ge \mathcal{N}+1}
		\left| \sqrt{I_n(0)}\right|  \,  \lfloor n \rfloor^{p} e^{r\ln^{\sigma}\lfloor n \rfloor}
		\nonumber\\
		&\le  2  \gamma^{-4/15} c^{-8/5} \left\|  \sqrt{I(0)}\right\|_{p}
		\le 2 h^{0.75},
	\end{align}
where the last inequality holds since $h = c^{-2}$, $c \geq \epsilon^{-0.5}$, and $\gamma$ depends on $\epsilon$.
	\begin{remark}
		We note that while setting $\mathcal{N} \sim c^{1/2}$ would yield the optimal decay rate in \eqref{260512-1}, it would lead to a larger resonance set. The lower scaling $\mathcal{N} \sim c^{2/5}$ chosen in  \eqref{260513-4} is a deliberate trade-off to balance the high-frequency approximation error and the subsequent measure estimate \eqref{260601-1}.
	\end{remark}
Therefore,  in view of (\ref{same_trans})-(\ref{260512-1}), we obtain
\begin{equation}\label{same_tran1}
	 \left\|\Phi^{\text{NLS}}(z', \bar{z}') -\Phi^{\text{NLS}}(z, \bar{z})\right\|_{p-4} \le 8 h^{0.75}.
\end{equation}

Combining \eqref{same_theta} with \eqref{same_tran1} finally yields the desired estimate \eqref{est_p-4_1}.

Its proof proceeds along the same lines as that of 
 \begin{align}\label{est_p}
		\left\|\Phi(z', \bar{z}') -\Phi^{\text{NLS}}(z, \bar{z})\right\|_{p} \leq 10.
	\end{align} 
Specifically, from \eqref{symid}, noting that both $\Phi$ and $\Phi^{\text{NLS}}$ are near-identity maps satisfying the inductive estimates, we obtain
	$$\left\|\Phi - \Phi^{\text{NLS}} \right\|_{p} \le \left\| \Phi - \mathrm{id} \right\|_{p} + \left\| \Phi^{\text{NLS}} - \mathrm{id} \right\|_{p} \le 2\epsilon_{0}^{0.4}.$$
	 Following a similar decomposition and the operator bounds used in the derivation of \eqref{est_p-4_1}, the estimate \eqref{est_p} follows immediately.

	Finally, by invoking the interpolation inequality
\begin{equation*}
	\| u \|_{p-4a} \leq \| u \|_{p-4}^a \| u \|_{p}^{1-a}, \quad a \in [0,1],
\end{equation*}
and substituting the estimates \eqref{est_p-4_1} and \eqref{est_p} into it, we conclude the proof of  \eqref{solu_conver}.

 	\textbf{Step 2. Measure estimate of common frequencies.}

First, we establish the estimate \eqref{mes_comm}. 

Before proceeding with the proof, we first introduce some notation for the sets involved.
Let $	\widetilde{\Omega}_{\mathcal{N}}\equiv\widetilde{\Omega}_{\mathcal{N}}(c)$ denote the intersection of the low-frequency components
 \begin{equation}\label{260513-3}
 \widetilde{\Omega}_{\mathcal{N}}(c) := \widetilde{\Pi}_{\mathcal{N}}'(c) \cap \widetilde{\Pi}^{\text{NLS}}_{\mathcal{N}},
 \end{equation}
 where $\widetilde{\Pi}_{\mathcal{N}}'(c)$ is defined in \eqref{260513-1} 
 and $\widetilde{\Pi}_{\mathcal{N}}^{\text{NLS}}$ is defined in a manner analogous to \eqref{260513-1} based on \eqref{Pi}.
 Subsequently, the sets $\Omega'(c)$ and $\Omega^{\text{NLS}}$ are constructed as Cartesian products
 \begin{equation}\label{Omega'}
 	\Omega'(c) := \widetilde{\Omega}_{\mathcal{N}}(c)
 	\times \widehat{\Pi}_{\mathcal{N}}'(c)
 	\subset \Pi'(c)
 	 \end{equation}
 	and
 	 \begin{equation}\label{Omega''}
 	\Omega^{\text{NLS}}(c) := \widetilde{\Omega}_{\mathcal{N}}(c)
  \times \widehat{\Pi}^{\text{NLS}}_{\mathcal{N}}\subset{\Pi}^{\text{NLS}}.
 \end{equation}
 We now define the non-resonant sets by excising the exceptional sets $\mathcal{R}'(c) := \mathcal{R}(c) - c^2$ and $\mathcal{R}^{\text{NLS}}$, which yields 
  \begin{equation}\label{Xi'}
  	{\Xi}'(c) := \Omega'(c) \setminus 
  	\left( {\mathcal{R}}'(c) \cup \mathcal{R}^{\text{NLS}}\right) 
  \end{equation}
  and
  \begin{equation}\label{Xi^NLS}
  	\Xi^{\text{NLS}}(c) :=
  	\Omega^{\text{NLS}}(c) \setminus
  	\left({\mathcal{R}}'(c) \cup \mathcal{R}^{\text{NLS}}\right).
  \end{equation}
  For brevity, these sets will often be written without explicit dependence on $c$ in the sequel.
  Next, we claim that the following upper bounds hold:
  \begin{equation}\label{260525-1}
  	\text{meas}\left( {\Pi}' \setminus {\Xi}' \right) 	= O(\gamma^{1/3})
  \end{equation}
  and
   	\begin{equation}\label{mes_Pi^NLS}
  	 		\text{meas}\left( {\Pi}^{\text{NLS}} \setminus \Xi^{\text{NLS}} \right)  = O(\gamma^{1/3}).
  	 	\end{equation} 

We first focus on the estimate \eqref{260525-1}.
By virtue of the relation \eqref{Xi'}, a simple subadditivity argument yields
\begin{equation}
	\text{meas}(\Pi' \setminus \Xi') \le \text{meas}(\Pi' \setminus \Omega') + \text{meas}(\mathcal{R}')
	+ \text{meas}( \mathcal{R}^{\text{NLS}}).
\end{equation}
In Lemmas \ref{Le_meas3} and \ref{Le_meas6}, we shall prove that the sets $\mathcal{R}' \subset \Pi'$ and $\mathcal{R}^{\text{NLS}} \subset \Pi^{\text{NLS}}$ satisfy
\begin{equation}\label{260513-5}
	\text{meas}(\mathcal{R}'), \; \text{meas}(\mathcal{R}^{\text{NLS}}) = O(\gamma^{1/3}),
\end{equation}
where $\text{meas}(\mathcal{R}') = \text{meas}(\mathcal{R})$ holds by translation invariance, and ${\Pi}'$ and $\Pi^{\text{NLS}}$ are naturally endowed with a probability measure.
Consequently, it remains to verify that
 \begin{equation}\label{260525-3}
	\text{meas}\left( {\Pi}' \setminus \Omega' \right) = O(\gamma^{1/3}).
\end{equation}

 To this end, utilizing the product structure from \eqref{Omega'}, we decompose $\Xi'$ 
  into their respective low- and high-frequency components:
\begin{equation}
	\Xi' = \widetilde{\Xi}'_{\mathcal{N}} \times \widehat{\Xi}'_{\mathcal{N}} .
\end{equation}
Combined with \eqref{260513-3} and \eqref{Omega'}, we observe the decomposition
\begin{equation}\label{Pi'}
	{\Pi}' \setminus \Omega' = \left( \widetilde{\Pi}'_{\mathcal{N}} \setminus \widetilde{\Omega}_{\mathcal{N}} \right) \times \widehat{\Pi}'_{\mathcal{N}}
	= \left( \widetilde{\Pi}'_{\mathcal{N}} \setminus \widetilde{\Pi}_{\mathcal{N}}^{\text{NLS}} \right) \times \widehat{\Pi}_{\mathcal{N}}'.
\end{equation}
It follows from this product structure and the normalization of the measure on $\widehat{\Pi}'_{\mathcal{N}}$ that
\begin{eqnarray}\label{260530-1}
	\text{meas}(\Pi' \setminus \Omega') = \text{meas}\left(\widetilde{\Pi}'_{\mathcal{N}} \setminus \widetilde{\Pi}_{\mathcal{N}}^{\text{NLS}}\right).
\end{eqnarray}
Recall that
$$\widetilde{\Pi}_{\mathcal{N}}' = \prod_{|n|\le \mathcal{N}}  \left[ a_{n}^1, a_{n}^2\right] \quad \text{and} \quad \widetilde{\Pi}_{\mathcal{N}}^{\text{NLS}} = \prod_{|n|\le \mathcal{N}}  \left[ b_{n}^1, b_{n}^2\right]  ,$$
with the interval endpoints given by 
\begin{align*}
	a_n^1 &= c\sqrt{c^2+n^2}-c^2, \quad  a_n^2 = a_n^1 + \frac{c}{3\sqrt{c^2+n^2}}, \\
	b_n^1 &= \frac{1}{2}n^2, \quad b_n^2 = b_n^1 + \frac{1}{3}.
\end{align*}
To estimate the measure of the set difference, we consider the intersection $\widetilde{\Pi}'_{\mathcal{N}} \cap \widetilde{\Pi}_{\mathcal{N}}^{\text{NLS}}$. Since $a_n^1 \le b_n^1$ and $a_n^2 \le b_n^2$ for all $|n| \le \mathcal{N}$, the intersection forms a regular product of intervals
$$\widetilde{\Pi}'_{\mathcal{N}} \cap \widetilde{\Pi}_{\mathcal{N}}^{\text{NLS}} = \prod_{|n|\le \mathcal{N}} \left[ b_n^1, a_n^2\right].$$
Recalling that the measure on each component interval of $\displaystyle \widetilde{\Pi}'_{\mathcal{N}}$ is normalized by its length $\displaystyle \text{len}_n = a_n^2 - a_n^1$, the product measure satisfies
$$\text{meas}\left(\widetilde{\Pi}'_{\mathcal{N}} \cap \widetilde{\Pi}_{\mathcal{N}}^{\text{NLS}}\right) = \prod_{|n|\le \mathcal{N}} \frac{a_n^2 - b_n^1}{a_n^2 - a_n^1} = \prod_{|n|\le \mathcal{N}} \left( 1 - \frac{b_n^1 - a_n^1}{a_n^2 - a_n^1} \right).$$
Applying the Weierstrass product inequality 
$$\displaystyle \prod_{|n|\le \mathcal{N}}(1-x_n) \ge 1 - \sum_{|n|\le \mathcal{N}} x_n,$$
 we obtain
\begin{align*}
	\text{meas}\left(\widetilde{\Pi}'_{\mathcal{N}} \setminus \widetilde{\Pi}_{\mathcal{N}}^{\text{NLS}}\right) = 1 - \text{meas}\left(\widetilde{\Pi}'_{\mathcal{N}} \cap \widetilde{\Pi}_{\mathcal{N}}^{\text{NLS}}\right)
	\le \sum_{|n|\le \mathcal{N}} \frac{b_n^1 - a_n^1}{a_n^2 - a_n^1}.
\end{align*}
Directly evaluating this summation yields
\begin{align*}
		\text{meas}\left(\widetilde{\Pi}'_{\mathcal{N}} \setminus \widetilde{\Pi}_{\mathcal{N}}^{\text{NLS}}\right)
		\le \sum_{|n| \le \mathcal{N}} \frac{n^4}{2(c^2+n^2)} \cdot \frac{3\sqrt{c^2+n^2}}{c}
		\le 2c^{-2}(2\mathcal{N}+1)^5. 
\end{align*}
Benefiting from the choice of  $\mathcal{N} \sim c^{2/5}$ in \eqref{260513-4}, we obtain
\begin{equation}\label{260601-1}
		\text{meas}\left( \widetilde{\Pi}'_{\mathcal{N}} \setminus \widetilde{\Pi}^{\text{NLS}}_{\mathcal{N}} \right) 
		= O(\gamma^{1/3}).
\end{equation}
Combining this with \eqref{Pi'} and \eqref{260530-1} establishes the full-space estimate \eqref{260525-3}, thereby completing the estimate \eqref{260525-1}. 

Similarly, an analogous argument applied to the NLS case yields the estimate \eqref{mes_Pi^NLS}.

Consequently, the estimates \eqref{260525-1} and \eqref{mes_Pi^NLS}
 collectively yield \eqref{mes_comm} with $\gamma' = \gamma^{1/3}$.
%
%

We now proceed to establish  \eqref{mes_comm1}. 

Given an increasing sequence $\{c_m\}_{m\ge 1}$ satisfying \eqref{seq_cn}, we take 
$$\mathcal{N}_m \equiv \mathcal{N}(c_m) \to \infty\quad\text{as}\quad m \to \infty$$ according to \eqref{260513-4}.
Then
we characterize the limiting set $\Xi_\infty$ by passing to the limit as $m \to \infty$ in \eqref{Xi'}:
$$\begin{aligned}
	\lim_{m \to \infty} \Xi'(c_m) &= \lim_{m \to \infty} \left[ \left( \widetilde{\Omega}_{\mathcal{N}_m}(c_m) \times \widehat{\Pi}'_{\mathcal{N}_m}(c_m) \right) \setminus \left( \mathcal{R}'(c_m) \cup \mathcal{R}^{\text{NLS}} \right) \right] \\
	&= \widetilde{\Omega}_\infty \setminus \left( \mathcal{R}'_\infty \cup \mathcal{R}^{\text{NLS}} \right),
\end{aligned}$$
where
 $$\widetilde{\Omega}_\infty :=
  \lim_{m \to \infty} \widetilde{\Omega}_{\mathcal{N}_m}(c_m)\quad\text{and}\quad \mathcal{R}'_\infty := \lim_{m \to \infty} \mathcal{R}'(c_m).$$ 
  Here, the limit is taken in the sense of set convergence on $\mathbb{R}^{\mathbb{Z}}$. Since $\mathcal{N}_m \to \infty$, the high-frequency sets $\widehat{\Pi}'_{\mathcal{N}_m}(c_m)$  shrink to the empty index set and are thus  absorbed into $\widetilde{\Omega}_{\infty}$.

 A parallel argument for the NLS case using \eqref{Xi^NLS} shows that these limiting sets coincide, allowing us to define
 $$\Xi_\infty := \lim_{m \to \infty} \Xi'(c_m) = \lim_{m \to \infty} \Xi^{\text{NLS}}(c_m).$$
 Recalling \eqref{260513-6} and invoking the continuity of the measure, we obtain
$$\text{meas}(\Pi'_\infty \setminus \Xi_\infty) \le \limsup_{m \to \infty} \text{meas}(\Pi'(c_m) \setminus \Xi'(c_m)).$$
By the uniformity of the estimate \eqref{260525-3} with respect to $m$, passing to the limit yields
$$\text{meas}(\Pi'_\infty \setminus \Xi_\infty)= O(\gamma^{1/3}) = O(\gamma').$$

Similarly, using \eqref{mes_Pi^NLS}, the same argument yields 
$$\text{meas}(\Pi^{\text{NLS}} \setminus \Xi_\infty)
 = O(\gamma').$$

This concludes the proof of \eqref{mes_comm1}.
\end{proof}

	\section{Measure Estimate}\label{sec 5}
	In this section, we will  demonstrate that almost all of $\omega\in\Pi$ and $\omega^{\text{NLS}}\in\Pi^{\text{NLS}}$ satisfy the  nonresonance conditions (\ref{NR_KG1}), (\ref{NR_KG2}) and \eqref{NR_NLS}, respectively. Precisely,  we will prove the following Lemmas hold:
	
	\begin{lemma}\label{Le_meas3}
		For fixed $c\in [1,\infty)$ and a sufficiently samll $\gamma>0$,
		there exists a subset $\mathcal{R}\equiv\mathcal{R}(c)\subset\Pi\equiv\Pi(c)$ 
		satisfying
		\begin{equation}\label{090603-1}
		     \mbox{meas}\  {\mathcal{R}} 
			= O( \gamma^{1/3}),
		\end{equation}
		such that for any $\omega\in\Pi\setminus  {\mathcal{R}}$,
		the nonresonance conditions (\ref{NR_KG1}) and  (\ref{NR_KG2}) hold.
	\end{lemma}

		\begin{lemma}\label{Le_meas6}
		For a sufficiently samll $\gamma>0$,
		there exists a subset $\mathcal{R}^{\text{NLS}}\subset\Pi^{\text{NLS}}$  satisfying
		\begin{equation}\label{090603-11}
			\mbox{meas}\  {\mathcal{R}^{\text{NLS}}} 
			= O( \gamma^{1/3}),
		\end{equation}
		such that for any $\omega^{\text{NLS}}\in\Pi^{\text{NLS}}\setminus  {\mathcal{R}^{\text{NLS}}}$, 
		the nonresonance conditions \eqref{NR_NLS} hold.
	\end{lemma}
To achieve this, we will employ the following lemmas. Note that for any $\omega\in\Pi$ given by \eqref{Pi_c} and $\omega^{\text{NLS}}\in\Pi^{\text{NLS}}$ given by \eqref{Pi}, we can express $\omega= (\omega_{n})_{n\in\mathbb{Z}}$
 and $\omega^{\text{NLS}}= (\omega^{\text{NLS}}_{n})_{n\in\mathbb{Z}} $ as
	\begin{equation}\label{omega'}
		\omega_{n}=c\sqrt{c^2+n^2}+r_n, \quad  r_n\in\left[ 0,\frac{c}{3\sqrt{c^2+n^2}}\right]
	\end{equation}
	and
		\begin{equation}\label{omega'_NLS}
		\omega^{\text{NLS}}_{n}=\frac{1}{2}n^2+r^{\text{NLS}}_n, \quad  r^{\text{NLS}}_n \in\left[ 0,\frac{1}{3}\right].
	\end{equation}

	Firstly, we will prove 	that almost all of $\omega\in\Pi$ and $\omega^{\text{NLS}}\in\Pi^{\text{NLS}}$ satisfy the following nonresonance conditions \eqref{041501} and \eqref{041501-1}, respectively,
	 which is parallel the measure estimate method in \cite{Bourgain2005JFA}.
	
	\begin{lemma}\label{Le_meas1}
		For fixed $c\in [1,\infty)$ and a sufficiently samll $\gamma>0$,
		there exists two subsets $\mathcal{R}^1\equiv\mathcal{R}^1(c)\subset\Pi$ 
		and $\mathcal{R}^{\text{NLS},1}\subset\Pi^{\text{NLS}}$ given by \eqref{Pi_c} satisfying
		\begin{equation}\label{090603}
			\mbox{meas}\,  {\mathcal{R}^1}, 
			\mbox{meas}\, \mathcal{R}^{\text{NLS},1} = O( \gamma^{1/3}),
		\end{equation}
		 such that for any $\omega\in\Pi\setminus \mathcal{R}^1$ and $\omega^{\text{NLS}}\in\Pi^{\text{NLS}}
		 \setminus \mathcal{R}^{\text{NLS},1}$, one has 
		\begin{align}
			&\left|\sum_{ n\in\mathbb{Z} }
			{\ell}_{n} \, \omega_{n}\right|
			\ge  \gamma^{1/3} 
			\left(\prod_{ n\in\mathbb{Z},\ell_n\neq 0}
			\frac{1}{|\ell_{n}|^2\, \lfloor n \rfloor^3}\right),\label{041501}\\
			&\left|\sum_{ n\in\mathbb{Z} }
			{\ell}_{n} \, \omega^{\text{NLS}}_{n}\right|
			\ge \gamma^{1/3} 
			\left(\prod_{ n\in\mathbb{Z},\ell_n\neq 0}
			\frac{1}{|\ell_{n}|^2\, \lfloor n \rfloor^3}\right)\label{041501-1}.
		\end{align}
		for any $\ell\in\mathbb{Z}^{\mathbb{Z}}$ with $0<|\ell|<\infty$.
	\end{lemma}
	\begin{proof}
		Let
		$$\mathcal{R}^1:=\bigcup_{{\ell}\in\mathbb{Z}^{\mathbb{Z}},\, 0<|\ell|<\infty}
		\mathcal{R}^1_{\ell}
		\qquad\text{and}\qquad
		\mathcal{R}^{\text{NLS},1}:=\bigcup_{{\ell}\in\mathbb{Z}^{\mathbb{Z}},\, 0<|\ell|<\infty}
		\mathcal{R}^{\text{NLS},1}_{\ell},$$
		where the resonant sets
		$\mathcal{R}_{\ell}^1$ and $\mathcal{R}^{\text{NLS},1}_{\ell}$ are defined by
		\begin{align}
			&\mathcal{R}^1_{\ell}=\left\{ \omega \in \Pi:
			\left|\sum_{ n\in\mathbb{Z} }
			{\ell}_{n} \, \omega_{n}\right|
			<  \gamma^{1/3} 
			\left(\prod_{ n\in\mathbb{Z},\ell_n\neq 0}
			\frac{1}{|\ell_{n}|^2\, \lfloor n \rfloor^3}\right)  \right\},\label{res1}\\
			&\mathcal{R}^{\text{NLS},1}_{\ell}=\left\{ \omega^{\text{NLS}} \in \Pi^{\text{NLS}}:
			\left|\sum_{ n\in\mathbb{Z} }
			{\ell}_{n} \, \omega^{\text{NLS}}_{n}\right|
			<  \gamma^{1/3} 
			\left(\prod_{ n\in\mathbb{Z},\ell_n\neq 0}
			\frac{1}{|\ell_{n}|^2\, \lfloor n \rfloor^3}\right)
			\right\}.\label{res1_NLS}		
			\end{align}
		It is easy to see that for each $\omega\in\Pi\setminus \mathcal{R}^1$
		 and $\omega^{\text{NLS}}\in\Pi^{\text{NLS}}\setminus 
		  {\mathcal{R}^{\text{NLS},1}}$,
		  the estimates (\ref{041501}) and (\ref{041501-1}) hold respectively. Then
		it suffices to prove (\ref{090603}).

		Let
		\begin{equation*}
			m(\ell)
			:=\min\left\{\lfloor n\rfloor \in\mathbb{Z}:
			{\ell}_{{n}}\neq 0\right\}.
		\end{equation*}
		In view of  \eqref{Pi_c} and \eqref{res1}, one has
		\begin{align*} 
			\mbox{meas}\ \mathcal{R}^1_{\ell}
			\le 6 \gamma^{1/3}\cdot  m(\ell)
			\left(\prod_{ n\in\mathbb{Z},\ell_n\neq 0}
			\frac{1}{|\ell_{n}|^2\, \lfloor n \rfloor^3}\right)
		\end{align*}
		and then
		\begin{align}	\label{032902}
			\mbox{meas}\ \mathcal{R}^1
			&\leq 6 \gamma^{1/3}
			\sum_{{\ell}\in\mathbb{Z}^{\mathbb{Z}},\, 0<|\ell|<\infty}
			m(\ell)  \left(\prod_{ n\in\mathbb{Z},\ell_n\neq 0}
			\frac{1}{|\ell_{n}|^2\, \lfloor n \rfloor^3}\right)
			\nonumber\\
			&=6 \gamma^{1/3}  \sum_{s\in\mathbb{Z}}s
			\sum_{(\ell_n)_{|n|\ge s},\ell_n\in\bar{\mathbb{Z}}}
			\left(\prod_{ n\in\mathbb{Z},|n|\ge s}
			\frac{1}{|\ell_{n}|^2\, \lfloor n \rfloor^3}\right)
			\nonumber\\
			&\le {C_1}\gamma^{1/3} 
			\sum_{s\in\mathbb{Z}}s^{-2}
			\prod_{n\in\mathbb{Z},|n|>s} 	
			\left(\sum_{ \ell_n\in \bar{\mathbb{Z}} }
			\frac{1}{|\ell_{n}|^2\, \lfloor n \rfloor^3}\right) \nonumber\\
			&\leq {C_1}\gamma^{1/3} \sum_{s\in\mathbb{Z}}s^{-2}
			\prod_{n\in\mathbb{Z},|n|>s} 
			\left(1+ \frac{{C_1}}{n^3} \right)
			\leq {C}_2\gamma^{1/3}, 
		\end{align}
	where $\bar{\mathbb{Z}}:= \mathbb{Z}\setminus\{0\}$ and $C_1,C_2>0$ are absolute constants.
	
	Similarly, following the proof of \eqref{032902}, we also obtain
	\begin{equation}\label{250806-9}
			\mbox{meas}\, \mathcal{R}^{\text{NLS},1}
			=O(\gamma^{1/3}).
	\end{equation}
	Thus, in view of \eqref{032902} and \eqref{250806-9}, we finish the proof of (\ref{090603}). 
	\end{proof}
	
	Next, we show that almost all of $\omega\in\Pi$ satisfy the following nonresonance conditions \eqref{NR4} and \eqref{NR3}. 
	For any $\ell\in\mathbb{Z}^{\mathbb{Z}}$ with $3\le |\ell|<\infty$, let
	$$
A(\ell):=2 \prod_{|n|\leq n_1^{*}(\ell),\,\ell_{n}\neq 0}
|\ell_{n}|\, \lfloor n \rfloor $$
and
	$$
B(\ell):=2 \prod_{|n|\leq n_3^{*}(\ell),\,\ell_{n}\neq 0}
|\ell_{n}|\, \lfloor n \rfloor.$$

		\begin{lemma}\label{Le_meas7}
		For fixed $c\in [1,\infty)$ and a sufficiently samll $\gamma>0$,
		there exists a subset $\mathcal{R}^2\equiv\mathcal{R}^2(c)\subset\Pi
		$  satisfying
		\begin{equation}\label{090603-2}
			\mbox{meas}\  {\mathcal{R}^2} = O( \gamma^{1/3}),
		\end{equation}
		such that for any $\omega\in\Pi\setminus \mathcal{R}^2$, one has 
		\begin{equation}\label{NR4}
			\left| \sum_{ n\in\mathbb{Z} }
			{\ell}_{n} \, \omega_{n} \right|
			\ge  \gamma^{1/3}c^2	\left(
			\prod_{n \in\mathbb{Z}, \, \ell_{n}\neq 0}
			\frac{1}{|\ell_{n}|^5\, \lfloor n \rfloor^6}\right)
		\end{equation}
		for any  $\ell\in\mathbb{Z}^{\mathbb{Z}}$ with $3\le |\ell|<\infty$ and $\alpha(\ell)\neq 0$.
	\end{lemma}
	\begin{proof}
			Let
		$$\mathcal{R}^2:=
		\bigcup_{ {\ell}\in\mathbb{Z}^{\mathbb{Z}}
			,\,   3\le |\ell|<\infty,\alpha(\ell)\neq 0 }
		\mathcal{R}^2_{\ell},$$
		where the resonant set
		$\mathcal{R}^2_{\ell}$ is defined by
		\begin{align}\label{260601-2}
			\mathcal{R}^2_{\ell}=\left\{ \omega \in \Pi:
			\left| \sum_{ n\in\mathbb{Z} }
			{\ell}_{n} \, \omega_{n} \right|
			<  \gamma^{1/3}c^2	\left(
			\prod_{n\in \mathbb{Z}}
			\frac{1}{|\ell_{n}|^5\, \lfloor n \rfloor^6}\right) 
			\right\}.
		\end{align}
		Then it is easy to see that for each $\omega\in\Pi\setminus \mathcal{R}^2$, the estimates (\ref{NR4}) hold and it suffices to prove (\ref{090603-2}).

			\textbf{Case 1.} $c> \sqrt{2} A(\ell).$	
			
				By using \eqref{omega'}, we rewrite
			\begin{equation}\label{250806-8}
				\omega_{n}
				=c^2+\omega_{n}'
				\quad\text{with}\quad
				\left|\omega_n'\right|\leq n^2+1,
			\end{equation}
			which implies 
			\begin{align*}
				\left| 	\sum_{ n\in\mathbb{Z}}
				{\ell}_{n} \, \omega_{n}' \right|
				\le A^2(\ell).
			\end{align*}	
				Then, since \(\alpha(\ell) \neq 0\),  we can subsequently deduce that
			\begin{align*}
				\left|	\sum_{ n\in\mathbb{Z} }
				{\ell}_{n} \, \omega_{n} \right| 
				\geq  \left|\alpha(\ell)\cdot c^{2}\right|
				-\left| \sum_{ n\in\mathbb{Z} }
				{\ell}_{n} \, \omega_{n}' \right|
				\geq  c^{2}-A^2(\ell) \geq \frac{1}{2}c^2,
			\end{align*}	
			which implies the resonant set $\mathcal{R}^2_{\ell}$ is empty.

			\textbf{Case 2.} $c\le \sqrt{2} A(\ell).$

			In this case, since $c$ is bounded from above by a constant depending only on $\ell$, the factor $c^2$ in the definition \eqref{260601-2} of $\mathcal{R}_\ell^2$ can be readily absorbed by relaxing the decay exponents of the denominators. Then, following the standard argument as in  \eqref{090603}, we directly obtain \eqref{090603-2}.
	\end{proof}

	\begin{lemma}\label{Le_meas2}
		For fixed $c\in [1,\infty)$ and a sufficiently samll $\gamma>0$,
		there exists a subset
		$\mathcal{R}^3\equiv\mathcal{R}^3(c)\subset\Pi$  satisfying
		\begin{equation}\label{090603-21}
			\mbox{meas}\  {\mathcal{R}^3} = O( \gamma^{1/3}),
		\end{equation}
		such that for any $\omega\in\Pi\setminus \mathcal{R}^3$, one has 
		\begin{equation}\label{NR3}
			\left|\left( \sum_{ |n| \leq  n_3^{*}(\ell)}
			{\ell}_{n} \, \omega_{n}\right)  + bc \right|
			\ge  \gamma^{1/3}	\left(
			\prod_{|n| \leq  n_3^{*}(\ell), \, \ell_{n}\neq 0}
			\frac{1}{|\ell_{n}|^5\, \lfloor n \rfloor^6}\right) 
		\end{equation}
		for any $b\in\mathbb{Z}$ with $ |b|<\infty$ and $\ell\in\mathbb{Z}^{\mathbb{Z}}$ with $3\le |\ell|<\infty$, provided that $\ell$ is chosen such that
		$2\sqrt{2} B^2(\ell)\ge c$.
		
	\end{lemma}
	\begin{proof}
		Let
		$$\mathcal{R}^3:=
		\bigcup_{
			{\ell}\in\mathbb{Z}^{\mathbb{Z}}
			,\,   3\le |\ell|<\infty,2\sqrt{2} B^2(\ell)\ge c
			\atop{b\in\mathbb{Z},\, |b|<\infty }}
		\mathcal{R}^3_{\ell,b},$$
		where the resonant set
		$\mathcal{R}^3_{\ell,b}$ is defined by
		\begin{align*}
			\mathcal{R}^3_{\ell,b}=\left\{ \omega \in \Pi:
			\left|\left( \sum_{ |n| \leq  n_3^{*}(\ell)}
			{\ell}_{n} \, \omega_{n}\right)  + bc \right|
			<  \gamma^{1/3}	\left(
			\prod_{|n| \leq  n_3^{*}(\ell), \, \ell_{n}\neq 0}
			\frac{1}{|\ell_{n}|^5\, \lfloor n \rfloor^6}\right) 
			\right\}.
		\end{align*}
		Then it is easy to see that for each $\omega\in\Pi\setminus \mathcal{R}^3$, the estimates (\ref{NR3}) hold and it suffices to prove (\ref{090603-21}).

		Note that by using \eqref{omega'}
		\begin{align*}
			\left|  \sum_{|n| \leq  n_3^{*}(\ell)}
			{\ell}_{n} \, \omega_{n} \right| 
			&\le\sum_{|n| \leq  n_3^{*}(\ell)}
			\left| {\ell}_{n}\right| \left( c\sqrt{c^2+n^2} +1\right) \nonumber\\
			&\le\sum_{ |n| \leq  n_3^{*}(\ell)}
			\left| {\ell}_{n}\right| \cdot 2c^2 |n|
			\le c^2 B(\ell),
		\end{align*}
		then if $|b|> c B(\ell)+1$ one has
		\begin{equation*}
			\left| \left( \sum_{ |n| \leq  n_3^{*}(\ell)}
			{\ell}_{n} \, \omega_{n}\right)  + bc \right|\ge
			|b|c- \left| \sum_{ |n| \leq  n_3^{*}(\ell)}
			{\ell}_{n} \, \omega_{n} \right| \ge 1,
		\end{equation*}
		which is not small. Therefore, we only need to consider 
		$|b|\le c B(\ell) +1 \le 4B^3(\ell) $.  And  we obtain
		\eqref{090603-21} by following the proof of \eqref{090603}.
	\end{proof}

	Finally, we will complete the proof of Lemma \ref{Le_meas3} by using Lemma \ref{Le_meas1}-\ref{Le_meas2}. 	
	In fact, according to  Lemma \ref{Le_meas1}, the core of the proof for Lemma \ref{Le_meas3} lies in showing that the first and second largest indices  $n_1^*(\ell),n_2^*(\ell)$ can be controlled by the indices $\left( n_i^*(\ell)\right)_{i\ge3}$,
	 with this upper bound being independent of \( c \).

	\textbf{Proof of Lemma \ref{Le_meas3}}
	\begin{proof}
		
		Let 
		$$\mathcal{R}'=\bigcup_{{\ell}\in\mathbb{Z}^{\mathbb{Z}}, \, 3\le |\ell|\le \infty }\mathcal{R}'_{\ell}
		\qquad\text{and}\qquad
		\mathcal{R}^*=\bigcup_{{\ell}\in\mathbb{Z}^{\mathbb{Z}}, \, 3\le |\ell|\le \infty }\mathcal{R}^*_{\ell},
		$$ 
		where the resonant sets $\mathcal{R}'_{\ell}$ and  $\mathcal{R}_{\ell}^*$ 
	are defined by\\
		\textbf{(1)} for all $\ell$,
		\begin{align*} 
			\mathcal{R}'_{\ell}&:=\left\{ \omega \in\Pi:
			\left|	\sum_{{n}\in\mathbb{Z}}	{\ell}_{n} \, \omega_{n} 	\right|
			<  {\gamma} \left( 
			\prod_{|n| \leq  n_3^{*}(\ell),\, \ell_n\neq 0}  
			\frac{1}{ |\ell_{n}|^5\,\lfloor n\rfloor^6}\right)^{5}
			\right\};
		\end{align*}
			\textbf{(2)} when $\alpha(\ell)\neq 0$,
		\begin{align*}
			\mathcal{R}_{\ell}^*&:=\left\{ \omega \in\Pi:
			\left|	\sum_{{n}\in\mathbb{Z}}	{\ell}_{n} \, \omega_{n} 	\right|
			<  {\gamma c^2} \left( 
			\prod_{|n| \leq  n_3^{*}(\ell),\, \ell_n\neq 0}  
			\frac{1}{ |\ell_{n}|^5\,\lfloor n\rfloor^6}\right)^{6}
			\right\}.
		\end{align*} 
		Furthermore, let
		$$ \mathcal{R}= \mathcal{R}' \bigcup  \mathcal{R}^*,$$
		then it is easy to see that for each $\omega\in\Pi\setminus \mathcal{R}$
		 the estimates (\ref{NR_KG1}) and (\ref{NR_KG2})
		  hold, and  it suffices to prove (\ref{090603-1}) in the following two cases.
		
		\textbf{Case 1.}	$  n_2^*(\ell)=n_3^*(\ell)$.
		
		\textbf{Subcase 1.1.} $ n_1^*(\ell)= n_2^*(\ell)$.
		
		In this case,  since $n_1^*(\ell) = n_2^*(\ell) = n_3^*(\ell)$, one has 
		\begin{equation*}
			\prod_{|n| \leq  n_3^{*}(\ell),\, \ell_n\neq 0}  
			\frac{1}{ |\ell_{n}|^5\,\lfloor n\rfloor^6}=\prod_{|n| \leq  n_1^{*}(\ell),\, \ell_n\neq 0}  
			\frac{1}{ |\ell_{n}|^5\,\lfloor n\rfloor^6}=	\prod_{n\in\mathbb{Z},\, \ell_n\neq 0}  
			\frac{1}{ |\ell_{n}|^5\,\lfloor n\rfloor^6}.
		\end{equation*}  
		Consequently, by applying a similar measure argument as in the proofs of Lemma \ref{Le_meas1}-\ref{Le_meas7} (specifically using (\ref{090603}) and \eqref{090603-2}), we obtain (\ref{090603-1}).

		\textbf{Subcase 1.2.} $ n_1^*(\ell)>n_2^*(\ell)$.
		
		In this case, one has 
		\begin{equation*}
			\prod_{|n| \leq  n_3^{*}(\ell),\, \ell_n\neq 0}  
			\frac{1}{ |\ell_{n}|^5\,\lfloor n\rfloor^6}=\prod_{|n| \leq  n_2^{*}(\ell),\ \ell_n\neq 0}  
			\frac{1}{ |\ell_{n}|^5\,\lfloor n\rfloor^6}.
		\end{equation*}
		To control the index $n_1^*(\ell)$, we invoke the momentum conservation condition \eqref{mome}, which yields 
		\begin{align*}
			|n_1^*(\ell) |
			\le \sum_{ |n| \le n_3^*(\ell) }  |\ell_{n}| |n|
			\le \prod_{ |n| \le n_3^*(\ell), \, \ell_n\neq 0}   |\ell_{n}| \,\lfloor n\rfloor .
		\end{align*}
		Consequently, the desired estimate \eqref{090603-1} can be established from the proofs of \eqref{090603} and \eqref{090603-2}.

		\textbf{Case 2.} $ n_2^*(\ell)>n_3^*(\ell)$.

		Rewrite
		\begin{align*}
		\sum_{ n\in\mathbb{Z} }	{\ell}_{n} \, \omega_{n} 
		&=\left(\sum_{ n\in\mathbb{Z},\, |n| \leq  n_3^{*}(\ell)}
		{\ell}_{n} \, \omega_{n}\right)
		+\varepsilon_1 \omega_{n_1} 
		+\varepsilon_2 \omega_{n_2}, 
		\end{align*}
		where $\varepsilon_1,\varepsilon_2\in\{-1,1\}$,
		$|n_1|=n_1^*(\ell)$ and $|n_2|=n_2^*(\ell)$.

			Furthermore, we always assume that
		\begin{equation}\label{250806-2}
			n_2^*(\ell)> B^2(\ell)=4 \prod_{|n|\leq n_3^{*}(\ell),\,\ell_{n}\neq 0}
			|\ell_{n}|^2\, \lfloor n \rfloor^2,
		\end{equation} 
		since the case $n_2^*(\ell) \le  B^2(\ell)$ can be treated in the same manner as \textbf{Case 1}.

		Next, we only need to consider \eqref{250806-2}.

		\textbf{Subcase 2.1.} $\varepsilon_1\cdot\varepsilon_2=1.$
		
		Without loss of generality, assume
		$
		\varepsilon_{1}=\varepsilon_{2}=-1,
		$ so that
		\begin{equation*}
			\left(\sum_{ n\in\mathbb{Z},\, |n| \leq  n_3^{*}(\ell)}
			{\ell}_{n} \, \omega_{n}\right)
			+\varepsilon_1 \omega_{n_1} 
			+\varepsilon_2 \omega_{n_2} 
			=\left(\sum_{ n\in\mathbb{Z},\, |n| \leq  n_3^{*}(\ell)}
			{\ell}_{n} \, \omega_{n}\right)
			- \omega_{n_1}- \omega_{n_2} .
		\end{equation*}
		
		Applying \eqref{250806-8} yields
		\begin{align}\label{mes10}		
			\left| 	\sum_{ n\in\mathbb{Z},\, |n| \leq  n_3^{*}(\ell)}
			{\ell}_{n} \, \omega_{n}' \right|
			\le B^2(\ell).
		\end{align}	
		Setting
			\begin{equation}\label{041702}
			\alpha'(\ell):=\sum_{ n\in\mathbb{Z},\, |n| \leq  n_3^{*}(\ell)}	{\ell}_{n},
		\end{equation}
		it follows from (\ref{mes10}) and the triangle inequality that if
		$$\left| \alpha'(\ell)\cdot c^2- \omega_{n_1}- \omega_{n_2} \right|> B^2(\ell)+1,$$
		then
		\begin{align*}
			\left| \left(\sum_{ n\in\mathbb{Z},\, |n| \leq  n_3^{*}(\ell)}
			{\ell}_{n} \, \omega_{n}\right)
			- \omega_{n_1} - \omega_{n_2} \right| \geq 1,
		\end{align*}
		which is not small.
		Thus, one just focuses on the case 
		\begin{equation}\label{count1}
			\left| \alpha'(\ell)\cdot c^2- \omega_{n_1}- \omega_{n_2} \right|\leq B^2(\ell)+1.
		\end{equation}
		It follows that 
		the number of $(n_1,n_2)\in\mathbb{Z}^2$ satisfying (\ref{count1}) is less than
		$ 32B^7(\ell) $.
		Consequently, the desired estimate \eqref{090603-1} also can be established from the proofs of \eqref{090603} and \eqref{090603-2}.

				\textbf{Subcase 2.2.} $\varepsilon_1\cdot\varepsilon_2=-1.$
				
				Without loss of generality, assume 
				$ \varepsilon_{1}=1,\, \varepsilon_{2}=-1,$ so that
				\begin{equation*}
					\left(\sum_{ n\in\mathbb{Z},\, |n| \leq  n_3^{*}(\ell)}
					{\ell}_{n} \, \omega_{n}\right)
					+\varepsilon_1 \omega_{n_1} 
					+\varepsilon_2 \omega_{n_2} 
					=\left(\sum_{ n\in\mathbb{Z},\, |n| \leq  n_3^{*}(\ell)}
					{\ell}_{n} \, \omega_{n}\right)
					+\omega_{n_1}- \omega_{n_2} .
				\end{equation*}
				When 
				$$\text{sgn} (n_1 )\cdot \text{sgn} (  n_2) =-1,$$  the momentum conservation condition \eqref{mome} implies that
				\begin{equation}\label{041701}
					|n_1|+|n_2|=|n_1 - n_2| \le B(\ell).
				\end{equation}
				It follows that the number of $(n_1,n_2)\in\mathbb{Z}^2$ satisfying (\ref{041701}) is less than $4B^2(\ell)$, thereby establishing \eqref{090603-1}.
				Consequently, the remaining analysis can be restricted to the case where 
				$$\text{sgn}(n_1)\cdot\text{sgn}(n_2) = 1,$$ for which the condition \eqref{mome} yields
				\begin{equation}\label{mome3}
					|n_1|-|n_2|=|n_1 - n_2| \le B(\ell).
				\end{equation}

				\textbf{Subcase 2.2.1.} $c> 2\sqrt{2} B^2(\ell).$

				
				
			In this case, by a direct calculation we have 
				\begin{equation}\label{mes19}
					2B^2(\ell) -1 \le \left| \omega_{n_1}-\omega_{n_2}\right|  \le  c B(\ell) +1.
				\end{equation}

			On one hand, when $\alpha(\ell)=0$, it follows that $\alpha'(\ell) = 0$ (see (\ref{041702})). Utilizing (\ref{mes10}) and (\ref{mes19}), we derive
				\begin{align*}
					\left|	\sum_{ n\in\mathbb{Z} }
					{\ell}_{n} \, \omega_{n}
					\right| 
					&\geq  \left| 	\omega_{n_1}-\omega_{n_2}  \right|
					-\left| \sum_{ n\in\mathbb{Z},\ |n| \leq  n_3^{*}(\ell) }
					{\ell}_{n} \, \omega_{n}' \right|
					\geq 	2B^2(\ell) -1 -B^2(\ell) \geq 1.
				\end{align*}
				
				On the other hand, when $\alpha(\ell)\neq 0$, we have $\alpha'(\ell)\neq0$. Applying (\ref{mes10}) and (\ref{mes19}) once again, we can deduce that
				\begin{align*}
					\left|	\sum_{ n\in\mathbb{Z} }
					{\ell}_{n} \, \omega_{n}
					\right| 
					&\geq  \left|\alpha'(\ell)\cdot c^{2}\right|- 
					\left| 	\omega_{n_1}-\omega_{n_2}  \right|
					-\left| \sum_{ n\in\mathbb{Z} ,|n| \leq  n_3^{*}(\ell)}
					{\ell}_{n} \, \omega_{n}' \right| \nonumber\\
					&\geq  c^{2}-c B(\ell)  -1-B^2(\ell) \geq \frac{1}{2}c^2.
				\end{align*}	
				
				Consequently, the combination of these two scenarios implies that the resonant sets $\mathcal{R}_{\ell}'$ and  $\mathcal{R}_{\ell}^*$ are empty, thereby verifying \eqref{090603-1} trivially.

				\textbf{Subcase 2.2.2.}
				$c\le 2\sqrt{2} B^2(\ell).$
				
			In this case, for any ${\ell}\in\mathbb{Z}^{\mathbb{Z}}$ with $0<|\ell|<\infty$,
			note that the resonant set \( \mathcal{R}_{\ell} \) 
			is a subset of the following resonant set
			\begin{align*}
				\mathcal{R}_{\ell}''&:=\left\{ \omega \in\Pi:
				\left|	\sum_{{n}\in\mathbb{Z}}	{\ell}_{n} \, \omega_{n} 	\right|
				<  {\gamma} \left( 
				\prod_{|n| \leq  n_3^{*}(\ell),\, \ell_n\neq 0}  
				\frac{1}{ |\ell_{n}|^5\,\lfloor n\rfloor^6}\right)^{5}
				\right\}
			\end{align*}
			and let
			$$
			\mathcal{R}'':= \bigcup_{{\ell}\in\mathbb{Z}^{\mathbb{Z}}, \, 3\le |\ell|\le \infty }\mathcal{R}''_{\ell}.
			$$ 
		To establish equation (\ref{090603-1}), it suffices to demonstrate that
			 \begin{equation}\label{090603-1'}
			 	\mbox{meas}\  {\mathcal{R}''} 
			 	= O( \gamma^{1/3}),
			 \end{equation}
			And the measure estimate method employed is similar to the one in Lemma \ref{Le_meas1} of  \cite{CY2021}. For the sake of completeness, we provide a detailed proof below.
				
				We rewrite by \eqref{omega'}
				\begin{equation*}
					\omega_n=c\left( c+|n|+\theta_n\right) +r_n
				\end{equation*}
				with
				\begin{equation}\label{theta_n}
					\theta_n=\frac{-2c|n|}{(c+|n|)+\sqrt{c^2+n^2}},
				\end{equation}
				then
				\begin{align}\label{n1-n2}
					\omega_{n_1}- \omega_{n_2}  
					= c\left( |n_1|-|n_2|\right) 	
					+c\left( \theta_{n_1}-\theta_{n_2}\right) +\left( r_{n_1}-r_{n_2}\right).
				\end{align}
				
				In view of \eqref{omega'}, \eqref{mome3} and \eqref{theta_n}, we derive respectively 
				\begin{equation*}\label{omega_n1-n2}
					\left| r_{n_1}-r_{n_2}\right| \le \frac{c}{3\sqrt{c^2+n_1^2}}+\frac{c}{3\sqrt{c^2+n_2^2}}
					\le \frac{c}{|n_2|}
				\end{equation*}
				and
				\begin{align*}\label{theta_n1-n2}
					\left|  \theta_{n_1}-\theta_{n_2}\right| \le 
					\frac{2c( |n_1|-|n_2|)}{(c+|n_2|)+\sqrt{c^2+n_2^2}} 
					\le  \frac{cB(\ell)}{|n_2|}.
				\end{align*}
				Therefore, noting that $c\le 2\sqrt{2} B^2(\ell)$, one obtains
				\begin{equation}\label{upper_n1-n2}
					\left|  c\left( \theta_{n_1}-\theta_{n_2}\right) 
					+\left(  r_{n_1}-r_{n_2} \right) \right|
					\le  \frac{9B^5(\ell)}{|n_2|}.
				\end{equation}
				
				On one hand, when $|n_2|> 18 \gamma^{-1/3} B^{11}(\ell)$, invoking \eqref{NR3}, \eqref{n1-n2}, \eqref{upper_n1-n2}  
			and the triangle inequality yields that for any $\omega \in \Pi \setminus \mathcal{R}^2$ given in Lemma \ref{Le_meas2},
				\begin{align*}
					\left| \left(\sum_{ n\in\mathbb{Z},\, |n| \leq  n_3^{*}(\ell)}
					{\ell}_{n} \, \omega_{n}\right)
					+ \omega_{n_1}- \omega_{n_2}\right| 
					\ge  \frac{ \gamma^{1/3}}{2}	\left(
					\prod_{|n| \leq  n_3^{*}(\ell), \, \ell_{n}\neq 0}
					\frac{1}{|\ell_{n}|^5\, \lfloor n \rfloor^6}\right)
				\end{align*}
				which is not small.

			On the other hand, when $|n_2|\le 18 \gamma^{-1/3} B^{11}(\ell) $, we  obtain $|n_1|\le  19 \gamma^{-1/3} B^{11}(\ell)$ by virtue of (\ref{mome3}), which finishes the proof of (\ref{090603-1}).
			\end{proof}

		\textbf{Proof of Lemma \ref{Le_meas6}}	
		\begin{proof}
			Let 
			$$\mathcal{R}^{\text{NLS}}=\bigcup_{{\ell}\in\mathbb{Z}^{\mathbb{Z}}, \, 3\le |\ell|\le \infty }\mathcal{R}^{\text{NLS}}_{\ell},
			$$ 
			where the resonant set $\mathcal{R}^{\text{NLS}}_{\ell}$
			is defined by
			\begin{align*}
				\mathcal{R}^{\text{NLS}}_{\ell}&:=\left\{ \omega^{\text{NLS}} \in\Pi^{\text{NLS}}:
				\left|	\sum_{{n}\in\mathbb{Z}}	{\ell}_{n} \, \omega^{\text{NLS}}_{n} 	\right|
				<  {\gamma} \left( 
				\prod_{|n| \leq  n_3^{*}(\ell),\, \ell_n\neq 0}  
				\frac{1}{ |\ell_{n}|^5\,\lfloor n\rfloor^6}\right)^{5}
				\right\}.
			\end{align*}
			Then it is easy to see that for each
			$\omega^{\text{NLS}}\in\Pi^{\text{NLS}}\setminus \mathcal{R}^{\text{NLS}}$
			the estimates  \eqref{NR_NLS}
			hold, and  it suffices to prove (\ref{090603-11}).
		
			According to the proof of Lemma \ref{Le_meas3}, we only need to consider the case where 
			$$  n_2^*(\ell)>n_3^*(\ell)
			\quad\text{and}\quad
			n_2^*(\ell)\geq B^2(\ell).$$
			Rewrite
			\begin{align}
				\sum_{ n\in\mathbb{Z} } {\ell}_{n} \,  \omega^{\text{NLS}}_{n} 
				=\left(\sum_{ n\in\mathbb{Z},\, |n| \leq  n_3^{*}(\ell)}
				{\ell}_{n} \, \omega^{\text{NLS}}_{n}\right)
				+\varepsilon_1 \omega^{\text{NLS}}_{n_1} 
				+\varepsilon_2 \omega^{\text{NLS}}_{n_2},\label{250806-1}
			\end{align}
			where $\varepsilon_1,\varepsilon_2\in\{-1,1\}$,
			$|n_1|=n_1^*(\ell)$ and $|n_2|=n_2^*(\ell)$.
			The following proof is divided into two cases.
			
				\textbf{Case 1.}
			$|n_1|\neq |n_2|$ or $ \varepsilon_1\cdot\varepsilon_2=1.$

			By using \eqref{omega'_NLS}, one has
			\begin{equation}
				\left| \varepsilon_1 \omega^{\text{NLS}}_{n_1} 
				+\varepsilon_2 \omega^{\text{NLS}}_{n_2} \right|
			 \ge \frac{1}{2}\left( |n_1|+|n_2|\right) -1
				\ge \frac{1}{2} |n_2|,
			\end{equation}
			then
			\begin{align*}
				\left| \sum_{ n\in\mathbb{Z} }	{\ell}_{n} \, \omega^{\text{NLS}}_{n}\right| 
				&\ge	\left| \varepsilon_1 \omega^{\text{NLS}}_{n_1} 
				+\varepsilon_2 \omega^{\text{NLS}}_{n_2}\right| -
				\left|\left(\sum_{ n\in\mathbb{Z},\, |n| \leq  n_3^{*}(\ell)}
				{\ell}_{n} \, \omega^{\text{NLS}}_{n}\right) \right| \nonumber\\
				&\ge \frac{1}{2} B^2(\ell)
				-\prod_{|n|\leq n_3^{*}(\ell),\,\ell_{n}\neq 0}
				|\ell_{n}|\, \lfloor n \rfloor^2\ge 1,
			\end{align*}
			which is not small. 
			
				\textbf{Case 2.}
			$|n_1| = |n_2|$ and $ \varepsilon_1\cdot\varepsilon_2=-1.$
			
			It follows from the momentum conservation condition \eqref{mome}
			that
			\begin{equation*}
				|n_1|=|n_2|\le B(\ell).
			\end{equation*}
		We then complete the proof of (\ref{090603-1}) in a manner analogous to the proof of Lemma \ref{Le_meas3}.
		\end{proof}

			\section{Appendix}\label{sec 6}

	\subsection{Proof of Lemma \ref{063004}}\label{Le3.5}
	
	\begin{proof}
		We first focus on the estimate \eqref{vector_1}
		 for the mapping $p \rightarrow p$.
		In light of (\ref{042501}), for any fixed \( j \in \mathbb{Z} \), proving \eqref{vector_1} reduces to estimating an upper bound for
		$$
		\left|\frac{\partial R}{\partial z_j} e^{r\ln ^\sigma\lfloor j \rfloor} \lfloor j \rfloor^{p}
		\right|.
		$$

		In view of (\ref{pertur1}), for any $a,k,k^{\prime}\in\mathbb{N}^{\mathbb{Z}}$, it follows that
		\begin{equation}\label{260408-1}
			\frac{\partial R}{\partial z_j}=
			\sum_{a, k, k^{\prime} \in \mathbb{N}^{\mathbb{Z}}} 
			R_{a k k^{\prime}}
			\left(\prod_{n \neq j} I_n(0)^{a_n} z_n^{k_n} \bar{z}_n^{k_n^{\prime}}\right)
			\left(k_j I_j(0)^{a_j} z_j^{k_j-1} \bar{z}_j^{k_j^{\prime}}\right),
		\end{equation}
		where the coefficient $R_{a k k^{\prime}}$
		satisfies the following inequality according to definition \eqref{norm+_1}
%
		\begin{equation*}
		\left( \prod_{n\in\mathbb{Z}} \left\langle n/c\right\rangle^{\frac{1}{2}
			\left(2a_{n}+k_{n}+k_{n}' \right)} \right) 	\left|R_{a k k^{\prime}}\right| 
			\leq
			\|R \|_{2\rho}^{c,+}\cdot
			e^{2\rho\left(\sum_{n \in \mathbb{Z}}\left(2 a_n+k_n+k_n^{\prime}\right) \ln ^\sigma\lfloor n\rfloor-2 \ln ^\sigma\left\lfloor n_1^*\right\rfloor\right)}.
		\end{equation*}
	Taking into account the fact that
		$\left\langle n/c\right\rangle\ge 1$, 
	we then obtain
			\begin{equation}\label{2410121}
				\left|R_{a k k^{\prime}}\right| \leq
			\|R \|_{2\rho}^{c,+}\cdot
			e^{2\rho\left(\sum_{n \in \mathbb{Z}}\left(2 a_n+k_n+k_n^{\prime}\right) \ln^\sigma\lfloor n\rfloor-2 \ln ^\sigma\left\lfloor n_1^*\right\rfloor\right)}.
		\end{equation}
		By substituting \eqref{2410121} into the equation for the derivative \eqref{260408-1} and assuming \(\|z\|_{p} \le 1\), it holds that
			\begin{align}
			\left|\frac{\partial R}{\partial z_j} e^{r\ln ^\sigma\lfloor j \rfloor} \lfloor j \rfloor^{p}
			\right|
			\leq & 	\|R \|_{2\rho}^{c,+}
		 \sum_{a, k, k^{\prime} \in \mathbb{N}^{\mathbb{Z}}} 
			  k_j \cdot
			\left( \prod_{n\in\mathbb{Z}} 
			\lfloor n\rfloor^{-p \left(2 a_n+k_n+k_n^{\prime}\right)} \right) \cdot
			\lfloor j\rfloor^{2p} \nonumber\\
			& \times e^{-r\sum_{n \in \mathbb{Z}}\left(2 a_n+k_n+k_n^{\prime}\right) \ln ^\sigma\lfloor n\rfloor} \cdot
			e^{2r \ln^\sigma\lfloor j \rfloor}  \nonumber\\
			&\times e^{2\rho\left(\sum_{n \in \mathbb{Z}}\left(2 a_n+k_n+k_n^{\prime}\right) \ln ^\sigma\lfloor n\rfloor-2 \ln ^\sigma\left\lfloor n_1^*\right\rfloor\right)}.\label{260408-3}
		\end{align}
	Furthermore, the fact that $\lfloor j\rfloor\le \lfloor n_1^* \rfloor $ along with \eqref{mome} implies
		\begin{align*}
		\left( \prod_{n\in\mathbb{Z}} 
		\lfloor n\rfloor ^{-p \left(2 a_n+k_n+k_n^{\prime}\right)} \right) \cdot
		\lfloor j\rfloor ^{2p}
		\le 
		\left( \prod_{i\ge 1} 
		\lfloor n_i^* \rfloor ^{-p } \right) 	
		\cdot	\lfloor n_1^* \rfloor^{2p} 
		\le 1.
		\end{align*}
		To complete the proof of \eqref{vector_1}, it remains to verify the following claim
		\begin{align}\label{eq:claim_C}
			 (\star)&:=\sum_{a, k, k^{\prime} \in \mathbb{N}^{ \mathbb{Z} } } k_j  
			 e^{(2\rho-r)\left(\sum_{n \in \mathbb{Z}}\left(2 a_n+k_n+k_n^{\prime}\right) \ln ^\sigma\lfloor n\rfloor\right)-4 \rho \ln ^\sigma\left\lfloor n_1^*\right\rfloor+2 r \ln ^\sigma\lfloor j \rfloor}  \nonumber \\
			& \leq	
			\exp \left\{\frac{80}{\rho^{\prime}-\rho} \cdot \exp\left\{\left(\frac{20}{\rho^{\prime}-\rho}\right)^{\frac{1}{\sigma-1}}\right\}\right\}
			\tag{C}.
		\end{align}
		Furthermore, it should be noted that the verification of the estimate \eqref{vector_2} similarly reduces to the proof of the aforementioned claim \eqref{eq:claim_C}. The justification for this reduction is as follows.
		

		 Regarding \eqref{vector_2}  for the mapping $p \rightarrow p-4$:
		For any fixed \( j \in \mathbb{Z} \), it reduces to estimating an upper bound for
		$$	\left|\frac{\partial R}{\partial z_j} e^{r\ln ^\sigma\lfloor j \rfloor} \lfloor j \rfloor^{p-4}
		\right|.	$$
		By invoking the definition \eqref{norm+_3} again, we obtain a more precise coefficient estimate
\begin{align}\label{260408-2}
	\left|R_{a k k^{\prime}}\right| \leq
	h   n^{*4}_1 \cdot \|R \|_{2\rho}^{*,+}\cdot
	e^{2\rho\left(\sum_{n \in \mathbb{Z}}\left(2 a_n+k_n+k_n^{\prime}\right) \ln ^\sigma\lfloor n\rfloor-2 \ln ^\sigma\left\lfloor n_1^*\right\rfloor\right)},
\end{align}		
the factor $h n^{*4}_1$ stems from the estimate $(n^*_1/\sqrt{c})^4 = h n^{*4}_1$, where we have identified $h = 1/c^2$.
By substituting \eqref{260408-2} into the derivative formula \eqref{260408-1}, we obtain an expression structurally identical to \eqref{260408-3}. The only difference is the presence of the additional factor $h n_1^{*4}$ and the replacement of $j^{2p}$ by $j^{2p-4}$. Crucially, since $|j| \leq |n_1^*|$, the product of these extra terms satisfy:
$$h n_1^{*4} \cdot \left( \prod_{n \in \mathbb{Z}} |n|^{-p(2a_n+k_n+k'_n)} \right) \cdot |j|^{2p-4} \leq h.$$ 	
Consequently, the entire estimate for \eqref{vector_2}  is again reduced to the verification of the same core claim \eqref{eq:claim_C}.
		
	We now proceed to establish the proof of the claim \eqref{eq:claim_C}. For clarity, the argument is divided into the following two cases based on the range of the index $|j|$:
		
		{\bf Case 1.} $\lfloor j \rfloor \leq\left\lfloor n_3^*\right\rfloor$.
		
		Firstly, note that
		$$
		e^{-4 \rho \ln^\sigma \left\lfloor n_1^{*}\right\rfloor+2r \ln^\sigma\left\lfloor j \right\rfloor} \leq e^{-4 \rho \ln^\sigma \left\lfloor n_3^{*}\right\rfloor +2 r \ln ^\sigma \left\lfloor n_3^{*}\right\rfloor}=e^{2(r-2\rho) \ln ^\sigma \left\lfloor n_3^{*}\right\rfloor} 
		$$
		and
		\begin{align}
			e^{(2\rho-r) \sum_{n\in\mathbb{Z}}\left(2 a_n+k_n+k_n^{\prime}\right) \ln^{\sigma}\lfloor n \rfloor}
			\leq e^{\frac{1}{3}(2\rho-r) \sum_{i\geq1} \ln^\sigma\left\lfloor n_i^{*}\right\rfloor} \cdot e^{2(2\rho-r) \ln^\sigma \left\lfloor n_3^{*}\right\rfloor}.\nonumber
		\end{align}
		Therefore, one deduces
		$$ (\star) \leq \sum_{a, k, k^{\prime} \in \mathbb{N}^{\mathbb{Z}}} k_j
		 e^{\frac{2\rho-r}{3} \sum_{i \geq 1} \ln ^\sigma\left\lfloor n_i^*\right\rfloor}.	$$
		Hence, in view of (\ref{2410074}), we have
		$$
		\begin{aligned}
		 (\star)& \leq 
			  \sum_{a, k, k^{\prime} \in \mathbb{N}^{\mathbb{Z}}}\left(\sum_{i \geq 1} \ln ^\sigma\left\lfloor n_i^*\right\rfloor\right) e^{\frac{2\rho-r}{3} \sum_{i \geq 1} \ln ^\sigma\left\lfloor n_i^*\right\rfloor} \\
			&\leq   \frac{12}{e(r-2\rho)} \sum_{a, k, k^{\prime} \in \mathbb{N}^{\mathbb{Z}}} e^{\frac{2\rho-r}{4} \sum_{i \geq 1} \ln ^\sigma\left\lfloor n_i^*\right\rfloor} \\
			&\quad\left(\text {due to (\ref{Lem6.3Con}) with $p=1$ and $\delta=\frac{r-2\rho}{12}$ }\right) \\
			&=  \frac{12}{e(r-2\rho)}
			\sum_{a, k, k^{\prime} \in \mathbb{N}^{\mathbb{Z}}} e^{\frac{2\rho-r}{4} \sum_{n\in\mathbb{Z}} \left(2a_n+k_n+k_n^{\prime}\right)\ln ^\sigma\left\lfloor n\right\rfloor} \\
			& \leq \exp \left\{\frac{40}{\rho^{\prime} -\rho} \cdot \exp \left\{\left(\frac{8}{\rho^{\prime}-\rho}\right)^{\frac{1}{\sigma-1}}\right\}\right\},
		\end{aligned}
		$$
		where the last inequality is based on  Claim \ref{a5} with 
		$\delta=\frac{\rho'-\rho}{2} \le\frac{r-2\rho}{4}$.
		
		{\bf Case 2.} $\lfloor j \rfloor>\left\lfloor n_3^*\right\rfloor$.
		
		If $k_j \geq 3$, then one has $\lfloor j \rfloor \leq\left\lfloor n_3^*\right\rfloor$, which is in {\bf Case 1}. Hence we always assume $k_j \leq 2$. Then  we have
		\begin{equation*}
		-4 \rho \ln ^\sigma\left\lfloor n_1^*\right\rfloor+2 r \ln ^\sigma\lfloor j \rfloor\le 	2 (r-2\rho)
		 \ln ^\sigma\left\lfloor n_1^*\right\rfloor 
		\end{equation*}
		and then
		\begin{align}\label{2410152}
			 (\star) &\leq 2\sum_{a, k, k^{\prime} \in \mathbb{N}^{\mathbb{Z}}} 
			e^{(2\rho-r)\left(\sum_{n \in \mathbb{Z}}\left(2 a_n+k_n+k_n^{\prime}\right) \ln ^\sigma\lfloor n\rfloor-2 \ln ^\sigma\left\lfloor n_1^*\right\rfloor \right)}
			\nonumber\\
			&\leq 2\sum_{a, k, k^{\prime} \in \mathbb{N}^{\mathbb{Z}}}  
			e^{\frac{2\rho-r}{2} \sum_{i \geq 3} \ln ^\sigma\left\lfloor n_i^*\right\rfloor }.
		\end{align}

		If $\left\{n_i^*\right\}_{i \geq 1}$ are given, then $\left\{2 a_n+k_n+k_n^{\prime}\right\}_{n \in \mathbb{Z}}$ are specified, and hence $\left(a, k, k^{\prime}\right)$ is specified up to a factor of
		$$
		\prod_n\left(1+l_n^2\right),\quad\text{where}\;
		l_n=\#\left\{j: n_j^*=n\right\}.
		$$
		
		Since $\lfloor j \rfloor>\left\lfloor n_3^*\right\rfloor$, then $j\in\left\{n_1^*, n_2^*\right\}$. Hence, if $\left\{n_i^*\right\}_{i \geq 3}$ and $j$ are given, then $n_1^*$ and $n_2^*$ are uniquely determined. Then, one has
		\begin{align}
	 (\star) &
		\leq 10\sum_{\left\{n_i^*\right\}_{i \geq 3}} \prod_{\substack{n \in \mathbb{Z} \\n \leq n_3^*}}\left(1+l_n^2\right) e^{\frac{2\rho-r}{2} \sum_{i \geq 3} \ln ^\sigma\left\lfloor n_i^*\right\rfloor}\nonumber\\
			& \leq 10\left(\sum_{\left\{n_i^*\right\}_{i \geq 3}} e^{-\frac{\rho^{\prime}-\rho}{2} \sum_{i \geq 3} \ln ^\sigma\left\lfloor n_i^*\right\rfloor}\right) \sup _{\left\{n_i^*\right\}_{i \geq 3}}\left(\prod_{n \leq n_3^*}\left(1+l_n^2\right) e^{-\frac{\rho^{\prime}-\rho}{2} \sum_{i \geq 3} \ln ^\sigma\left\lfloor n_i^*\right\rfloor}\right) \nonumber\\
			&\quad\text{(due to $\rho'-\rho\le \frac{r-2\rho}{2} $ )}\nonumber\\
			& \leq	\exp \left\{\frac{80}{\rho^{\prime}-\rho} \cdot \exp\left\{\left(\frac{20}{\rho^{\prime}-\rho}\right)^{\frac{1}{\sigma-1}}\right\}\right\},\nonumber
		\end{align}
		where the last inequality is based on Claim \ref{a5} with $\delta=\frac{\rho'-\rho}{2} $ and (\ref{24100813}) with $p=2$ and $\delta=\frac{\rho^{\prime}-\rho}{4}$.

	Taking into account the above two cases, we have thus completed the proof of \eqref{eq:claim_C}. 
	\end{proof}

	\subsection{Proof of Lemma \ref{010} and \ref{E1} }
	\label{sec6.1}
	
   	\begin{proof}[Proof of Lemma \ref{010}]
    	
   	We first focus on the estimate \eqref{042704}. 
   	
   	For the Poisson bracket $\{R_1, R_2\}$, let
	$$
	R_1(z, \bar{z})=\sum_{a, k, k^{\prime} \in \mathbb{N}^{\mathbb{Z}} } R_{a k k^{\prime}} \mathcal{M}_{a k k^{\prime}}
	$$
	and
	$$
	R_2(z, \bar{z})=\sum_{A, K, K^{\prime} \in \mathbb{N}^{\mathbb{Z}}} R_{A K K^{\prime}} \mathcal{M}_{A K K^{\prime}},
	$$
		where the so-called monomials
	$$
	\mathcal{M}_{a k k^{\prime}}=\prod_{n \in \mathbb{Z}} I_n(0)^{a_n} z_n^{k_n} \bar{z}_n^{k_n^{\prime}}
	$$
	satisfy the zero momentum conservation \eqref{mome}.
	Then one has
	$$
	\left\{R_1,R_2 \right\}=\sum_{a, k, k^{\prime}, A, K, K^{\prime} \in \mathbb{N}^{\mathbb{Z}}} R_{a k k^{\prime}} R_{A K K^{\prime}}\left\{\mathcal{M}_{a k k^{\prime}}, \mathcal{M}_{A K K^{\prime}}\right\},
	$$
	where
	$$
	\begin{aligned}
		\left\{\mathcal{M}_{a k k^{\prime}}, \mathcal{M}_{A K K^{\prime}}\right\}=&\mathrm{i} \sum_{j\in {\mathbb{Z}}}\left(\prod_{n \neq j} I_n(0)^{a_n+A_n} z_n^{k_n+K_n} \bar{z}_n^{k_n^{\prime}+K_n^{\prime}}\right)\\
		&\qquad\; \times\left(\left(k_j K_j^{\prime}-k_j^{\prime} K_j\right) I_j(0)^{a_j+A_j} z_j^{k_j+K_j-1} \bar{z}_j^{k_j^{\prime}+K_j^{\prime}-1}\right).
	\end{aligned}
	$$
	
	The coefficient of the monomial $\mathcal{M}_{\alpha \kappa \kappa^{\prime}}$ in $	\left\{R_1,R_2 \right\}$ is given by
	\begin{equation}\label{coefficient}
		R_{\alpha \kappa \kappa^{\prime}}=\mathrm{i} \sum_{j\in  {\mathbb{Z}}} \sum_* \sum_{* *}\left(k_j K_j^{\prime}-k_j^{\prime} K_j\right) R_{a k k^{\prime}} R_{A K K^{\prime}},
	\end{equation}
	where
	\begin{equation}\label{coefficient-1}
			\sum_*=\sum_{\substack{a, A \in \mathbb{N}^{\mathbb{Z}} \\ a+A=\alpha}}\qquad \text{and}\qquad
		\sum_{* *}=\sum_{\substack{k, k^{\prime}, K, K^{\prime} \in \mathbb{N}^{\mathbb{Z}} \\ k+K-e_j=\kappa, k^{\prime}+K^{\prime}-e_j=\kappa^{\prime}}}
	\end{equation}	
	with $e_j$ is the $j$-unit vector in $\mathbb{N}^{\mathbb{Z}}$. 
	Let
	$$
	N_1^*=\max \left\{|n|: A_n+K_n+K_n^{\prime} \neq 0\right\},
	$$
	$$
	n_1^*=\max \left\{|n|: a_n+k_n+k_n^{\prime} \neq 0\right\}
	$$
	and
	$$
	\nu_1^*=\max \left\{|n|: \alpha_n+\kappa_n+\kappa_n^{\prime} \neq 0\right\} .
	$$
	
	In view of Definition \ref{def_norm} and Lemma \ref{005}, one has
	\begin{align}\label{2410079}
		&\left|\left(\prod_{n \in \mathbb{Z}} \left\langle\frac{n}{c}\right\rangle ^{\frac{1}{2}\left(2 a_n+k_n+k_n^{\prime}\right)}\right) 
		R_{a k k^{\prime}}\right|\nonumber\\ \leq&\left\|R_1\right\|_{\rho-\delta_1}^{c,+} e^{\rho\left(\sum_{n \in \mathbb{Z}}\left(2 a_n+k_n+k_n^{\prime}\right) \ln^\sigma\lfloor n\rfloor-2 \ln ^\sigma\left\lfloor n_1^*\right\rfloor\right)} e^{-\frac{1}{2} \delta_1\left(\sum_{i \geq 3} \ln ^\sigma\left\lfloor n_i^*\right\rfloor\right)}
	\end{align}
	and
	\begin{align}\label{24100710}
		&\left|\left(\prod_{n \in \mathbb{Z}} \left\langle\frac{n}{c}\right\rangle ^{\frac{1}{2}\left(2 A_n+K_n+K_n^{\prime}\right)}\right)R_{A K K^{\prime}}\right|\nonumber\\ \leq&\left\|R_2\right\|_{\rho-\delta_2}^{c,+} e^{\rho\left(\sum_{n \in \mathbb{Z}}\left(2 A_n+K_n+K_n^{\prime}\right) \ln ^\sigma\lfloor n\rfloor-2 \ln ^\sigma\left\lfloor N_1^*\right\rfloor\right)} e^{-\frac{1}{2} \delta_2\left(\sum_{i \geq 3} \ln ^\sigma\left\lfloor N_i^*\right\rfloor\right)}.
	\end{align}
	In light of the relations in \eqref{coefficient-1} and $ \left\langle{n}/{c}\right\rangle\ge 1$, we observe that
	\begin{align}\label{250519-1}
		\left|\left(\prod_{n \in \mathbb{Z}} \left\langle\frac{n}{c}\right\rangle ^{\frac{1}{2}\left(2 \alpha_n+\kappa_n+\kappa_n^{\prime}\right)}\right) R_{\alpha \kappa \kappa^{\prime}}
		\right|
		\leq \left(\prod_{n \in \mathbb{Z}} \left\langle\frac{n}{c}\right\rangle  ^{\frac{1}{2}\left(2 a_n+k_n+k_n^{\prime}+2 A_n+K_n+K_n^{\prime}\right)}\right)
		\left|	R_{\alpha \kappa \kappa^{\prime}} \right|.
	\end{align}
	Furthermore, by substituting \eqref{coefficient}, (\ref{2410079}) and (\ref{24100710}) into (\ref{250519-1}), we obtain
	\begin{align}\label{260412-3}
		&\quad\;	\left|\left(\prod_{n \in \mathbb{Z}} \left\langle\frac{n}{c}\right\rangle ^{\frac{1}{2}\left(2 \alpha_n+\kappa_n+\kappa_n^{\prime}\right)}\right) R_{\alpha \kappa \kappa^{\prime}}
		\right| \nonumber\\
		&\leq \left\|  R_1 \right\|_{\rho-\delta_1}^{c,+}
		\left\| R_2 \right\|_{\rho-\delta_2}^{c,+} 
		\sum_{j \in \mathbb{Z}} \sum_* \sum_{* *}\left|k_j K_j^{\prime}-k_j^{\prime} K_j\right| e^{\rho\left(\sum_{n \in \mathbb{Z}}\left(2 a_n+k_n+k_n^{\prime}\right) \ln ^\sigma\lfloor n\rfloor-2 \ln ^\sigma\left\lfloor n_1^*\right\rfloor\right)} \nonumber\\
		&\quad \times e^{\rho\left(\sum_{n \in \mathbb{Z}}\left(2 A_n+K_n+K_n^{\prime}\right) \ln ^\sigma\lfloor n\rfloor-2 \ln ^\sigma\left\lfloor N_1^*\right\rfloor\right)}e^{-\frac{1}{2} \delta_1\left(\sum_{i \geq 3} \ln ^\sigma\left\lfloor n_i^*\right\rfloor\right)} e^{-\frac{1}{2} \delta_2\left(\sum_{i \geq 3} \ln ^\sigma\left\lfloor N_i^*\right\rfloor\right)}\nonumber \\
		&=\left\| R_1 \right\|_{\rho-\delta_1}^{c,+}
		\left\| R_2 \right\|_{\rho-\delta_2}^{c,+} e^{\rho\left(\sum_{n \in \mathbb{Z}}\left(2 \alpha_n+\kappa_n+\kappa_n^{\prime}\right) \ln ^\sigma\lfloor n\rfloor-2 \ln ^\sigma\left\lfloor \nu_1^*\right\rfloor\right)} \nonumber\\
		&\quad \times \sum_{j\in \mathbb{Z}} \sum_* \sum_{* *}\left|k_j K_j^{\prime}-k_j^{\prime} K_j\right| e^{2 \rho\left(\ln ^\sigma\lfloor j \rfloor+\ln ^\sigma\left\lfloor \nu_1^*\right\rfloor-\ln ^\sigma\left\lfloor n_1^*\right\rfloor-\ln ^\sigma\left\lfloor N_1^*\right\rfloor\right)}\nonumber\\
		& \quad\times e^{-\frac{1}{2} \delta_1\left(\sum_{i \geq 3} \ln ^\sigma\left\lfloor n_i^*\right\rfloor\right)} e^{-\frac{1}{2} \delta_2\left(\sum_{i \geq 3} \ln ^\sigma\left\lfloor N_i^*\right\rfloor\right)},
	\end{align}
	where the last inequality follows from \eqref{coefficient-1}.
	Then, dividing by the corresponding weight and taking the supremum over $\alpha, \kappa, \kappa'$, we arrive at
	\begin{align*}
		\left\|	\left\{R_1,R_2 \right\} \right\| _{\rho}^{c,+}
		\leq I  \cdot \left\|  R_1 \right\|_{\rho-\delta_1}^{c,+}
		\left\| R_2 \right\|_{\rho-\delta_2}^{c,+} ,
	\end{align*}
	where
	\begin{align*}
		I=\sum_{j \in \mathbb{Z}} \sum_* &\sum_{* *}\left|k_j K_j^{\prime}-k_j^{\prime} K_j\right| 
		e^{2 \rho\left(\ln^\sigma\lfloor j \rfloor+\ln ^\sigma\left\lfloor \nu_1^*\right\rfloor-\ln ^\sigma\left\lfloor n_1^*\right\rfloor-\ln ^\sigma\left\lfloor N_1^*\right\rfloor\right)} \nonumber\\
		&\qquad \times e^{-\frac{1}{2} \delta_1\left(\sum_{i \geq 3} \ln ^\sigma\left\lfloor n_i^*\right\rfloor\right)} e^{-\frac{1}{2} \delta_2\left(\sum_{i \geq 3} \ln ^\sigma\left\lfloor N_i^*\right\rfloor\right)}.
	\end{align*}
To complete the proof of the estimate \eqref{042704}, it remains to show 
	\begin{equation}\label{2410081}
		I \leq \frac{1}{\delta_2} \exp \left\{\frac{1000}{\delta_1} \cdot \exp \left\{\left(\frac{10}{\delta_1}\right)^{\frac{1}{\sigma-1}}\right\}\right\}.
	\end{equation}

Before proceeding with the proof of \eqref{2410081}, let us briefly pause to address the estimate \eqref{250804-3}. Similar to the previous case  \eqref{042704}, the verification of \eqref{250804-3} essentially reduces to proving the estimate \eqref{2410081} for the following reasons.

    
   	In view of Definition \ref{def_norm} and Lemma \ref{005} again, one has
    	\begin{align}\label{24100710'}
    	&\left|  \left( \frac{\sqrt{c}}{N_1^*}\right) ^4 
    	R_{A K K^{\prime}}\right|\nonumber\\ \leq&\left\|R_2\right\|_{\rho-\delta_2}^{*,+} e^{\rho\left(\sum_{n \in \mathbb{Z}}\left(2 A_n+K_n+K_n^{\prime}\right) \ln ^\sigma\lfloor n\rfloor-2 \ln ^\sigma\left\lfloor N_1^*\right\rfloor\right)} e^{-\frac{1}{2} \delta_2\left(\sum_{i \geq 3} \ln ^\sigma\left\lfloor N_i^*\right\rfloor\right)}.
    \end{align}
    Note that $\nu_1^* \geq N_2^*$ according to the relations in \eqref{coefficient-1}, it then follows that
    $$\left( {\sqrt{c}}/{\nu_1^*} \right)^4 \leq \left( {\sqrt{c}}/{N_2^*} \right)^4.$$
    Furthermore, by the momentum conservation condition \eqref{mome}, one has
      \begin{equation*}\label{250412-2}
    	N_1^* \le  \sum_{i \geq 2}  N_i^*
    	\le \left(\sum_{n \in \mathbb{Z}}2 A_n+K_n+K_n^{\prime}\right) N_2^*
    	\le   \left( 2	\sum_{i \geq 3} \ln^\sigma\left\lfloor N_i^*\right\rfloor  \right) N_2^*,
    \end{equation*}
    where the last inequality follows from an estimate that will be established in \eqref{2410075}. Therefore, combining the above two estimates, we obtain the core inequality chain
      \begin{equation*}
    	\left( \frac{\sqrt{c}}{\nu_1^*}\right) ^4 
    	\le     \left( 2	\sum_{i \geq 3} \ln^\sigma\left\lfloor N_i^*\right\rfloor  \right)^4 	
    	\left( \frac{\sqrt{c}}{N_1^*}\right) ^4.
    \end{equation*}
    Consequently, by invoking the property $\langle n/c \rangle \geq 1$ and combining the above estimates, it follows that
    \begin{equation}\label{260412-4}
    	\left|    \left( \frac{\sqrt{c}}{\nu_1^*}\right) ^4
    R_{\alpha \kappa \kappa^{\prime}} \right|\le 
    \left( 2	\sum_{i \geq 3} \ln^\sigma\left\lfloor N_i^*\right\rfloor  \right)^4
    	\left( \frac{\sqrt{c}}{N_1^*}\right) ^4
    	\left(\prod_{n \in \mathbb{Z}} \left\langle\frac{n}{c}\right\rangle^{\frac{1}{2}\left(2 a_n+k_n+k_n^{\prime}\right)}\right) 
    	\left| R_{\alpha \kappa \kappa^{\prime}} \right|.
    \end{equation}
    Following the argument for  \eqref{260412-3} and incorporating \eqref{2410079}, \eqref{24100710'} and \eqref{260412-4}, we can deduce that
    $$\|\{R_1, R_2\}\|_\rho^{*,+} \leq I_1 \cdot \|R_1\|_{\rho-\delta_1}^{c,+} \|R_2\|_{\rho-\delta_2}^{*,+},$$
    where the factor $I_1$ is analogous to $I$, with the primary difference being the presence of an additional factor
    $$\left( 2	\sum_{i \geq 3} \ln^\sigma\left\lfloor N_i^*\right\rfloor  \right)^4.$$
    This factor can be bounded by a constant following the same procedure as in the proof of \eqref{2410081}.
    
	With this reduction in mind, we now proceed to establish (\ref{2410081}). To this end, we first collect several preliminary facts.

	$\bullet$
	In view of the conservation laws in $(\ref{coefficient-1})$, the indices satisfy 
	\begin{align}\label{2410077}
		\left( \sum_{n \in \mathbb{Z}} 2 \alpha_n+\kappa_n+\kappa_n^{\prime}\right)  
		= \left(\sum_{n \in \mathbb{Z}}
		2 a_n+k_n+k_n^{\prime}\right) + 
		\left(\sum_{n \in \mathbb{Z}}2 A_n+K_n+K_n^{\prime}\right)-2.
	\end{align}
	
	$\bullet$ 
	If $j \notin \operatorname{supp}\left(k+k^{\prime}\right) \cap \operatorname{supp}\left(K+K^{\prime}\right)$, then
	$
	k_j K_j^{\prime}-k_j^{\prime} K_j=0.
	$
	Hence we always assume $j\in \operatorname{supp}\left(k+k^{\prime}\right) \cap \operatorname{supp}\left(K+K^{\prime}\right)$, which implies
	\begin{equation*}\label{2410071}
		|j|\leq \min \left\{n_1^*, N_1^*\right\}.
	\end{equation*}
Moreover, a straightforward comparison of the involved modes yields
	\begin{equation}\label{2410072}
		\nu_1^* \leq \max \left\{n_1^*, N_1^*\right\}.
	\end{equation}
	Combining the above estimates, one has
	\begin{equation}\label{2410073}
		\ln ^\sigma\lfloor j \rfloor+\ln ^\sigma\left\lfloor \nu_1^*\right\rfloor-\ln ^\sigma\left\lfloor n_1^*\right\rfloor-\ln ^\sigma\left\lfloor N_1^*\right\rfloor \leq 0 .
	\end{equation}
	
	$\bullet$ Noting that the rearranged weighted sums satisfy
	\begin{align}\label{2410074}
		\sum_{i \geq 1} \ln ^\sigma\left\lfloor n_i^*\right\rfloor  
		& =	\sum_{n \in \mathbb{Z}}\left(2 a_n+k_n+k_n^{\prime}\right) \ln ^\sigma\lfloor n\rfloor
		\nonumber \\
		 &\geq \sum_{n \in  \mathbb{Z}}
		\left(2 a_n+k_n+k_n^{\prime}\right) 
		\geq \sum_{n \in \mathbb{Z}}
		\left(k_n+k_n^{\prime}\right).
	\end{align}
	Similarly, for the higher-order terms ($i \ge 3$), a application of Remark \ref{Re2.2} yields
	\begin{align}\label{2410075}
		\sum_{i \geq 3} \ln ^\sigma\left\lfloor n_i^*\right\rfloor & \geq \left(\sum_{n \in \mathbb{Z}} 2 a_n+k_n+k_n^{\prime} \right)-2 \nonumber\\
		& \geq \frac{1}{2} \left( \sum_{n \in \mathbb{Z}} 
		2 a_n+k_n+k_n^{\prime} \right) 
		 \geq \frac{1}{2} \left( \sum_{n \in \mathbb{Z}} k_n+k_n^{\prime} \right) .
	\end{align}
	Combining the above estimates, we arrive at
	\begin{equation}\label{2410076}
		\sum_{n \in \mathbb{Z}}
		\left(k_n+k_n^{\prime}\right)\left(K_n+K_n^{\prime}\right) \leq 2\left(\sum_{i \geq 1} \ln ^\sigma\left\lfloor n_i^*\right\rfloor\right)\left(\sum_{i \geq 3} \ln ^\sigma\left\lfloor N_i^*\right\rfloor\right)
	\end{equation}
	and
		\begin{equation}\label{2410076'}
		\sum_{n \in \mathbb{Z}}
		\left(k_n+k_n^{\prime}\right)\left(K_n+K_n^{\prime}\right) \leq 4\left(\sum_{i \geq 3} \ln ^\sigma\left\lfloor n_i^*\right\rfloor\right)\left(\sum_{i \geq 3} \ln ^\sigma\left\lfloor N_i^*\right\rfloor\right).
	\end{equation}


To establish \eqref{2410081}, we analyze the following two cases based on the relationship between $\nu_1^*$ and $N_1^*$:
	
	{\bf Case 1.} $\nu_1^* \leq N_1^*$.
	
	{\bf Case 1.1.} $|j| \leq n_3^* \Rightarrow
	\lfloor j \rfloor \leq\left\lfloor n_3^*\right\rfloor$.
	
	Using $0<\delta_1 \leq \rho / 4$, one has
	\begin{equation}\label{2410083}
		e^{2 \rho\left(\ln^\sigma\lfloor j \rfloor-\ln ^\sigma\left\lfloor n_1^*\right\rfloor\right)} \leq e^{2 \rho\left(\ln ^\sigma\left\lfloor n_3^*\right\rfloor-\ln ^\sigma\left\lfloor n_1^*\right\rfloor\right)} \leq e^{\frac{1}{2} \delta_1\left(\ln ^\sigma\left\lfloor n_3^*\right\rfloor-\ln ^\sigma\left\lfloor n_1^*\right\rfloor\right)}
	\end{equation}
	and then
	\begin{equation}\label{2410084}
		e^{2 \rho\left(\ln ^\sigma\lfloor j \rfloor+\ln ^\sigma\left\lfloor \nu_1^*\right\rfloor-\ln ^\sigma\left\lfloor n_1^*\right\rfloor-\ln ^\sigma\left\lfloor N_1^*\right\rfloor\right)} e^{-\frac{1}{2} \delta_1\left(\sum_{i \geq 3} \ln ^\sigma\left\lfloor n_i^*\right\rfloor\right)} \leq e^{-\frac{1}{6} \delta_1\left(\sum_{i \geq 1} \ln ^\sigma\left\lfloor n_i^*\right\rfloor\right)},
	\end{equation}
	which follows from (\ref{2410083}).
	
	Note that once $j, a, k, k'$ are fixed, the values of $A, K, K'$ are uniquely determined. Consequently, in view of  (\ref{2410076}) and (\ref{2410084}), we obtain
	\begin{align*}
		I \leq 2 \sum_{a, k, k^{\prime} \in \mathbb{N}^{\mathbb{Z}} } \left(\sum_{i \geq 1} \ln ^\sigma\left\lfloor n_i^*\right\rfloor\right)\left(\sum_{i \geq 3} \ln ^\sigma\left\lfloor N_i^*\right\rfloor\right) e^{-\frac{1}{6} \delta_1\left(\sum_{i \geq 1} \ln ^\sigma\left\lfloor n_i^*\right\rfloor\right)} e^{-\frac{1}{2} \delta_2\left(\sum_{i \geq 3} \ln ^\sigma\left\lfloor N_i^*\right\rfloor\right)} .
			\end{align*}

\begin{claim}[\cite{Cong2023TheEO}, Lemma 4.5]\label{Lem6.3}
	For $p \geq 1$ and $\delta \in (0, 1)$, let $g_{p,\delta}(x) := x^p e^{-\delta x}$. Then,
	\begin{equation}\label{Lem6.3Con}
		\max _{x \geq 1} g_{p, \delta}(x) \leq\left(\frac{p}{e \delta}\right)^p.
	\end{equation}
\end{claim}
Applying Claim \ref{Lem6.3} with $p = 1$ and $\delta \in \{\delta_1/12, \delta_2/2\}$, we have
		\begin{align*}
		I	\leq & \frac{48}{e^2 \delta_1 \delta_2} \sum_{a, k, k^{\prime} \in \mathbb{N}^{\mathbb{Z}}} e^{-\frac{1}{12} \delta_1\left(\sum_{i \geq 1} \ln ^\sigma\left\lfloor n_i^*\right\rfloor\right)}  \nonumber\\
		\leq & \frac{48}{e^2 \delta_1 \delta_2}\left(\sum_{a \in \mathbb{N}^{\mathbb{Z}}} e^{-\frac{1}{12} \delta_1 \sum_{n \in \mathbb{Z}} 2 a_n \ln^\sigma\lfloor n\rfloor}\right)\left(\sum_{k \in \mathbb{N}^{\mathbb{Z}}} e^{-\frac{1}{12} \delta_1 \sum_{n \in \mathbb{Z}} k_n \ln ^\sigma\lfloor n\rfloor}\right)^2 \\
		\leq& \frac{1}{\delta_2} \exp \left\{\frac{550}{\delta_1} \cdot \exp \left\{\left(\frac{48}{\delta_1}\right)^{\frac{1}{\sigma-1}}\right\}\right\},
	\end{align*}
	where the last inequality follows from Claim \ref{a5}.
	
	{\bf Case 1.2.} $|j|>n_3^* \Rightarrow |j|\in\left\{n_1^*, n_2^*\right\}$.
	
	
	From (\ref{2410073}) and \eqref{2410076'}, it follows that
	\begin{align}\label{260411-1}
		I &\leq 4 \sum_{a, k, k^{\prime} \in \mathbb{N}^{\mathbb{Z}}}\left(\sum_{i \geq 3} \ln ^\sigma\left\lfloor n_i^*\right\rfloor\right)\left(\sum_{i \geq 3} \ln ^\sigma\left\lfloor N_i^*\right\rfloor\right)  e^{-\frac{1}{2}\left(\delta_1\sum_{i \geq 3} \ln ^\sigma\left\lfloor n_i^*\right\rfloor
		 	+\delta_2 \sum_{i \geq 3} \ln ^\sigma\left\lfloor N_i^*\right\rfloor \right)}\nonumber\\
		&\le \frac{64}{e^2 \delta_1 \delta_2}  \sum_{a, k, k^{\prime} \in \mathbb{N}^{\mathbb{Z}}}
	    e^{-\frac{1}{4}\left(\delta_1\sum_{i \geq 3} \ln ^\sigma\left\lfloor n_i^*\right\rfloor
			+\delta_2 \sum_{i \geq 3} \ln ^\sigma\left\lfloor N_i^*\right\rfloor \right)},
	\end{align}
	where the last inequality is obtained by applying Claim \ref{Lem6.3} with $p=1$ and $\delta \in \{\delta_1/4, \delta_2/4\}$.
In light of (\ref{2410077}) and  \eqref{2410075}, it follows that
	\begin{equation}\label{2410086}
		\sum_{n \in \mathbb{Z}} 2 \alpha_n+\kappa_n+\kappa_n^{\prime} 
	 \leq 2\left(\sum_{i \geq 3} \ln ^\sigma\left\lfloor n_i^*\right\rfloor+ \ln ^\sigma\left\lfloor N_i^*\right\rfloor\right) .
	\end{equation}
	By choosing $\delta := \min\{\delta_1, \delta_2\}$, we obtain
	\begin{align*}
      &\quad\; e^{-\frac{1}{4}\left(\delta_1\sum_{i \geq 3} \ln ^\sigma\left\lfloor n_i^*\right\rfloor
			+\delta_2 \sum_{i \geq 3} \ln ^\sigma\left\lfloor N_i^*\right\rfloor \right)}\\
		&\le  e^{-\frac{1}{8} \delta_1\sum_{i \geq 3} \ln ^\sigma\left\lfloor n_i^*\right\rfloor  }
			e^{-\frac{1}{8}\delta\left(\sum_{i \geq 3} \ln ^\sigma\left\lfloor n_i^*\right\rfloor
		+ \ln ^\sigma\left\lfloor N_i^*\right\rfloor \right)}\\
		&\le  e^{-\frac{1}{8} \delta_1\sum_{i \geq 3} \ln ^\sigma\left\lfloor n_i^*\right\rfloor  }
		e^{-\frac{1}{16}\delta\left(\sum_{n \in \mathbb{Z}}2 \alpha_n+\kappa_n+\kappa_n^{\prime}\right)},
	\end{align*}
	where the last inequality follows from (\ref{2410086}). Consequently, combined with \eqref{260411-1}, this leads to
	\begin{align}\label{2410088}
		I &\leq   \frac{64}{e^2 \delta_1 \delta_2}  \sum_{a, k, k^{\prime} \in \mathbb{N}^{\mathbb{Z}}}
	e^{-\frac{1}{8} \delta_1\sum_{i \geq 3} \ln ^\sigma\left\lfloor n_i^*\right\rfloor  }
	e^{-\frac{1}{16}\delta\left(\sum_{n \in \mathbb{Z}}2 \alpha_n+\kappa_n+\kappa_n^{\prime}\right)}  .
	\end{align}

	Now we claim some facts:
	
	$\bullet$ Firstly, note that $\left\{n_1^*, n_2^*\right\} \cap \operatorname{supp} \mathcal{M}_{\alpha \kappa \kappa^{\prime}} \neq \emptyset$, thus $n_1^*$ or $n_2^*$ ranges in a set of cardinality no more than
	\begin{equation}\label{2410089}
		\# \operatorname{supp} \mathcal{M}_{\alpha \kappa \kappa^{\prime}} \leq \sum_{n \in \mathbb{Z}}\left(2 \alpha_n+\kappa_n+\kappa_n^{\prime}\right).
	\end{equation}
	
	$\bullet$ Secondly, if $\left\{n_i^*\right\}_{i \geq 3}$ and $n_1^*$ (resp. $n_2^*$) are specified, then $n_2^*$ (resp. $n_1^*$) is determined uniquely since the momentum conservation. 
	
	$\bullet$ Thirdly, if $\left\{n_i^*\right\}_{i \geq 1}$ are given, then $\left\{2 a_n+k_n+k_n^{\prime}\right\}_{n \in \mathbb{Z}}$ are specified, and hence $\left(a, k, k^{\prime}\right)$ is specified up to a factor of
	\begin{equation}\label{24100810}
		\prod_{n \in \mathbb{Z}}\left(1+\binom{l_n+2}{2}\right)\leq\prod_{n \in \mathbb{Z}}\left(1+l_n^2\right),
		\quad
			l_n=\#\left\{i: n_i^*=n\right\}.
	\end{equation}
	Since $n_1^* \geq n_2^*>n_3^*$, one has
	\begin{equation}\label{24100811}
		\prod_{n \in \mathbb{Z} , n=n_1^*, n_2^*}\left(1+l_n^2\right) \leq 5.
	\end{equation}
	
	Following (\ref{2410088})-(\ref{24100811}), we obtain
	\begin{align*}\label{24100812}
		I	&\le \frac{64}{e^2 \delta_1\delta_2} \sum_{\left\{n_i^*\right\}_{i \geq 1}} 
		\prod_{\substack{m \in \mathbb{Z} \\
		m \leq n_1^*}}\left(1+l_m^2\right) e^{-\frac{1}{8} \delta_1\left(\sum_{i \geq 3} \ln ^\sigma\left\lfloor n_i^*\right\rfloor\right)}  e^{-\frac{1}{16} \delta \sum_{n \in \mathbb{Z}}\left(2 \alpha_n+\kappa_n+\kappa_n^{\prime}\right)} \nonumber\\
		&\leq  \frac{320}{e^2 \delta_1\delta_2} \sum_{\left\{n_i^*\right\}_{i \geq 3}} \prod_{\substack{m \in \mathbb{Z} \\
		m \leq n_3^*}}\left(1+l_m^2\right) e^{-\frac{1}{8} \delta_1\left(\sum_{i \geq 3} \ln ^\sigma\left\lfloor n_i^*\right\rfloor\right)} \\
	    & \quad \times \left(\sum_{n \in \mathbb{Z}} 2 \alpha_n+\kappa_n+\kappa_n^{\prime} \right) e^{-\frac{1}{16} \delta \sum_{n \in \mathbb{Z}}\left(2 \alpha_n+\kappa_n+\kappa_n^{\prime}\right)} \nonumber\\
		&\leq \frac{6000}{e^3\delta_1 \delta_2 \delta} \sum_{\left\{n_i^*\right\}_{i \geq 3}} \prod_{\substack{m \in \mathbb{Z} \\
		m \leq n_3^*}}\left(1+l_m^2\right) e^{-\frac{1}{8} \delta_1\left(\sum_{i \geq 3} \ln ^\sigma\left\lfloor n_i^*\right\rfloor\right)} \quad
		(\mbox{by Claim \ref{Lem6.3}})\nonumber\\
		&\leq\frac{1}{\delta_2}\exp\left\{\frac{1000}{\delta_1}\cdot\exp\left\{\left(\frac{100}{\delta_1}\right)^{\frac{1}{\sigma-1}}\right\}\right\},
	\end{align*}
where the last inequality follows from  Claim \ref{a5} and the following claim:
\begin{claim}[\cite{Cong2023TheEO}, Lemma 4.8]\label{Lem6.7}
	Assuming $\sigma>2, \delta \in(0,1), p=1,2$ and $a=\left(a_n\right)_{n \in \mathbb{Z}} \in \mathbb{N}^{\mathbb{Z}}$, then we have
	\begin{equation}\label{24100813}
		\prod_{n \in \mathbb{Z}}\left(1+a_n^p\right) e^{-2 \delta a_n\ln ^\sigma\lfloor n\rfloor} \leq \exp \left\{3 p\left(\frac{p}{\delta}\right)^{\frac{1}{\sigma-1}} \cdot \exp \left\{\left(\frac{1}{\delta}\right)^{\frac{1}{\sigma}}\right\}\right\}.
	\end{equation}	
\end{claim}

{\bf Case 2.} $\nu_1^*>N_1^*$.

In view of (\ref{2410072}), one has $\nu_1^*=n_1^*$. Hence, $n_2^*$ is determined by $n_1^*$ and $\left\{n_i^*\right\}_{i \geq 3}$. Similar as Case 1.2, we have
$$
I  \leq\frac{1}{\delta_2}\exp\left\{\frac{1000}{\delta_1}\cdot\exp\left\{\left(\frac{100}{\delta_1}\right)^{\frac{1}{\sigma-1}}\right\}\right\}.
$$    

Combining Case 1 and Case 2, we conclude the proof of the estimate \eqref{2410081}.
\end{proof}

	\begin{proof}[Proof of Lemma \ref{E1}]
	We begin by establishing the estimate \eqref{260413-1}. To this end, let us expand the composition $R \circ \Phi$ into a Taylor series
\begin{equation}\label{2410091}
	R \circ \Phi=\sum_{n \geq 0} \frac{1}{n!} R^{(n)},
\end{equation}
where the sequence of Hamiltonians is defined recursively by $$R^{(0)} = R,\quad R^{(n)} = \{R^{(n-1)}, F\}.$$
We proceed to estimate the norm $\|R^{(n)}\|_\rho^+$ by applying the estimate \eqref{042704} from Lemma \ref{010}  iteratively, with the choices $\delta_1 = \frac{\delta}{2}$ and $\delta_2 = \frac{\delta}{2n}$:
\begin{align}\label{2410092}
	\|  R^{(n)} \|^{c,+}_{\rho+\delta} 
	&\leq  
	C\left( \frac{\delta}{2},\frac{\delta}{2n},\sigma \right) 	\|F\|^{c,+}_{\rho+\frac{\delta}{2}}
	\left\|R^{(n-1)}\right\|^{c,+}_{\rho+\delta-\frac{\delta}{2 n}}\nonumber \\
	& \leq\left( 
	C\left( \frac{\delta}{2},\frac{\delta}{2n},\sigma \right)  
	\|F\|^{c,+}_{\rho+\frac{\delta}{2}-\frac{\delta}{2 n}}\right)^2
   \left\|R^{(n-2)}\right\|^{c,+}_{\rho+\delta-\frac{2 \delta}{2 n}} \nonumber\\
	&\quad \cdots \cdots  \nonumber\\
	&\leq\left(
 	C\left(\frac{\delta}{2},\frac{\delta}{2n},\sigma \right) 
	\|F\|^{c,+}_{\rho} \right)^n
	\|R\|^{c,+}_{\rho+\frac{\delta}{2}} \nonumber\\
	&\leq\left(
	C\left(\frac{\delta}{2},\frac{\delta}{2n},\sigma \right) 
	\|F\|^{c,+}_{\rho} \right)^n
	\|R\|^{c,+}_{\rho},
\end{align}
where the constant is given by
\begin{equation*}
	C\left( \frac{\delta}{2},\frac{\delta}{2n},\sigma \right) = \frac{2 n}{\delta} \exp \left\{\frac{2000}{\delta} \cdot \exp \left\{\left(\frac{200}{\delta}\right)^{\frac{1}{\sigma-1}}\right\}\right\}=:\frac{2 n}{\delta} C_1\left({\delta},\sigma \right).
\end{equation*}
Next, we observe that
\begin{equation*}
	\frac{1}{n!}	  C^n\left( \frac{\delta}{2},\frac{\delta}{2n},\sigma \right) 
	 = \frac{1}{n!}\left(\frac{2 n}{\delta}\right)^n 
	C^n_1\left({\delta},\sigma \right) 
	\le \left( \frac{2e }{\delta}
    	C_1\left({\delta},\sigma \right)\right)^n ,
\end{equation*}
where the last inequality is a consequence of the bound $n^n < n! e^n$.
Consequently, combining (\ref{2410091}) and (\ref{2410092}), it follows that for any $ i\ge 0$,
\begin{align*}
	\left\|  \sum_{n \geq i}\frac{1}{n!} R^{(n)}   \right\| ^{c,+}_{\rho+\delta}
	&\leq \sum_{n \geq i} \frac{1}{n!}\left(
	C\left(\frac{\delta}{2},\frac{\delta}{2n},\sigma \right) 
	\|F\|^{c,+}_{\rho } \right)^n
	\|R\|^{c,+}_{\rho } \nonumber \\
	&\leq \sum_{n \geq i}\left(\frac{2 e}{\delta} C_1\left({\delta},\sigma \right)
	\|F\|^{c,+}_{\rho}\right)^n\|R\|^{c,+}_{\rho}\nonumber \\
	&\leq 2\left( \frac{2 e}{\delta} C_1\left({\delta},\sigma \right) \|F\|^{c,+}_{\rho} \right)^i 
	\|R\|^{c,+}_{\rho}.
\end{align*}
This completes the proof of \eqref{260413-1}. 

Following an argument similar to that for \eqref{260413-1} and incorporating the bound \eqref{250804-3} in Lemma \ref{010},
we obtain
\begin{equation*}
	 	\left\| \sum_{n \geq i}\frac{1}{n!} R^{(n)} \right\|^{*,+}_{\rho+\delta}
	 \leq C(\delta, \sigma)  \left( \|F\|_{\rho}^{c,+}\right)^i \|R\|_{\rho}^{*,+}
\end{equation*}
and
\begin{equation*}
	\left\|  \sum_{n \geq i}\frac{1}{n!} R^{(n)}  \right\|^{*,+}_{\rho+\delta}
	\leq C(\delta, \sigma) \left(  \|F\|_{\rho}^{*,+}\right) ^i
	 \|R\|_{\rho}^{c,+},
\end{equation*}
respectively. This completes the proof of the estimate \eqref{260413-4}.
		\end{proof}

	\subsection{Proof of Lemma \ref{051301} }\label{sec6.2}
		\begin{proof}
According to the definition of the weighted $\ell^\infty$ norm (Definition \ref{def_norm} and \ref{083103}), without loss of generality, we shall only provide the proofs for (\ref{N6}) and (\ref{N7}).

			Before the proof begins, we rewrite the Hamiltonian $R$ in (\ref{pertur}) as
			\begin{align*}
				\nonumber R(z,\bar z)
				=R^0+R^1+R^2,
			\end{align*}
			where
			\begin{align}
				&\label{1219001} R^0=\sum_{\substack{a, k, k^{\prime} \in \mathbb{N}^{\mathbb{Z}} \\ \operatorname{supp} k \cap \operatorname{supp} k^{\prime}=\emptyset}} R_{a k k^{\prime}} \mathcal{M}_{a k k^{\prime}},\\
				&\label{1219002} R^1=\sum_{m \in \mathbb{Z}} \left(\sum_{\substack{a, k, k^{\prime} \in \mathbb{N}^{\mathbb{Z}} \\
						\operatorname{supp} k \cap \operatorname{supp} {k^{\prime}}=\emptyset}} 
					R_{a k k^{\prime}}^{(m)} J_m\mathcal{M}_{a k k^{\prime}}\right), \\
				&\label{1219003} R^2=\sum_{m_1, m_2 \in \mathbb{Z}} \left(\sum_{\substack{
				a, k, k^{\prime} \in \mathbb{N}^{\mathbb{Z}}\\
					\text { no assumption }
				}} R_{a k k^{\prime}}^{\left(m_1, m_2\right)}J_{m_1} J_{m_2} \mathcal{M}_{a k k^{\prime}}\right),
			\end{align}
			where the so-called monomials
			$$
			\mathcal{M}_{a k k^{\prime}}=\prod_{n \in \mathbb{Z}} I_n(0)^{a_n} z_n^{k_n} \bar{z}_n^{k_n^{\prime}}=:\prod_{n \in \mathbb{Z}} I_n(0)^{a_n} z_n^{k_n} \bar{z}_n^{k_n^{\prime}}
			$$
			satisfy the zero momentum conservation \eqref{mome}.
			Here, we denote $I_n(0) \equiv I_n(0)$ for notational convenience, and $R_{akk'}$, $R_{akk'}^{(m)}$, and $R_{akk'}^{(m_1,m_2)}$ represent the corresponding coefficients. 
			 Moreover, we also rewrite the Hamiltonian $R$ in (\ref{pertur1}) as
			\begin{equation}\label{121903}
				R(z,\bar z)=\sum_{\alpha ,\kappa,\kappa'\in\mathbb{N}^{\mathbb{Z}}}
				R_{\alpha \kappa\kappa'} \mathcal{M}_{\alpha \kappa \kappa^{\prime}}.
			\end{equation}

			Now we prove the inequality (\ref{N6}).

			First, consider the term $\mathcal{M}_{a k k^{\prime}}$ in \eqref{1219001}
			for fixed $a, k, k' \in \mathbb{N}^{\mathbb{Z}}$ satisfying $k_n \cdot k'_n = 0$ for all $n \in \mathbb{Z}$. It is evident that $\mathcal{M}_{akk'}$ originates from specific components of the terms $\mathcal{M}_{\alpha\kappa\kappa'}$ without restrictions on $\kappa$ and $\kappa'$. We express $\mathcal{M}_{\alpha\kappa\kappa'}$ as
			$$
			\mathcal{M}_{\alpha \kappa \kappa^{\prime}}=\mathcal{M}_{\alpha b l l^{\prime}}=\prod_{n \in \mathbb{Z}} I_n(0)^{\alpha_n} I_n^{b_n} z_n^{l_n} \bar{z}_n^{l_n^{\prime}},
			$$	
			where	
			$$
			b_n=\min \left\{\kappa_n, \kappa_n^{\prime}\right\}, \quad l_n=\kappa_n-b_n, \quad l_n^{\prime}=\kappa_n^{\prime}-b_n
			$$	
			and then $l_n \cdot l_n^{\prime}=0$ for all $n \in \mathbb{Z}$.
		By substituting the expansion $I_n = I_n(0) + J_n$ into the product $\prod_{n \in \mathbb{Z}} I_n^{b_n}$, we can express the term as a sum of monomials of the following forms:
			$$	\prod_{n \in \mathbb{Z}} I_n(0)^{b_n},$$
			which contributes to  $R^0$;
			for $m\in\mathbb{Z}$,
			$$		\left(b_m I_m(0)^{b_m-1} J_m\right)\left(\prod_{n \neq m} I_n(0)^{b_n}\right),$$
			which contributes to  $R^1$;
			and	for $m,m_1,m_2\in\mathbb{Z}$,
			$$	\left(\prod_{n<m} I_n(0)^{b_n}\right)\left(
			C_{r+2}^2 I_m(0)^r J_m^2 I_m^{b_m-r-2}\right)\left(\prod_{n>m} I_n^{b_n}\right),\quad r\le b_m-2, $$	
			\begin{align*}
				&\left(\prod_{n<m_1} I_n(0)^{b_n}\right)
				\left(b_{m_1}I_{m_1}(0)^{b_{m_1}-1} J_{m_1}\right)
				\left(\prod_{m_1<n<m_2} I_n(0)^{b_n}\right)\\
				&\quad\times \left(C_{r+1}^1 I_{m_2}(0)^r J_{m_2} I_{m_2}^{b_{m_2}-r-1}\right)\left(\prod_{n>m_2} I_n^{b_n}\right),\quad r\le b_{m_2}-1.
			\end{align*}		
			which contribute to  $R^2$.
	Next, we estimate the bounds of the coefficients for the terms mentioned above.

	First, we consider the term $\mathcal{M}_{akk'}$ in \eqref{1219001}. 
	For any given $n$, we have
			\begin{equation}\label{Makk}
				I_n(0)^{a_n} z_n^{k_n} \bar{z}_n^{k_n^{\prime}}=\sum_{b_n=\min \left\{\kappa_n, \kappa_n^{\prime}\right\}} I_n(0)^{\alpha_n+b_n} z_n^{\kappa_n-b_n} \bar{z}_n^{\kappa_n^{\prime}-b_n},
			\end{equation}
			where
			\begin{equation}\label{a=02}
				\alpha_n+b_n=a_n,\quad\kappa_n-b_n=k_n, \quad \kappa_n^{\prime}-b_n=k_n^{\prime}.
			\end{equation}
			Consequently, it follows that
			\begin{equation}\label{a=03}
				2a_n+k_n+k_n^{\prime}
				=2 \alpha_n+\kappa_n+\kappa_n^{\prime}
					\quad\text{and} \quad
				n_1^*(a, k, k^{\prime}) =  n_1^*(\alpha,\kappa, \kappa^{\prime}).
			\end{equation}

			For fixed $a, k, k'$, once $0 \le \alpha_n \le a_n$ is specified, the indices $b_n, \kappa_n, \kappa'_n$ are uniquely determined by (\ref{a=02}). Since $\alpha_n$ has at most $(1 + a_n)$ choices for each $n$, we obtain
				\begin{align}\label{Conb0}
				\nonumber&\left|\left(\prod_{n \in \mathbb{Z}} \left\langle \frac{n}{c}\right\rangle
				^{\frac{1}{2}\left(2 a_n+k_n+k_n^{\prime}\right)}\right) R_{a k k^{\prime}}\right|\\
				\leq& \left|\left(\prod_{n \in \mathbb{Z}} \left\langle \frac{n}{c}\right\rangle^{\frac{1}{2}\left(2 \alpha_n+\kappa_n+\kappa_n^{\prime}\right)}\right) R_{\alpha \kappa \kappa^{\prime}}\right|
				\left(\prod_{n \in \mathbb{Z}} 1+a_n\right)
				\nonumber\\
				\leq&\|R\|_\rho^{c,+} e^{\rho\left(\sum_{n \in \mathbb{Z}}\left(2 \alpha_n+\kappa_n+\kappa_n^{\prime}\right) \ln ^\sigma\lfloor n\rfloor-2 \ln ^\sigma\left\lfloor n_1^*(\alpha,\kappa, \kappa^{\prime})\right\rfloor\right)}
				\left(\prod_{n \in \mathbb{Z}} 1+a_n \right),
			\end{align}
			where the last inequality follows from (\ref{norm_1}).
			Note that from \eqref{a=03}, one has
			\begin{align*}
				&\quad 	\sum_{n \in \mathbb{Z}}\left(2 \alpha_n+\kappa_n+\kappa_n^{\prime}\right) \ln ^\sigma\lfloor n\rfloor-2 \ln ^\sigma\left\lfloor n_1^*(\alpha,\kappa, \kappa^{\prime}) \right\rfloor \\
				&= \sum_{n \in \mathbb{Z}}\left(2 a_n+k_n+k_n^{\prime}\right) \ln ^\sigma\lfloor n\rfloor-2 \ln ^\sigma\left\lfloor n_1^*(a,k,k')\right\rfloor.
			\end{align*}
		Combining this identity with (\ref{Conb0}), we deduce
			\begin{align}\label{250519}
				\left\|R^0\right\|^c_{\rho+\delta} & \leq\|R\|_\rho^{c,+}
				\sup_{a,k,k'\in\mathbb{N}^{\mathbb{Z} }} e^{-\delta\left(\sum_{n \in \mathbb{Z}}\left(2 a_n+k_n+k_n^{\prime}\right) \ln ^\sigma\lfloor n\rfloor-2 \ln ^\sigma\left\lfloor n_1^*\right\rfloor\right)}
				\left(\prod_{n \in \mathbb{Z}} 1+a_n \right)
				\nonumber \\
				&  \leq\|R\|_\rho^{c,+} \sup_{a,k,k'\in\mathbb{N}^{\mathbb{Z} }}   e^{-\frac{1}{2} \delta \sum_{i \geq 3} \ln ^\sigma\left\lfloor n_i^*\right\rfloor} 
				\left(\prod_{n \in \mathbb{Z}} 1+a_n \right),
			\end{align}
			where the last inequality follows from Lemma \ref{005}.
			Consequently, the problem of obtaining the desired estimate is reduced to establishing the following upper bound for any $a, k, k' \in \mathbb{N}^{\mathbb{Z}}$:
			\begin{equation}\label{260414-2}
					\mathcal{Q} :=
				\left( \prod_{n \in \mathbb{Z}} 1 + a_n \right)  e^{-\frac{1}{2}\delta \sum_{i \ge 3} \ln^\sigma |n_i^*|} \le \exp \left\{3\left(\frac{6}{\delta}\right)^{\frac{1}{\sigma-1}} \cdot \exp \left\{\left(\frac{6}{\delta}\right)^{\frac{1}{\sigma}}\right\}\right\}. 
			\end{equation}
 To this end, we consider the following three cases:
			
			{\bf Case 1.} $n_1^*=n_2^*=n_3^*$.
			
			One has
			$$
			\begin{aligned}
				\mathcal{Q}  & \leq \left( \prod_{n \in \mathbb{Z}} 1 + a_n \right)  e^{-\frac{1}{6}\delta \sum_{i \ge 1} \ln^\sigma |n_i^*|}  \\
				& \leq 
				\left(\prod_{n \in \mathbb{Z}}  1+a_n\right) 
				e^{-\frac{1}{3} \delta a_n \ln ^\sigma\lfloor n\rfloor}   \\
				& \leq \exp \left\{3\left(\frac{6}{\delta}\right)^{\frac{1}{\sigma-1}} \cdot \exp \left\{\left(\frac{6}{\delta}\right)^{\frac{1}{\sigma}}\right\}\right\},
			\end{aligned}
			$$
			where the last inequality is based on Claim \ref{Lem6.7} (with $p=1$ and $\delta$ replaced by $\delta/6$).
			
			{\bf Case 2.} $n_1^*>n_2^*=n_3^*$.
			
			In this case, one has $a_{n_1^*}=0$, otherwise $n_1^*=n_2^*$. Then we have
			
			$$
			\begin{aligned}
					\mathcal{Q} & =
					\left( \prod_{ |n|\le n_2^* } 
					1 + a_n \right) 
					 e^{-\frac{1}{2}\delta \sum_{i \ge 3} \ln^\sigma |n_i^*|} \\
				& \leq 	\left( \prod_{ |n|\le n_2^* } 
				1 + a_n \right) 
				e^{-\frac{1}{4}\delta \sum_{i \ge 2} \ln^\sigma |n_i^*|}  \\
				& \leq \exp \left\{3\left(\frac{4}{\delta}\right)^{\frac{1}{\sigma-1}} \cdot \exp \left\{\left(\frac{4}{\delta}\right)^{\frac{1}{\sigma}}\right\}\right\},
			\end{aligned}
			$$
			where the last inequality is based on Claim \ref{Lem6.7} (with $p=1$ and $\delta$ replaced by $\delta/4$).
			
			{\bf Case 3.} $n_2^*>n_3^*$.
			
			In this case, one has 
			$$
			\prod_{|n|=n_1^*, n_2^*}1+a_n  \leq 2.
			$$
			Hence we have
			$$
			\begin{aligned}
					\mathcal{Q} & \leq 2\left( \prod_{ |n|\le n_3^* } 
					1 + a_n \right) 
					e^{-\frac{1}{2}\delta \sum_{i \ge 3} \ln^\sigma |n_i^*|} \\
				& \leq 2 \exp \left\{3\left(\frac{2}{\delta}\right)^{\frac{1}{\sigma-1}} \cdot \exp \left\{\left(\frac{2}{\delta}\right)^{\frac{1}{\sigma}}\right\}\right\},
			\end{aligned}
			$$
			where the last inequality is based on Claim \ref{Lem6.7}  (with $p=1$ and $\delta$ replaced by $\delta/2$).
			
		Combining the estimates from the three cases above with \eqref{250519}, we obtain
			\begin{equation}\label{2410103}
				\left\|R^0\right\|^c_{\rho+\delta} \leq\|R\|_\rho^{c,+} \exp \left\{3\left(\frac{6}{\delta}\right)^{\frac{1}{\sigma-1}} \cdot \exp \left\{\left(\frac{6}{\delta}\right)^{\frac{1}{\sigma}}\right\}\right\}.
			\end{equation}
			
Next, we consider the term $J_m \mathcal{M}_{akk'}$ in \eqref{1219002}. For fixed $a, k, k'$ satisfying $k_n \cdot k'_n = 0$ for all $n \in \mathbb{Z}$, this term can be decomposed into two parts:
			$$
			J_m\mathcal{M}_{akk^{\prime}}=\left( J_m I_m(0)^{a_m} z_m^{k_m} \bar{z}_m^{k'_m}\right)\cdot
			 \left(  \prod_{n \neq m} I_n(0)^{a_n} z_n^{k_n} \bar{z}_n^{k'_n} \right) =: \textcircled{1} \cdot \textcircled{2}.
			$$
		The term $J_m \mathcal{M}_{a k k^{\prime}}$ also comes from some parts of the terms $\mathcal{M}_{\alpha \kappa \kappa^{\prime}}$ with no assumption for $\kappa$ and $\kappa^{\prime}$. 
			
			On the one hand, we analyze part $\textcircled{2}$. For any $n \neq m$, we have
			$$
			I_n(0)^{a_n} z_n^{k_n} \bar{z}_n^{k_n^{\prime}}=\sum_{b_n=\min\{\kappa_n, \kappa_n^{\prime}\}} I_n(0)^{\alpha_n+b_n} z_n^{\kappa_n-b_n} \bar{z}_n^{\kappa_n^{\prime}-b_n}.
			$$
			Following the same logic as in (\ref{a=02}) and (\ref{a=03}), the indices $b_n, \kappa_n, \kappa'_n$ are uniquely determined once $0 \le \alpha_n \le a_n$ is fixed.

			On the other hand, we consider part $\textcircled{1}$. For $n = m$, it follows that
			$$
			J_m I_m(0)^{a_m} z_m^{k_m} \bar{z}_m^{k_m^{\prime}}=\sum_{b_n=\min\{\kappa_n, \kappa_n^{\prime}\}} b_m J_m I_m(0)^{\alpha_m+b_m-1} z_m^{\kappa_m-b_m} \bar{z}_m^{\kappa_m^{\prime}-b_m}.
			$$
			Hence, we have
			\begin{equation}\label{24100961}
					\alpha_m+b_m-1=a_m, \quad\kappa_m-b_m=k_m, \quad \kappa_m^{\prime}-b_m=k_m^{\prime}
			\end{equation}
			and then
			\begin{equation}\label{2410097}
				2 \alpha_m+\kappa_m+\kappa_m^{\prime}=2 a_m+k_m+k_m^{\prime}+2
				\quad\text{and}\quad 
				n_1^*(a+e_m, k, k^{\prime}) = n_1^*(\alpha,\kappa, \kappa^{\prime}).
			\end{equation}	
		Moreover, if $0 \le \alpha_m \le a_m$ is chosen (noting that $\alpha_m = a_m + 1 \iff b_m = 0$), the indices $b_m, \kappa_m, \kappa'_m$ are likewise determined by  (\ref{24100961}).
			Then, one has 
				\begin{align}
					& \left|
					\left(\left\langle \frac{m}{c}\right\rangle 
					\prod_{n \in \mathbb{Z}} 	  \left\langle \frac{n}{c}\right\rangle^{\frac{1}{2}\left(2 a_n+k_n+k_n^{\prime}\right)}\right) 		
					  R_{a k k^{\prime}}^{(m)}\right|\nonumber \\
					\leq & \left(\prod_{n \in \mathbb{Z}} 	\left\langle \frac{n}{c}\right\rangle^{\frac{1}{2}\left(2 \alpha_n+\kappa_n+\kappa_n^{\prime}\right)}\right)  
					 \left(\prod_{n \in \mathbb{Z}} 1+a_n\right)  b_m   \left| R_{\alpha \kappa \kappa^{\prime}}\right|\nonumber \\
					\leq &  \|R\|_\rho^{c,+} 
				e^{\rho\left(\sum_{n \in \mathbb{Z}}\left(2 \alpha_n+\kappa_n+\kappa_n^{\prime}\right) \ln ^\sigma\lfloor n\rfloor-2 \ln ^\sigma\left\lfloor n_1^*(\alpha,\kappa, \kappa^{\prime})\right\rfloor\right)}
					\left(\prod_{n \in \mathbb{Z}} 1+a^2_n \right),\label{260414-1}
				\end{align}
		where the last inequality follows from   (\ref{norm_1}).
		Note that from \eqref{2410097}, one has
		\begin{align*}
				&\quad 	\sum_{n \in \mathbb{Z}}\left(2 \alpha_n+\kappa_n+\kappa_n^{\prime}\right) \ln ^\sigma\lfloor n\rfloor-2 \ln ^\sigma\left\lfloor n_1^*(\alpha,\kappa, \kappa^{\prime}) \right\rfloor \\
				&=
			\sum_{n \in \mathbb{Z}}\left(2 a_n+k_n+k_n^{\prime}\right) \ln ^\sigma\lfloor n\rfloor+2 \ln ^\sigma\left\lfloor m\right\rfloor-2 \ln ^\sigma\left\lfloor n_1^*\right\rfloor,
		\end{align*}
			where $n_1^*=n_1^*(a+e_m,k,k') $.
			Combining this identity with \eqref{260414-1},
			 we deduce
			\begin{align}\label{2410101}
				\left\|R^1\right\|^c_{\rho+\delta} & \leq\|R\|_\rho^{c,+}
				\sup_{a,k,k'\in\mathbb{N}^{\mathbb{Z} }} e^{-\delta\left(\sum_{n \in \mathbb{Z}}\left(2 a_n+k_n+k_n^{\prime}\right) \ln ^\sigma\lfloor n\rfloor -2 \ln^\sigma\left\lfloor n_1^*\right\rfloor\right)}
				\left(\prod_{n \in \mathbb{Z}} 1+a_n^2 \right)
				 \nonumber\\
				&  \leq\|R\|_\rho^{c,+} \sup_{a,k,k'\in\mathbb{N}^{\mathbb{Z} }} 
				\left(\prod_{n \in \mathbb{Z}}1+a_n^2 \right) e^{-\frac{1}{2} \delta \sum_{i \geq 3} \ln ^\sigma\left\lfloor n_i^*\right\rfloor}\quad(\text{by Lemma \ref{005}})\nonumber\\
				& \leq  \|R\|_\rho^{c,+} \exp \left\{6\left(\frac{12}{\delta}\right)^{\frac{1}{\sigma-1}} \cdot \exp \left\{\left(\frac{6}{\delta}\right)^{\frac{1}{\sigma}}\right\}\right\},
			\end{align}
			where the last inequality follows from the same argument as in \eqref{260414-2} , combined with Claim \ref{Lem6.7} for the case $p=2$.
			
			Finally, we consider the term $J_{m_1} J_{m_2} \mathcal{M}_{a k k^{\prime}}$ in \eqref{1219003}. 
			Following the proof of \eqref{2410101}, we can deduce
				\begin{align}\label{2410102}
				\left\|R^2\right\|^c_{\rho+\delta} 
				&  \leq\|R\|_\rho^{c,+} \sup_{a,k,k'\in\mathbb{N}^{\mathbb{Z} }} 
				\left(\prod_{n \in \mathbb{Z}}1+a_n^2\right) e^{-\frac{1}{2} \delta \sum_{i \geq 3} \ln ^\sigma\left\lfloor n_i^*\right\rfloor}\nonumber\\
				& \leq  \|R\|_\rho^{c,+} \exp \left\{6\left(\frac{12}{\delta}\right)^{\frac{1}{\sigma-1}} \cdot \exp \left\{\left(\frac{6}{\delta}\right)^{\frac{1}{\sigma}}\right\}\right\}.
			\end{align}
			
	In summary, by synthesizing the estimates for the three cases above—namely (\ref{2410103}), (\ref{2410101}) and (\ref{2410102}), we finish the proof of (\ref{N6}).

			Next, we shall prove inequality (\ref{N7}). 
			
			For any given $\alpha ,\kappa,\kappa'\in\mathbb{N}^{\mathbb{Z}}$, we consider the term $\mathcal{M}_{\alpha\kappa\kappa^{\prime}}$ in (\ref{121903}). In view of (\ref{1219001})-(\ref{1219003}), it is easy to see that $\mathcal{M}_{\alpha\kappa\kappa^{\prime}}$ comes from some parts of the terms
		\begin{equation*}
			J^b\mathcal{M}_{a k k^{\prime}}=\left\{\begin{array}{l}
			\prod_{n \in \mathbb{Z}} I_n(0)^{a_n} z_n^{k_n} \bar{z}_n^{k_n^{\prime}} ,\quad\mbox{for}\ |b|=0,\\[7pt]
				J_{m}I_m(0)^{a_m} z_m^{k_m} \bar z_m^{k_m^{\prime}}
			\prod_{	n\neq m}
			 I_n(0)^{a_n} z_n^{k_n} \bar z_n^{k_n^{\prime}},
				 \quad\mbox{for}\ |b|=1,\\[7pt]
				\prod_{m=m_1,m_2}J_{m}I_{m}(0)^{a_{m}} z_{m}^{k_{m}} \bar z_{m}^{k_{m}^{\prime}}
			\prod_{	n\neq m_1, m_2	} I_n(0)^{a_n} z_n^{k_n}\bar{z}_n^{k^{\prime}_n},\quad\mbox{for}\ |b|=2.
					\end{array}\right.
				\end{equation*}
				Now we discuss the relationship between the coefficients of $\mathcal{M}_{\alpha\kappa\kappa^{\prime}}$ and $J^b\mathcal{M} _{akk^{\prime}}$. Precisely,
				
				$\bullet$ If $|b| = 0$, we have $a_n = \alpha_n$, $k_n = \kappa_n$, and $k'_n = \kappa'_n$.
				
				$\bullet$ If $|b| = 1$ and $n \neq m$, one has $a_n = \alpha_n$, $k_n = \kappa_n$, and $k'_n = \kappa'_n$. If $n = m$, we obtain two possible forms
				\begin{equation}\label{2410221}
					I_m(0)^{a_m+1}  z_m^{k_m} \bar{z}_m^{k_m^{\prime}}\quad{\mbox{ and }}\quad
					I_m(0)^{a_m}  z_m^{k_m+1} \bar{z}_m^{k_m^{\prime}+1}.
				\end{equation}
				For the former in (\ref{2410221}), we have $a_m + 1 = \alpha_m$, $k_m = \kappa_m$, and $k'_m = \kappa'_m$; for the latter, we have $a_m = \alpha_m$, $k_m + 1 = \kappa_m$, and $k'_m + 1 = \kappa'_m$. In either case, we notice that $m$ is uniquely determined once $\alpha, \kappa$, and $\kappa'$ are given. Moreover, since $m \in \text{supp } \alpha$, there are at most $\sum_{n\in\mathbb{Z}} \alpha_n$ choices for $m$. Noting that $\alpha_n \leq a_n + b_n$, it follows that
				$$\sum_{n \in \mathbb{Z}}\alpha_n\leq\sum_{n \in \mathbb{Z}} a_n+b_n.$$
				
				$\bullet$ If $|b| = 2$, using the same argument as in the case $|b| = 1$, we find that $m_1$ and $m_2$ have at most 
				\begin{equation}\label{1219004}
					\left(\sum_{n \in \mathbb{Z}}  a_n+b_n  \right)^2
				\end{equation} 
				possible choices.

				Consequently, the coefficient of $J^b \mathcal{M}_{akk'}$ increases by at most the factor in (\ref{1219004}),
				%
				then one has
				$$
				\begin{aligned}
					\|R\|_{\rho+\delta}^{c,+} &\leq\|R\|^c_\rho\sup_{a,k,k'\in\mathbb{N}^{\mathbb{Z} }}
					\left(\sum_{n \in \mathbb{Z}} a_n+b_n \right)^2
					 e^{-\delta\left(\sum_{n \in \mathbb{Z}}\left(2 a_n+k_n+k_n^{\prime}\right) \ln ^\sigma\lfloor n\rfloor-2\ln^\sigma\left\lfloor n_1^*\right\rfloor\right)} \\
					& \leq\|R\|^c_\rho\sup_{a,k,k'\in\mathbb{N}^{\mathbb{Z} }}\left(2 \sum_{i \geq 3} \ln ^\sigma\left\lfloor n_i^*\right\rfloor\right)^2 e^{-\frac{1}{2} \delta \sum_{i \geq 3} \ln ^\sigma\left\lfloor n_i^*\right\rfloor} \quad(\mbox{by Lemma \ref{005}})\\
					& \leq \frac{64}{e^2 \delta^2}\|R\|^c_\rho,
				\end{aligned}
				$$
				where the last inequality  follows from Claim \ref{Lem6.3} with $p=2$.
			We complete the proof of (\ref{N7}).
			\end{proof}
		%
		%

			\section*{Acknowledgments}
			On behalf of all authors, the corresponding author states that there is no conflict of interest.

            \textbf{Data availability:} Data sharing not applicable to this article as no datasets were generated or analysed during the current study.
			
			
			\bibliographystyle{alpha}
			\bibliography{phase1}

		\end{document}